\documentstyle[11pt,twoside]{article}

\newcommand{\ncm}{\newcommand}
\ncm{\rncm}{\renewcommand}

\ncm{\lb}{\label}

\newtheorem{theorem}{Theorem }[section]
\newtheorem{lemma}[theorem]{Lemma }
\newtheorem{proposition}[theorem]{Proposition }
\newtheorem{corollary}[theorem]{Corollary }
\newtheorem{definition}[theorem]{Definition }

\ncm{\Theorem}[2]{\begin{theorem} \lb{Thm #1} {\sl #2} \end{theorem}}

\ncm{\thm}[1]{Theorem \ref{Thm #1}}

\ncm{\Definition}[2]{\begin{definition} \lb{Def #1}{\rm #2}\end{definition}}

\ncm{\defi}[1]{Definition \ref{Def #1}}

\ncm{\Lemma}[2]{\begin{lemma} \lb{Lem #1} {\sl #2} \end{lemma}}

\ncm{\lem}[1]{Lemma \ref{Lem #1}}

\ncm{\Proposition}[2]{\begin{proposition}\lb{Prop #1}{\sl #2}
		      \end{proposition}}

\ncm{\prop}[1]{Proposition \ref{Prop #1}}

\ncm{\Corollary}[2]{\begin{corollary}\lb{Cor #1} {\sl #2} \end{corollary}}

\ncm{\cor}[1]{Corollary \ref{Cor #1}}

\ncm{\Thm}{\Theorem}
\ncm{\Def}{\Definition}
\ncm{\Lem}{\Lemma}
\ncm{\Prop}{\Proposition}
\ncm{\Cor}{\Corollary}

\renewcommand{\theequation}{\mbox{\arabic{section}.\arabic{equation}}}

\ncm{\setc}[1]{\setcounter{equation}{#1}}
\rncm{\sec}{\setc{0}\section}
\ncm{\no}[1]{(\ref{#1})}
\ncm{\Eq}[1]{Eq.\,\no{#1}}

\newcommand{\be}{\begin{equation}}
\newcommand{\ee}{\end{equation}}
\ncm{\beq}{\begin{equation}}
\ncm{\eeq}{\end{equation}}

\ncm{\bleq}[1]{\beq\lb{#1}}

\ncm{\bea}{\begin{eqnarray}}
\ncm{\eea}{\end{eqnarray}}
\ncm{\beanon}{\begin{eqnarray*}}
\ncm{\eeanon}{\end{eqnarray*}}
\ncm{\ba}{\begin{array}}
\ncm{\ea}{\end{array}}

\ncm{\bsn}{\bigskip\noindent}

\ncm{\proof}{\bsn{\bf Proof:\ }}
\ncm{\qed}{\hspace*{0.4cm}\rule{0.24cm}{0.24cm}}

\def\AAA{(\A,\onne,\Delta,\e)}
\def\AAAA{(\A,\onne,\Delta,\e,S,\al,\be)}
\ncm{\A}{{\cal A}}
\ncm{\Aop}{{\cal A}_{op}}
\ncm{\AL}{\A_L}
\ncm{\AR}{\A_R}
\ncm{\hA}{\hat\A}
\ncm{\hM}{\hat\M}
\ncm{\AsR}{\A_{\s,R}}
\ncm{\AsL}{\A_{\s,L}}
\ncm{\hAsR}{\hat\A_{\s,R}}
\ncm{\As}{\A_\s}
\ncm{\Asp}{\A_{\s'}}
\ncm{\hAs}{\hat\A_\s}
\ncm{\hAsp}{\hA_{\s'}}
\ncm{\hAssp}{\hA_{\s\s'}}
\ncm{\hAsps}{\hA_{\s'\s}}
\ncm{\Asps}{\A_{\s'\s}}
\ncm{\AspR}{\A_{\s'R}}
\ncm{\Assp}{\A_{\s\s'}}
\ncm{\Ass}{\A_{\s\s}}
\ncm{\ARs}{\A_{R\s}}
\ncm{\ALs}{\A_{L\s}}
\ncm{\ALL}{\A_{LL}}
\ncm{\ALR}{\A_{LR}}
\ncm{\ARL}{\A_{RL}}
\ncm{\ARR}{\A_{RR}}
\ncm{\hAL}{\hat\A_L}
\ncm{\hAR}{\hat\A_R}
\ncm{\hARL}{\hat\A_{RL}}
\ncm{\hALR}{\hat\A_{LR}}
\ncm{\hALL}{\hat\A_{LL}}
\ncm{\hARR}{\hA_{RR}}
\ncm{\hAss}{\hat\A_{\s\s}}

\ncm{\B}{{\cal B}}
\ncm{\C}{{\cal C}}
\ncm{\CL}{{\cal C}_L}
\ncm{\CR}{{\cal C}_R}
\ncm{\Cs}{{\cal C}_\s}
\ncm{\Csp}{{\cal C}_{\s'}}
\ncm{\D}{{\cal D}}
\ncm{\E}{{\cal E}}
\ncm{\F}{{\cal F}}
\ncm{\G}{{\cal G}}
\def\GGG{(\G,\onne,\Delta,\e)}
\def\hG{\hat\G}
\rncm{\H}{{\cal H}}
\def\hH{\hat\H}
\ncm{\I}{{\cal I}}
\ncm{\J}{{\cal J}}
\ncm{\K}{{\cal K}}
\rncm{\L}{{\cal L}}
\ncm{\M}{{\cal M}}

\ncm{\N}{{\cal N}}
\ncm{\NL}{{\cal N}_L}
\ncm{\NR}{{\cal N}_R}
\ncm{\NLL}{{\cal N}_{LL}}
\ncm{\NLR}{{\cal N}_{LR}}
\ncm{\NRR}{{\cal N}_{RR}}
\ncm{\NRL}{{\cal N}_{RL}}
\ncm{\Nssp}{{\cal N}_{\s\s'}}
\ncm{\Nss}{{\cal N}_{\s\s}}
\ncm{\Ns}{{\cal N}_{\s}}
\ncm{\Nsp}{{\cal N}_{\s'}}
\ncm{\Nsps}{{\cal N}_{\s'\s}}
\ncm{\NsL}{{\cal N}_{\s L}}
\ncm{\NsR}{{\cal N}_{\s R}}
\ncm{\NRs}{{\cal N}_{R\s}}
\ncm{\NLs}{{\cal N}_{L\s}}

\def\SB{S_\B}
\ncm{\R}{{\cal R}}
\ncm{\T}{{\cal T}}
\ncm{\Z}{{\cal Z}}

\rncm{\AA}{{\bf A}}
\ncm{\BB}{{\bf B}}
\ncm{\FF}{{\bf F}}
\ncm{\RR}{{\bf R}}
\ncm{\LL}{{\bf L}}
\ncm{\TT}{{\bf T}}
\ncm{\CC}{{\bf C}}
\ncm{\ZZ}{{\bf Z}}
\ncm{\NN}{{\bf N}}

\def\SSS{(S,\,\al,\,\be)}

\ncm{\one}{{\bf 1}}
\def\onne{\one}
\ncm{\hone}{\hat\onne}
\ncm{\honne}{\hat\onne}
\ncm{\onneV}{{\bf 1}_V}
\ncm{\onneW}{{\bf 1}_W}
\ncm{\onneE}{{\bf 1}_\E}
\ncm{\onneA}{{\bf 1}_\A}
\ncm{\onneVxW}{{\bf 1}_{V\x W}}
\ncm{\onneVxExW}{{\bf 1}_{V\x\E\x W}}
\ncm{\onneExVxW}{{\bf 1}_{\E\x V\x W}}
\ncm{\onneWxVxE}{{\bf 1}_{W\x V\x\E}}
\ncm{\onneVxE}{{\bf 1}_{V\x \E}}
\ncm{\onneWxE}{{\bf 1}_{W\x \E}}
\ncm{\onneExW}{{\bf 1}_{\E\x W}}
\ncm{\onneExV}{{\bf 1}_{\E\x V}}

\def\Lra{\Leftrightarrow}
\def\Llra{\Longleftrightarrow}

\def\axiomsl{{\sc Axioms} $\L$ }
\def\axiomsr{{\sc Axioms} $\R$ }

\def\<cros{\,\raise1.5pt\hbox{$\scriptstyle\triangleright$}\!
	   \raise1.9pt\hbox{$\scriptscriptstyle < $}\,}

\def\cros{\,\raise1.9pt\hbox{$\scriptscriptstyle  > $}\!
	  \raise1.5pt\hbox{$\scriptstyle\triangleleft$}\,}
\def\>cros{\cros}

\def\la{\rightharpoonup}
\def\ra{\leftharpoonup}
\def\arr{\la}
\def\arl{\ra}

\def\lef{{\,\hbox{$\textstyle\triangleright$}\,}}

\def\re{\lef}
\def\reS{\,{_S}\!\re}
\def\li{{\,\hbox{$\textstyle\triangleleft$}\,}}
\def\liS{\li\!_S\,}
\def\reli{\,\,\raise1.5pt\hbox{$\scriptstyle\triangleright$}\!
	   \raise1.5pt\hbox{$\scriptstyle\triangleleft$}\,\,}

\def\Rep{\mbox{Rep}\,}
\def\com{\mbox{Cmod}\,}
\def\Ker{\mbox{Ker}\,}
\def\End{\mbox{End}\,}
\def\EndK{\mbox{End}_K\,}
\def\EndA{\mbox{End}_\A\,}

\def\EndAA{\mbox{End}_{\AL\cap\AR}\,}
\def\EndhAA{\mbox{End}_{\hAL\cap\hAR}\,}
\def\HomA{\mbox{Hom}_\A\,}
\def\HomK{\mbox{Hom}_K\,}
\def\EndSA{\mbox{End}_{S(\A)}\,}
\def\EndcA{\mbox{End}^\A\,}
\def\Aut{\mbox{Aut}\,}

\def\Ad{\mbox{Ad}\,}
\def\Ind{\mbox{Ind}\,}
\def\id{\mbox{id}\,}
\def\idA{\mbox{id}_\A}
\def\idV{\mbox{id}_V}
\def\idVW{\mbox{id}_{VW}}
\def\idW{\mbox{id}_W}

\def\Mat{\mbox{Mat}\,}
\def\Hom{\mbox{Hom}\,}

\def\bra{\langle}
\def\ket{\rangle}
\def\o{\otimes}
\def\e{\varepsilon}
\def\he{\hat\e}
\def\heL{\he_{L}}
\def\heR{\he_{R}}
\def\heLR{\he_{LR}}
\def\heRR{\he_{RR}}
\def\heLL{\he_{LL}}
\def\heRL{\he_{RL}}
\def\hes{\hat\e_\s}

\def\hesp{\hat\e_{\s'}}
\def\hessp{\hat\e_{\s\s'}}
\def\eL{\e_L}
\def\eR{\e_R}
\def\eLR{\e_{LR}}
\def\eLs{\e_{L\s}}
\def\eRs{\e_{R\s}}
\def\eRL{\e_{RL}}
\def\eLL{\e_{LL}}
\def\eRR{\e_{RR}}
\def\es{\e_\s}
\def\esR{\e_{\s,R}}
\def\esL{\e_{\s,L}}
\def\ess{\e_{\s\s}}
\def\esp{\e_{\s'}}
\def\essp{\e_{\s\s'}}
\def\esps{\e_{\s'\s}}
\def\eB{\e_\B}

\def\PL{\sqcap^L}
\def\PR{\sqcap^R}
\def\PLR{\sqcap^{L/R}}

\def\cd{\cdot}
\def\c{\circ}
\def\x{\times}
\def\al{\alpha}
\def\be{\beta}
\def\r{\rho}
\def\rV{\r_V}
\def\rW{\r_W}
\def\rVW{\r_{VW}}
\def\rVoW{\r_{V\o W}}
\def\s{\sigma}

\def\pLR{\pi_{LR}}
\def\psR{\pi_{\s R}}
\def\psL{\pi_{\s L}}
\def\piW{\pi_W}
\def\piV{\pi_V}
\def\piVxW{\pi_{V\x W}}

\newcommand{\om}{\omega}

\rncm{\l}{\lambda}

\def\Del{\Delta}
\def\hDel{\hat\Delta}
\ncm{\0}{_{(0)}}
\ncm{\1}{_{(1)}}
\ncm{\2}{_{(2)}}
\ncm{\3}{_{(3)}}
\ncm{\4}{_{(4)}}
\ncm{\5}{_{(5)}}
\ncm{\6}{_{(6)}}
\ncm{\7}{_{(7)}}
\ncm{\8}{_{(8)}}
\ncm{\9}{_{(9)}}

\def\m1{_{(-1)}}
\def\n2{_{(-2)}}

\begin{document}

{
\rncm{\thefootnote}{\fnsymbol{footnote}}
\title{ \bf Axioms for Weak Bialgebras}

\date{Mai 11, 1998}

\author{{\sc Florian Nill}
\\
{\small Inst. theor. Physik, Univ. M\"unchen, Theresienstr. 39,
D-80333 M\"unchen\footnotemark}  \\
{\small Inst. theor. Physik, FU-Berlin,
Arnimallee 14, D-14195 Berlin, Germany}}

\footnotetext
{$^*$ Present address. \qquad e-mail:  nill@physik.fu-berlin.de
\qquad florian.nill@physik.uni-muenchen.de \\
supported by DFG, SFB 288 {\em Differentialgeometrie und
  Quantenphysik}}

\maketitle

}

\vspace{-7cm}
\rightline{\normalsize q-alg/9805104}
\vspace{6.7cm}

\begin{abstract}
Let $(\A,\onne)$ be a finite dimensional
unital associative algebra over a field $K$, which is also equipped
with a coassociative counital coalgebra structure $(\Delta,\e)$.
$\A$ is called a {\em weak bialgebra} if the coproduct
$\Delta:\A\to\A\o\A$ satisfies $\Delta(ab)=\Delta(a)\Delta(b)$.
We do {\bf not} require $\Delta(\onne)=\onne\o\onne$ {\bf nor}
multiplicativity of the counit $\e:\A\to K$.
Instead, we propose a new set of counit axioms,
which are modelled so as to guarantee that $\Rep\A$ becomes a
monoidal category with unit object given
by the cyclic left $\A$-module $\E:=(\A\arr\hone)\subset\hat\A$, where
$\hone\equiv\e$ is the unit in the dual weak bialgebra $\hat\A$.
Under these {\em monoidality axioms} $\E$ and
$\bar\E:=(\honne\arl\A)$ become commuting unital subalgebras of
$\hA$, which are trivial if and only if $\e$ is multiplicative.
We also propose axioms for an antipode  $S:\A\to\A$,
such that the category $\Rep\A$ becomes rigid. $S$ is uniquely
determined, provided it exists. If a monoidal
weak bialgebra $\A$ has an antipode $S$, then its dual $\hA$ is
monoidal if and only if $S$ is a
bialgebra anti-homomorphism, in which case $S$ is also
invertible. In this way we obtain a definition of weak Hopf algebras
which in Appendix A will be shown to be equivalent to the one given
independently by G. B\"ohm and K. Szlach\'anyi. Special
examples are given by the face algebras of T. Hayashi and the
generalised Kac algebras of T. Yamanouchi.
\end{abstract}

\footnotesize
\tableofcontents

\normalsize
\newpage

\sec{Introduction}

In [MS,S] G. Mack and V. Schomerus have introduced the notion of
{\em weak} coproducts
$\Del:\G\to\G\o\G$
on quasi-Hopf algebras $\G$ by allowing $\Del$ to be non-unital,
$\Del(\one)\neq\one\o\one$.
Examples are semisimple quotients of quantum groups at $q=$
roots of unity.
The underlying motivation was to obtain symmetry candidates
$\G$ in  low dimensional quantum field theories.
Technically this may be understood as a Tannaka-Krein like reconstruction
program [Ma,H\"a], starting from the rigid monoidal category of
Doplicher-Haag-Roberts (DHR) endomorphisms on a local observable
algebra $\M$.
In this way one may successfully match the quantum field theoretic
fusion rules with {\em non-integer} (``statistical" or
``categorical" or ``q-") dimensions with those of $\Rep\G$.

The price to pay in this setting is the quasi-coassociativity
of the coproduct $\Del$. Thus, the dual $\hat\G$ of $\G$ is not
an object of the same kind. In particular $\hat\G$ is not even
an associative algebra, which makes it impossible to define an
analogue of the DHR-field algebra $\F=\M\cros\hat\G$.%
\footnote{There is, however, a sensible definition of the
double crossed product $\M\cros\hat\G\cros\G$ [HN1].}
On the other hand, in Ocneanu's approach of recovering ``quantum
symmetries" from (depth 2) Jones inclusions
[Oc,Da,Lo,Szy,EN,Ya,NW] one always expects a concept of
symmetry algebras $\A$, such that the dual $\hA$ is of the same
type (due to the two-step periodicity in any Jones tower).

\bsn
In this work I propose a new axiomatic approach to weak (Hopf)
bialgebras $\AAA$, which strictly meets this duality principle.
In particular, I start from the observation that dualizing the
property $\Del(\one)\neq\one\o\one$ suggests to allow
non-multiplicative counits as well, i.e.
$\e(ab)\neq\e(a)\e(b)$.
On the other hand, I don't give up coassociativity of the
coproduct $\Del$, such that the dual $(\hA,\hone,\hat\Del,\he)$
is of the same kind.
For simplicity - and to make this duality strategy manifest -
throughout I will restrict myself to finite dimensional algebras $\A$ over
a field $K$. A generalization to infinite dimensional settings
together with appropriate topological (like $C^*$- or
von-Neumann algebraic) structures should be a future goal.

\bsn
A first announcement of the present work has been given in 1994
[N2]. Subsequent discussions with H.-W. Wiesbrock and K.-H.
Rehren have soon lead to first applications in Jones theory and
quantum field theory [W,Re].
In 1996 G. B\"ohm and K. Szlach\'anyi [BSz,Sz] independently came up with
very similar ideas. The main progress of the present paper in
comparison with the BSz-approach is that here I propose so-called
{\em (co)monoidality axioms}, the necessity and
consequences of which are discussed individually and
without referring to antipode structures.
Also, the antipode axioms presented here are simpler than those of
B\"ohm-Szlach\'anyi and are motivated by a more general
analysis of {\em rigidity structures} on weak bialgebras.
In this way I will end up with a set of axioms for weak Hopf algebras, which
will be shown in Appendix A to be equivalent to those of
[BSz,Sz]. Also, the face algebras of T. Hayashi [Ha] and the
generalized Kac algebras of T. Yamanouchi [Ya] are special
kinds of weak Hopf algebras, see Sect. 5 and Sect. 8,
respectively.

Meanwhile G. B\"ohm, K. Szlach\'anyi and I
have exchanged and unified our ideas.
Parallel to this work
we present further common results on the theory of integrals and
$C^*$-structures on finite dimensional weak Hopf algebras in [BNS].
Moreover,
in [NSW] we develop a theory of (co)actions and
crossed products by weak Hopf
algebras and generalize Ocneanu's ideas by showing that any
reducible finite index and depth-2 Jones extension of von-Neumann algebras
with finite dimensional centers is given by a crossed product
with a weak Hopf algebra $\A$.

\bsn
In future work [HN2] we will also clarify the role of
our coassociative weak Hopf algebras as a symmetry in the
quasi-coassociative quantum field theoretic scenario of [MS,S].
Another exciting application will emerge from the
fact [N3] that the observable algebras of
a large class of physical quantum chain models naturally
acquire a weak Hopf algebra structure. In particular, the Hopf
spin models or lattice current algebras of [NSz,AFFS] on an
open chain (of even number of sites) are self-dual weak Hopf
algebras $\A$ and their periodic extensions by one link joining the
endpoints are given by the (weak) quantum double $\D(\A)$. For a first
sketch see also the remarks following Example 3 in Appendix D.
It will be challenging to identify the vacuum representation
of these models with the GNS-representation obtained from the
counit (i.e. the monoidal unit in $\Rep\D(\A)$) and relate their
DHR-theory with the braided rigid monoidal structure of $\Rep\D(\A)$.
As a further interesting conjecture one may suggest similar
applications
\footnote{based on infinite dimensional weak Hopf algebras}
in conformal quantum field theory, such that the quantum field
theoretic fusion rules are reproduced by a weak  Hopf algebra
structure on the obsevables.

\bsn
The plan of this paper is as follows.
In Section 2 we analyse the so-called {\em monoidality axioms}
for the counit $\e\equiv\hat\one\in\hA$
on a weak bialgebra $\A$, making $\Rep\A$ a
monoidal category with nontrivial cyclic unit object
$\A\arr\hat\one\subset\hA$. The dual analogues are the {\em
comonoidality axioms} studied in Section 3. In particular, we
will see that the canonical
left (right) action of $\hA$ on
$\one_\A$ induces nontrivial subspaces
$\AL:=\one_\A\arl\hA\subset\A$ and
$\AR:=\hA\arr\one_\A\subset\A$, which in a comonoidal weak bialgebra are
in fact commuting unital subalgebras of $\A$. Considering $\A$
as a left (right) comodule algebra over itself we also show
that $\A$ is comonoidal if and only if the left (right)
coinvariants of $\A$ are given by $\A_{L/R}$, respectively.
In Section 4 we describe the category $\com\A$ of right
$\A$-comodules and show that if $\A$ is comonoidal then
$\AR$ is the
unit object in $\com\A$. Moreover,
the endomorphism algebra of this comodule is shown to be given by
$\EndcA\AR=\AL\cap\AR$, acting by multiplication on $\AR$. For
weak Hopf algebras this has been noticed before in [Sz].
In Section 5 we generalize an observation of [Sz] by showing
that in bimonoidal weak bialgebras the subalgebras $\A_{L/R}$
are separable $K$-algebras.
In Section 6 we adapt ideas developped  for
quasi-Hopf algebras by Drinfel'd [Dr]
to formulate a theory of {\em rigidity structures} on monoidal
weak bialgebras. This will help to motivate our antipode axioms
in Section 7, where we will see that in
{\em bimonoidal} (i.e. monoidal and comonoidal) weak bialgebras
an antipode $S$ always provides a rigidity structure.
Thus, in Section 8 we define a weak Hopf algebra to be a
bimonoidal weak bialgebra with
antipode $S$. One of our main results here will be that a weak
bialgebra $\A$ with antipode $S$ is a weak Hopf algebra (i.e.
bimonoidal) if and only if $[\AL,\AR]=0$ and $S$ is a bialgebra
anti-morphism.

In Appendix A we relate the present approach to the axioms of
[BSz,Sz]. Appendix B gives more details on rigidity structues
in the spirit of [Dr]. In particular this leads to a proof that
on (finite dimensional!) monoidal weak bialgebras $\A$
rigidity maps $S:\A\to\A$ (i.e. ``quasi-antipodes") are always
invertible
\footnote{Presumably a similar proof also works for finite
dimensional quasi-Hopf algebras.}.
In Appendix C we analyse {\em minimal} (comonoidal) weak
bialgebras $\A=\AL\AR$,
which are defined to be generated by the commuting subalgebras $\AL$ and
$\AR$, as well as their {\em cominimal} dual analogues $\hA$.
In Appendix D we give several examples, most noteworthy a
two-sided crossed product construction of a minimal weak
bialgebra $\A=\AL\o\AR$ with a Hopf algebra $\G$. This example
puts a weak Hopf algebra structure on the Hopf algebraic quantum
chains considered in [NSz,AFFS], such that $\AL$ and $\AR$
become the left and right ``wedge algebras", respectively, of
these models.

\bsn
{\bf Note added:}
After finishing this paper I have been informed by L. Vainerman
that presumably the notion of a {\em quantum groupoid} [M, V,
NV] is equivalent to that of a weak Hopf algebra with
involutive antipode.

\sec{Monoidal Weak Bialgebras}

Throughout all spaces are assumed finite dimensional over a fixed
field $K$. The dual of a linear space $V$ is denoted as $\hat
V=\Hom_K(V,K)$ and the center of an algebra $\A$ is denoted by
$\C(\A)$.

\Definition{2.1}
{A {\em weak  bialgebra} $(\A,\onne,\Delta,\e)$ is an
associative unital algebra $(\A,\onne)$ together with a
coassociative coproduct $\Delta:\A\to\A\o\A$ and a counit $\e:\A
\to K$ for $\Delta$, such that
$\Delta (ab)=\Delta(a)\Delta(b),\ \forall a,b\in \A$.
}
As opposed to ordinary bialgebras we do
{\bf not} require $\Delta(\onne) =\onne\o\onne$
{\bf nor} its dual version $\e (ab)=\e(a)\e(b)$.
Clearly, the dual $\hA$ of $\A$ also is a weak
bialgebra $(\hA,\hone,\hDel,\he)$ with structure maps
\beanon
\bra\phi\psi\mid a\ket &:=& \bra\phi\o\psi\mid\Delta(a)\ket\\
\bra\hone\mid a\ket &:=& \e(a)\\
\bra\hDel (\phi)\mid a\o b\ket &:=& \bra\phi\mid ab\ket\\
\he(\phi) &:=& \bra\phi\mid\onne\ket
\eeanon
where $\phi,\psi\in\hA,\ a,b,\in\A$ and where
$\bra\cdot\mid\cdot\ket$ denotes the dual pairing $\hA\o\A\to K$.
We denote elements of $\A$ by $a,b,c,...$ and
elements of $\hA$ by $\phi,\psi,\xi,...\ $.
We also use standard Hopf algebra notations like
$\Delta(a)=a\1\o a\2,\ (\Delta\o id)(\Delta(a))
\equiv (id\o\Delta)(\Delta (a)) = a\1\o a\2 \o a\3$, etc.,
where a summation symbol and summation  indices are suppressed.
The canonical left and right actions of $\A$ on $\hA$
are denoted by
\bea
\lb{2.0a} a\arr \phi &:=& \phi\1 \bra\phi\2\mid a\ket\\
\lb{2.0b} \phi\arl a &:=& \bra\phi\1\mid a\ket \phi\2
\eea
and similarly for $\A$ and $\hA$ interchanged.
Acting in particular on $\one\in\A$ and $\hat\one\in\hA$,
respectively, we obtain nontrivial linear subspaces
$\A_{L/R}\subset\A$ and $\hA_{L/R}\subset\hA$ given by
$$
\ba{rccclcrcccl}
\AL &:=&\one\arl\hA &\subset& \A \quad &,&  \quad
\AR &:=&\hA\arr\one &\subset& \A
\\
\hAL &:=&\hat\one\arl\A &\subset& \hA \quad &,&  \quad
\hAR &:=&\A\arr\hat\one &\subset& \hA
\ea
$$
Let us now consider the category $\Rep \A$ of finite dimensional
unital representations $\pi_V:\A\to\End V$. We also use the
module language by writing
$\pi_V(a) v\equiv a\cdot v,\ a\in\A,\ v\in V$.
If $\A$ is a weak bialgebra then $\Rep\A$ is equipped with a strictly
associative tensor functor $\Rep\A\x\Rep\A\to\Rep\A$ given on the objects by
\bea\lb{2.0i}
V\x W &:=& \onneVxW (V\o W)\\
\lb{2.0ii}
\piVxW &:=& (\piV \o \piW)\circ\Delta
\eea
and on $\A$-linear morphisms by
\beq\lb{2.0iii}
f\x g:= (f\o g)\circ\onneVxW\ ,
\eeq
where $f\in \Hom_\A(V,V'),\ g\in \Hom_\A(W,W')$
and $\onneVxW :=(\pi_V \o \pi_W)(\Delta(\onne))$.
As a special object in $\Rep\A$ we consider the cyclic $\A$-submodule
$\E\equiv\hAR\subset\hA$
with left $\A$-action given by \Eq{2.0a}.
Our aim is to specify additional axioms for weak bialgebras
$\A$, such that $\Rep\A$ becomes a monoidal category with unit object $\E$.
To this end the following notions will be useful

\Definition{2.3}
{A weak bialgebra $(\A,\onne,\Delta,\e)$ is called
\bea
&\bullet& \mbox{{\em left-monoidal}, if}\
\e (abc) =\e (ab\1)\e(b\2 c),\ \forall a,b,c\in\A
\lb{2.3a}\\
&\bullet& \mbox{{\em right-monoidal}, if}\
\e(abc)=\e(ab\2)\e(b\1 c),\ \forall a,b,c\in\A
\lb{2.3b}\\
&\bullet& \mbox{{\em left-comonoidal}, if}\
(\Delta(\onne)\o\onne)(\onne\o\Delta(\onne))
= \onne\1\o\onne\2\o\onne\3
\lb{2.4a}\\
&\bullet& \mbox{{\em right-comonoidal}, if}\
(\onne\o\Delta(\onne))(\Delta(\onne)\o\onne)
= \onne\1\o\onne\2\o\onne\3
\lb{2.4b}
\eea
A weak bialgebra is called (co)monoidal, if it is left- and
right-(co)monoidal, and it is called {\em bimonoidal}, if it is
comonoidal and monoidal.
}
Clearly, $\A$ is (left-, right-) comonoidal if and only if $\hA$ is
(left-, right-) monoidal. We also note the equivalencies:
$\A$ left-(co)monoidal $\Lra \Aop^{cop}$ left-(co)monoidal
$\Lra \Aop$ right-(co)monoidal $\Lra \A^{cop}$ right-
(co)monoidal, where ``op" means opposite multiplication and
``cop" means opposite comultiplication.
If $\e$ is multiplicative, then $\A$ is always monoidal, and if
$\Del(\onne)=\onne\o\onne$, then $\A$ is always comonoidal.
The face algebras of T. Hayashi [Ha] provide
examples of comonoidal weak bialgebras. In fact, we will see in
\cor{3.7} that finite dimensional face algebras are also
bimonoidal.

The terminologies of \defi{2.3} will be motivated below by showing that if
$\A$ is monoidal then $\Rep\A$ is a monoidal category with
unit object given by $\E$. By duality, if $\A$ is comonoidal,
then the category $\com\A$ of $\A$-comodules becomes a
monoidal category, see Sect. 4

Let now $\AAA$ be a weak bialgebra.
For any representation $(\pi_V,V)$ of $\A$ we
introduce the K-linear maps
$L_V : V\to\E\x V$ and $R_V :V\to V\x\E$ given by
\bea
\lb{2.5} L_V(v) &:=& \onneExV (\hone \o v) \\
\lb{2.6} R_V(v) &:=& \onneVxE (v\o \hone)\ .
\eea
These maps satisfy ``naturality" in the sense that
$$
L_W\circ f = (\onneE\x f)\circ L_V
\quad \mbox{and}\quad
R_W \circ f= (f\x\onneE)\circ R_V
$$
for all $f\in \Hom_\A (V,W)$, as well as
\bea
\lb{2.7} \onneExVxW \circ (L_V\o \onneW) &=& L_{V\x W}\\
\lb{2.8} \onneWxVxE \circ (\onneW\o R_V) &=& R_{W\x V}\ .
\eea
Also, we have a kind of ``pre-" triangle identity in the sense that
\beq\lb{2.9}
\onneVxExW \circ (\onneV\o L_W) =\onneVxExW \circ (R_V\o \onneW)
\eeq
as maps $V\o W\to V\x  \E \x W$. Eqs. \no{2.7}-\no{2.9} follow easily from
\beq\lb{2.10}
\onne\1\o\onne\2\o\onne\3 =
\onne\1\o\onne\2\onne_{(1')}\o\onne\3\onne_{(2')} =
\onne\1\onne_{(1')} \o\onne\2\onne_{(2')}\o\onne\3
\eeq
Moreover, $L_V $ and $R_V$ are always injective with left inverses given by
\bea
\lb{2.11} \bar L_V &:& \E\o V\ni \phi\o v\mapsto
\bra\phi\mid\onne\1\ket\onne\2\cdot v\in V
\\
\lb{2.12}
\bar R_V &:& V\o\E\ni v\o\phi \mapsto
\bra\phi\mid\onne\2\ket\onne\1\cdot v\in V
\eea
More precisely we have

\Lemma {2.4}
{Let $(\A,\onne,\Delta,\e)$ be a weak bialgebra. Then for all
$V$ in $\Rep\A$ we have
\beanon
&\bar L_V\onneExV  = \bar L_V \qquad ,\qquad \bar R_V\onneVxE = \bar R_V&
\\
&\bar L_V L_V = \bar R_V R_V =\onneV &
\eeanon
}
\proof
By definition we have for all $\phi\in\E$ and $v\in V$
$$
(\bar L_V\onneExV)(\phi\o v) =
\bra\onne_{(1')} \arr\phi\mid\onne\1\ket\onne\2\onne_{(2')}\cdot v
= \bra\phi\mid\onne\1\onne_{(1')}\ket \onne\2\onne_{(2')}
\cdot v =\bar L_V(\phi\o v).
$$
Hence, we also get $\bar L_VL_V v=\bar L_V(\hone \o v)=v$.
The argument for $R_V$ is anologous. \qed

\bsn
We now show that $L_V$ and $R_V$ are $\A$-linear for all $V$ in
$\Rep\A$ if and only if $\A$ is monoidal, in which case
$L_V$ and $ R_V$ are also bijective.
More precisely, we have the following

\Theorem{2.5}
{Let $(\A,\onne,\Delta,\e)$ be a weak bialgebra. Then \\
i) $L_V$ is $\A$-linear for all $V$ in $\Rep\A$ if and only if
$\A$ is left-monoidal.
\\
ii) $R_V$ is $\A$-linear for all $V$ in $\Rep\A$ if and only if
$\A$ is right-monoidal.
\\
iii) If $\A$ is monoidal, then $L_V$ and $R_V$ are bijective
and we have the identities
\bea
\lb{2.14} &R_V\x \onneW = \onneV \x L_W&
\\
\lb{2.13} &L_V\x\onne_W =L_{V\x W}
\qquad,\qquad
\onneW\x R_V = R_{W\x V}&
\\
\lb{2.13a} &R_\E=L_\E&
\eea
}
\thm{2.5} implies
that for monoidal weak bialgebras $\A$ the category $\Rep\A$
becomes a strictly associative monoidal category with unit object
given by the $\A$-module $\E$.
Note that for ordinary bialgebras the $\A$-module $\E$ is ``trivial",
i.e. it coincides with the 1-dimensional representation given by the counit
$\e:\A\to K$. In our setting $\E$ need not even  be $\A$-irreducible.
If it is, then we call $\A$ {\em pure}, following [BSz].

To prove \thm{2.5} we have to introduce some formalism.
For $\phi\in\hA$ we introduce the maps $\phi_{L/R} :\A\to\hA$ given by
$$
\bra\phi_L(a)\mid b\ket = \bra a\mid\phi_R(b)\ket := \bra\phi\mid ab\ket
$$
for $a,b,\in\A$.
Note the obvious identities
\beq\lb{2.15a}
\ba{rclcrcl}
\phi_R(a) &=& a\arr \phi &\quad,\quad& \phi_L(b) &=& \phi\arl b
\\
\phi_R(ab) &=& a\arr\phi_R(b) &\quad,\quad& \phi_L(ab) &=& \phi_L(a)\arl b
\ea
\eeq
for all $a,b\in\A$ and $\phi\in\hA$.
In particular, $\e_{L/R}(\A)=\hA_{L/R}$ and
$\he_{L/R}(\hA)=\A_{L/R}$.
For $\sigma,\sigma'\in\{L,R\}$ we also use the notation
$$
\essp := \hes \circ \esp \in\EndK\A\qquad,\qquad
\hessp := \es\circ \hesp \in\EndK \hA
$$
where $\he\equiv \onne\in\A$ is the counit on $\hA$. For the
reader's convenience we give the explicit formulas
\beq\lb{2.15}
\ba{rclcrcl}
\eLL (a) &=& \e(a\onne\1)\onne\2 &\quad,\quad&
\eRR (a) &=& \onne\1\e(\onne\2a)
\\
\eLR (a) &=& \e(\onne\1 a)\onne\2 &\quad,\quad&
\eRL (a) &=& \onne\1\e(a\onne\2)
\ea
\eeq
From these one immediately verifies the following identities
\beq\lb{2.16}
\ba{rclcrcl}
a\2 \eLL (ba\1) &=& a\2\e(ba\1) &\quad,\quad&
\eRR(a\2b)a\1 &=& \e(a\2b) a\1
\\
\eLR(a\1 b)a\2 &=& \e(a\1 b)a\2 &\quad,\quad&
a\1\eRL(ba\2) &=& a\1\e(ba\2)
\ea
\eeq
for all $a,b\in\A$. Also note, that the maps $\phi_L$ and $\phi_R$ are
transposes of each other and therefore
\beq\lb{2.17}
(\essp)^t =\he_{-\sigma',-\sigma}
\eeq
where $-L=R$ and $-R=L$, and where the superscript $^t$ denotes
the transposed map.

\smallskip
To prepare the proof of \thm{2.5} we are now going to express
the left- and right- monoidality axioms of \defi{2.3} in
terms of properties the maps $\es,\ \essp$ and $\hessp$.
To this end consider the following list of additional axioms
for weak bialgebras, divided into two groups called \axiomsl
and \axiomsr, respectively.

\bea
\mbox{\axiomsl\qquad} &\quad&  \mbox{\quad\axiomsr}  \nonumber
\\
\nonumber
\\
\lb{2.23}
\phi\1 \o\heLL(\phi\2)=\phi\hone\1 \o \hone\2 && \phi\1\o\heLR(\phi\2)
=\hone\1\phi \o \hone\2
\\
\lb{2.21}
a\1\o\eL(a\2)=a\onne\1\o\eL(\onne\2) &&
\eL(a\1)\o a\2=\eL(\onne\1)\o a\onne\2
\\
\lb{2.19}
a\eRR(b) = \e(a\2 b)a\1\equiv\eRR(a\2b)a\1 &&
a\eLR(b) = \e(a\1b)a\2\equiv\eLR(a\1b)a\2\qquad
\\
\nonumber
\\
\lb{2.24}
\heRR (\phi\1)\o \phi\2 =\hone\1 \o \hone\2 \phi &&
\heRL(\phi\1)\o\phi\2 = \hone\1\o \phi\hone\2
\\
\lb{2.22}
\eR(a\1)\o a\2=\e_R(\onne\1)\o\onne\2 a &&
a\1\o \e_R(a\2)=\onne\1 a\o\e_R(\onne\2)
\\
\lb{2.20}
\eLL (b) a= a\2\e(ba\1)\equiv a\2\eLL(ba\1) &&
\eRL(b)a =a\1\e(ba\2)\equiv a\1\eRL(ba\2)
\eea
where Eqs. \no{2.19}, \no{2.21}, \no{2.20} and \no{2.22}
are supposed to hold for all $a,b,c,\in\A$,
respectively, and Eqs. \no{2.23} and \no{2.24} for all $\phi\in\hA$.
Note that the second identities in \no{2.19} and \no{2.20}
follow from \no{2.16}.
\Proposition{2.6}
{For a weak bialgebra $(\A,\onne,\Delta,\e)$ any one of the
list of \axiomsl (\axiomsr) is equivalent to $\A$ being left-
(right-) monoidal.}
\proof
It is enough to prove the ``left"-statements, since the \axiomsr
reduce to the \axiomsl in $\A^{cop}$. Also note that the axioms
\no{2.24}, \no{2.20} and \no{2.22} reduce to the axioms
\no{2.23}, \no{2.19} and \no{2.21}, respectively, in
$\Aop^{cop}$. Hence, it is enough to prove the equivalences
(\ref{2.3a}) $\Leftrightarrow$ (\ref{2.23}Left) $\Leftrightarrow$
 (\ref{2.19}Left) $\Leftrightarrow$ (\ref{2.21}Left).
To this end first note that (\ref{2.3a}) may be rewritten as
\beq \lb{2.25}
\hone\1\o \hone\2\o\hone\3 =\hone\1\o\hone\2\hone_{(1')} \o\hone_{(2')}
\eeq
implying for all $\phi\in\hA$
\beanon
\phi\hat\onne\1\o\hone\2 &\equiv& \phi\1\hone\1\o\he(\phi\2\hone\2)\hone\3\\
&=& \phi\1\o\he(\phi\2\hone\1)\hone\2\equiv\phi\1\o\he_{LL}(\phi\2)
\eeanon
and therefore (\ref{2.23}Left). Converseley, assume
(\ref{2.23}Left) holds, then
$\hDel(\hone)=\hone\1\o \heLL (\hone\2)$ and therefore
$(\hDel \o id)(\hDel(\hone))=\hone\1\o\hone\2\o\heLL(\hone\3)
=\hone_{(1')}\o\hone_{(2')} \hone\1\o\hone\2$
where the second equation follows by putting $\phi=\hone_{(2')}$
in (\ref{2.23}Left). Hence we have shown
(\ref{2.3a}) $\Leftrightarrow$ (\ref{2.23}Left). The equivalence
(\ref{2.23}Left) $\Leftrightarrow$ (\ref{2.19}Left) follows
by pairing both sides of (\ref{2.23}Left) with $a\o b$ and using
$(\heLL)^t=\eRR$.
Finally, the equivalence (\ref{2.19}Left) $\Leftrightarrow $
(\ref{2.21}Left) follows
from $a\eRR(b)=a\onne\1 \bra \eL(\onne\2)\mid b\ket$ and
$a\1 \e(a\2 b)=a\1\bra\eL(a\2)\mid b\ket$.
\qed

\bsn
For $\s,\s'\in\{L,R\}$ let us now introduce the subspaces
\bleq{Assp}
\Assp:=\essp(\A)\subset\As\subset\A
\eeq
It turns out, that if $\A$ is monoidal, then $\Assp\subset\A$ is a
unital subalgebra. More precisely, we have

\bsn
\Prop{2.6'}{ \quad\\[.1cm]
1.) If $\A$ is left-monoidal, then for all $a,b\in\A$
\beq
\begin{array}{rrclcrcl}
\mbox{i)}&\quad\Delta(\eLL(a)) &=& \eLL(a)\onne\1\o\onne\2
&\ , \ &
\Delta(\eRR(a)) &=& \onne\1\o\onne\2\eRR(a)
\\
\mbox{ii)}&\eLL(b)\eLL(a) &=& \eLL(b\eLL(a))
&\ ,\ &
\eRR(a)\eRR(b) &=& \eRR(\eRR(a)b)
\end{array}
\eeq
In particular, for $\s\in\{L,R\}$,
$(\ess)^2=\ess$ and $\Ass\subset\A$ is a unital
subalgebra.
\\
2.) If $\A$ is right-monoidal, then for all $a,b\in\A$
\beq
\begin{array}{rrclcrcl}
\mbox{i)}&\quad\Delta(\eRL(a)) &=& \onne\1\o\eRL(a)\onne\2
&\ , \ &
\Delta(\eLR(a)) &=& \onne\1\eLR(a)\o\onne\2
\\
\mbox{ii)}&\eRL(b)\eRL(a) &=& \eRL(b\eRL(a))
&\ ,\ &
\eLR(a)\eLR(b) &=& \eLR(\eLR(a)b)
\end{array}
\eeq
In particular, for $\s\neq\s'\in\{L,R\}$, $(\essp)^2=\essp$ and
$\Assp\subset\A$ is a unital subalgebra.
\\
3.) If $\A$ is monoidal, then $[\ALs,\ARs]=0$ for $\s=L$ and $\s=R$.
}

\proof
Part 2.) reduces to part 1.) in $\A^{cop}$ and the right
identities reduce to the left ones in $\A_{op}^{cop}$.
Let now $\A$ be  left-monoidal, then using \no{2.15}
$$
\Delta(\eLL(a)) =\e(a\onne\1)\onne\2\o\onne\3
= \e(a\onne_{(1')})\onne_{(2')}\onne\1\o\onne\2
= \eLL(a)\onne\1\o\onne\2
$$
where in the second identity we have used (\ref{2.22}left).
Applying this to (\ref{2.20}left) yields
$$
\eLL(b)\eLL(a) = \e(b\eLL(a)\onne\1)\onne\2 = \eLL(b\eLL(a))\ .
$$
The remaining statements of part 1.) follow from
$\essp(\onne)=\onne$.
Part 3.) follows, since 1i)+2i) imply for $a\in\ALL$ and
$b\in\ARL$
$$
\Del(ab)=(a\o b)\Del(\onne)=\Del(ba)\,.
$$
Hence, $ab=ba$ since $\Del$ is injective. The identity
$[\ALR,\ARR]=0$ follows in $\A_{op}$.
\qed


\bsn
Next, we study the counit axioms of [BSz,Sz].

\Lemma{2.7}
{Let $(\A,\onne,\Delta,\e)$ be a weak bialgebra and consider
the following BSz--Axioms
\bea\lb{L}
\mbox{\sc l):}&
\quad\e(ab) =\e(a\onne\1)\e (\onne\2 b),&\ \forall a,b\in\A
\\\lb{R}
\mbox{\sc r):}&
\quad\e(ab) =\e(a\onne\2)\e (\onne\1 b),&\ \forall a,b\in\A
\eea
Then the following equivalences hold for $\s,\s'\in\{L,R\},\ \s\neq\s'$
\beanon
\mbox{\sc l)}
&&\Lra\eL(ab)=\eL(\eLL(a)b),\ \forall a,b
\Lra\eR(ab)=\eR(a\eRR(b)),\ \forall a,b
\Lra\es\circ \hes\circ \es =\es
\\
\mbox{\sc r)}
&&\Lra\eL(ab)=\eL(\eRL(a)b),\ \forall a,b
\Lra\eR(ab)=\eR(a\eLR(b)),\ \forall a,b
\Lra\es\circ \hesp\circ \es =\es
\eeanon
}
\proof
The equivalencies {\sc r)} reduce to {\sc l)} in $\A^{cop}$.
The equivalencies {\sc l)} follow from the identities
$\eL(ab)=\e_L(a)\arl\nolinebreak b,\ \eR(ab)=a\arr\e_R(b)$ and
$\e(a\onne\1)\e(\onne_2 b) =\e(\eLL(a)b)=\e(a\eRR(b))$.
\qed

\bsn
The axioms {\sc (l)} and {\sc (r)} of \lem{2.7}
have been proposed as axioms for weak Hopf algebras in [BSz,Sz].
They imply the monoidality properties of \defi{2.3} only
under additional antipode  axioms, see [BSz,Sz] or
\lem{9.8} in Appendix A.
The property ({\sc l}) also appears as a counit axiom in
Hayashi's face algebra theory [Ha].
The monoidality axioms \no{2.3a} and \no{2.3b} are the ones used in
[BNS,NSW]
and they obviously always imply the BSz-axioms of \lem{2.7}.
As has been observed similarly in [BSz,Sz], \lem{2.7} also implies
the following

\Corollary{3.1}
{Under the BSz-axioms {\sc (l)} and {\sc (r)} of \lem{2.7} the
following bilinear forms are nondegenerate for all
$\s,\s'\in\{L,R\}$
\bea\lb{2.7.1}
\AsL\o\AspR\ni(a\o b) &\mapsto & \e(ab)\in K
\\\lb{2.7.2}
\AsL\o\E\ni(a\o\psi) &\mapsto & \bra\psi\mid a\ket\in K
\\\lb{2.7.3}
\hat\E\o\AsR\ni(\phi\o b) &\mapsto & \bra\phi\mid b\ket\in K
\\\lb{2.7.4}
\hat\E\o\E\ni(\phi\o \psi) &\mapsto & \hat\e(\phi\psi)\in K
\eea
where $\E:=\eR(\A)$ and $\hat\E:=\eL(\A)$. Morover, $\hat\E$
as a right $\A$-module is
dual to $\E$, i.e. for all
$\phi\in\hat\E,\ \psi\in\E$ and $a\in\A$
\bleq{2.7.5}
\hat\e(\phi(a\arr\psi))=\hat\e((\phi\arl a)\psi)
\eeq
}
\proof
The nondegeneracy of \no{2.7.1} - \no{2.7.4} follows
immediately from \lem{2.7}.
To prove \no{2.7.5} write $\phi=\eL(b)$ and $\psi=\eR(c)$.
Then $a\arr\psi=\eR(ac)$ and $\phi\arl a=\eL(ba)$ and therefore
$$
\hat\e(\phi(a\arr\psi))=\e(b\eRR(ac))=\e(\eLL(ba)c)=\hat\e((\phi\arl a)\psi)
\quad\qed
$$

\bsn
Note that \cor{3.1} in particular implies
\bleq{2.7.6}
\dim\Assp=\dim\E=\dim\hat\E,\quad\forall\s,\s'\in\{L,R\}.
\eeq
After these preparations we are now in the position to give the

\bsn
{\bf Proof of \thm{2.5}:}\\
Throughout we use that being finite dimensional $\A$ as a left
$\A$-module is itself an object in $\Rep\A$.
Hence, $L_V$ is $\A$-linear for all
$V$ in $\Rep\A$ if and only if
$$
(\onne\1\arr \hone) \o\onne\2 a=(a\1\arr \hone) \o a\2
$$
for all $a\in\A$, which is precisely the condition (\ref{2.22}Left).
Similarly, $R_V$ is $\A$-linear for all $V$ in $\Rep\A$ if and only if
$$
\onne\1 a\o(\onne\2 \arr \hone) =a\1\o (a\2\arr \hone)
$$
for all $a\in\A$, which is precisely the condition
(\ref{2.22}Right).
To prove part iii) let now $\A$ be monoidal implying all
Eqs. \no{2.23} - \no{2.22} as well as those of \prop{2.6'} and \lem{2.7}.
Using $\E=\eR(\A)$ we get for $\phi=\eR(b) \in \E$ and $v\in V$
\beanon
L_V \bar L_V (\phi\o v) &=& L_V(\heL(\phi) \cdot v) =\eR(\onne\1)\o\onne\2
\eLR (b)\cdot v\\
&=& \eR (\eLR(b)\1)\o \eLR (b)\2 \cdot v\\
&=& \eR(\onne\1 \eLR (b))\o\onne\2\cdot v\\
&=& \big[\onne\1\arr \eR(\eLR(b))\big]\o\big[\onne\2\cdot v\big]\\
&=&\big[\onne\1\arr\phi\big]\o\big[\onne\2\cdot v\big]
\equiv \onne_{\E\times V} (\phi\o v)
\eeanon
Here we have used Eq. (\ref{2.22}Left) in the second line, part
(2.i) of \prop{2.6'} in the third
line, \Eq{2.15a} in the
fourth line and \lem{2.7}({\sc r} in the last line. Repeating this
proof in $\A^{cop}$ yields $R_V\bar R_V =\onneVxE$.
Finally, $\A$-linearity implies
$\onneVxExW\circ (\onneV \o L_W)=\onneV\x L_W$ and
$\onneVxExW \circ (R_V\o \onneW)=R_V\x\onne_W$.
Hence, the triangle identity \no{2.14} follows from
Eq. \no{2.9}. The remaining Eqs. \no{2.13} and \no{2.13a} follow by standard
arguments for any monoidal category (see e.g. [Ka, Sec. XI.2.2]).
\qed

\bsn
In view of \thm{2.5} we will from now on denote
$$
\pi_\e:\A\to\EndK\E,\ \pi_\e(a)\phi:=a\arr\phi
$$
as the
``trivial" or unit representation of $\A$. The following Corollary
states that $\pi_\e$ may equivalently be realized as a left
$\A$-action on $\A_{\s,R},\ \s\in\{L,R\}$, such that
$\pi_\e|_{\A_{\s,R}}$ becomes the left multiplication on itself.

\Corollary{2.8'}%
{Let $\A$ be left (right) monoidal and for $\s=R\ (\s=L)$ let
$\pi_{\s,R}:\A\to\EndK\AsR$ be given by
$\pi_{\s,R}(a)b:=\esR(ab)$. Then $\pi_{\s,R}$ is a
representation of $\A$
satisfying $\pi_{\s,R}(a)b=ab,\ \forall a,b\in\AsR$, and
$\e_R:\AsR\to\E$ is an $\A$-linear isomorphism
with inverse $\hes:\E\to\AsR$, i.e. for all $a\in\A$
\bea\lb{pLR}
\e_R\circ\pi_{\s,R}(a)&=&\pi_\e(a)\circ\e_R|_{\AsR}
\\\lb{pLR'}
\hat\e_\s\circ\pi_\e(a) &=&\pi_{\s,R}(a)\circ\hat\es|_\E
\eea
}
\proof
This follows immediately from  \lem{2.7}.
\qed

\bsn
In Section 3 we will see that for monoidal weak bialgebras $\A$
the $\A$-submodule  $\E\equiv\eR(\A)$ is also a subalgebra of
$\hA$ and $\eR:\AsR\to\E$ is also an algebra isomorphism (for
$\s=L$) or anti-isomorphism (for $\s=R$).

Finally, we emphasize that in the present
context (i.e. without furher assumptions like e.g. existence of
an antipode) monoidality is {\bf not} a selfdual
concept for weak bialgebras.
The following Lemma provides the conditions under which a monoidal
weak bialgebra is also comonoidal (a comonoidal weak bialgebra is
also monoidal).

\Lemma{2.8}
{\quad\\[.1cm]
i) A left-monoidal weak bialgebra $\AAA$ is left-comonoidal if and only if
\beq\lb{2.26}
\Delta(\onne)=(id\o\e\o id)[(\Delta(\onne)\o\onne)(\onne\o \Delta(\onne)]
\eeq
ii) A right-monoidal weak bialgebra $\AAA$
is right-comonoidal if and only if
\beq\lb{2.27}
\Delta(\onne)=(id\o\e\o id)[(\onne\o\Delta (\onne))(\Delta(\onne)\o\onne)]
\eeq
iii) A left-comonoidal weak bialgebra $\AAA$ is left-monoidal if and only if
\beq\lb{2.28}
\e(ab)=\e(a\onne\1)\e(\onne\2 b)
\eeq
iv) A right-comonoidal weak bialgebra $\AAA$
is right-monoidal if and only if
\beq
\e(ab)=\e(a\onne\2)\e(\onne\1 b)
\eeq
}
\proof
Part ii) reduces to part i) in $\A_{op}$ and iii), iv) are the
dual versions of i), ii). To prove part i)
first note that by the counit property \no{2.4a} implies \no{2.26}.
On the other hand, \Eq{2.26} is equivalent to
$\he(\phi\psi)=\he(\phi\hone\1)\he(\hone\2\psi)$,
$\forall\phi,\psi\in\hA$, and therefore, by \lem{2.7}(a), to
$\heL\circ\eL\circ\heL=\heL$.
If in this case $\A$ is also left-monoidal then we may apply
$id\o\heL$ to Eq. (\ref{2.23}Left) to get
$$
\phi\1\o\heL(\phi\2)=\phi\hone\1\o\heL(\onne\2),~~\forall\phi\in\hA
$$
which is precisely the dual version of the condition (\ref{2.21}Left).
Hence, in this case $\hA$ is left-monoidal
and therefore $\A$ is left-comonoidal.
\qed

\bsn
\lem{2.8} implies that the face algebras
of T. Hayashi [Ha] are left monoidal, since by definition they are
comonoidal and satisfy \no{2.28}. In fact, they are even
bimonoidal, since in \cor{3.7} we
will see that comonoidal weak bialgebras are left monoidal if and
only if they are right monoidal.
A comonoidal weak bialgabra which is not monoidal
will be given in Example 1 of Appendix D.

\sec{Comonoidal Weak Bialgebras}

In the previous Section we have emphasized the monoidality
axioms for weak bialgebras $\A$ by relating them to the
monoidality properties of $\E\equiv\eR(\A)$ and
$\hat\E\equiv\eL(\A)$ as objects in $\Rep\A$ and $\Rep\A_{op}$,
respectively.
In this Section we pass to the dual point of view by
investigating the comonoidality axioms \no{2.4a} and \no{2.4b}
and relating them to algebraic properties of the linear
subspaces $\A_{L/R}\subset \A$.
Note that these are just the dual counterparts of $\E$ and
$\hat\E$, respectively, and that for $\s,\s'\in\{L,R\}$ we have
$\As\supset\Assp$.
We will show that for
comonoidal weak bialgebras $\As=\Assp$ and that the spaces $\As,\ \s=L,R$,
are in fact commuting unital subalgebras of $\A$ as in the weak
Hopf setting of [BSz].
We will see that these algebras coincide with the
``fixed point" subalgebras of $\A$ under the natural (left or
right, respectively) action of $\hA$ on $\A$.
We will also show, that $\es:\Asp\to\hAssp$ provides an algebra
isomorphism for $\s\neq\s'$ and an algebra anti-isomorphism for
$\s=\s'$.%
\footnote
{Recall from (the dual of) \prop{2.6'} that
$\hAssp\subset\hA$ is also a subalgebra, if $\A$ is comonoidal.
}
If $\A$ is bimonoidal, then also $\hAssp=\hAs$.

These results
play an important role in the theory of crossed
products by weak Hopf algebras [NSW]. In regular crossed products of
von-Neumann algebras $\M$ by weak Hopf algebra actions
$\re:\A\o\M\to\M$ one
requires $\A\re\onne_\M=\A_L\re\onne_\M\cong\A_L$ and
$\C(\M)=(\A_L\cap\AR)\re\onne_\M\cong\A_L\cap\AR$
implying the dual properties in $\M\cros\A$, see [NSW].
In this way the
algebras $\As$ and $\hAs$ appear as the lowest relative commutants in the
resulting reducible Jones tower [NSW].
 If $\A$ is a Frobenius weak Hopf algebra, then $\As$ parametrizes the
space of integrals in $\A$ [BNS].

\bsn
Let now $\AAA$ be a weak bialgebra.
We start with introducing the following four unital subalgebras
$\Nssp (\A) \subset \A,\ \sigma,\sigma'\in\{ L,R\}$, given by
\bea
\lb{3.3} \NLL(\A) &:=& \{ a\in\A\mid\Delta(a)=a\onne\1\o\onne\2\} \\
\lb{3.4} \NLR(\A) &:=& \{a\in\A\mid\Delta(a)=\onne\1 a\o\onne\2\}\\
\lb{3.5} \NRL(\A) &:=& \{a\in\A\mid\Delta(a)=\onne\1\o a\onne\2\}\\
\lb{3.6} \NRR(\A) &:=& \{a\in\A\mid\Delta(a)=\onne\1\o\onne\2a\}
\eea
The subalgebras $\Nssp(\hA)\subset \hA$ are defined accordingly.
These algebras may be considered as the
``left and right fixed point subalgebras"
of $\A$ under the canonical $\hA$-actions in the following sense
\bea
\lb{3.7}
\NLL(\A) &=& \{a\in\A\mid\phi \arr(ab)=a(\phi\arr b),
\ \forall \phi\in\hA,\ \forall b\in\A\}
\\ \lb{3.8}
\NLR(\A) &=& \{a\in\A\mid \phi\arr(ba)=(\phi\arr b)a,
\ \forall\phi\in\hA,\ \forall b\in\A\}
\\ \lb{3.9}
\NRL(\A) &=&\{ a\in\A\mid (ab)\arl \phi=a(b\arl\phi),
\ \forall\phi\in\hA,\ \forall b\in\A\}
\\ \lb{3.10}
\NRR(\A) &=&\{ a\in\A\mid (ba)\arl\phi =(b\arl\phi)a,
\ \forall\phi\in\hA,\ \forall b\in\A\}
\eea
Some immediate consequences of the above definitions are given
in the following

\Corollary {3.2}
{For any weak bialgebra $(\A,\onne,\Delta,\e)$ and for all
$\s,\s'\in\{L,R\}$ we have
$$
\ba{rl}
i)&
a\in\Nssp(\A) \Longrightarrow \essp(a)=a.
\\
&\mbox{In particular}\ \Nssp(\A)\subset \Assp,\
\mbox{and}\ \Nssp(\A)=\Assp\ \mbox{if $\A$ is monoidal.}
\\
ii)&
\NsL(\A)\cap\C(\A) = \NsR(\A)\cap\C(\A) =: \Cs(\A)
\\
iii)& [\NLs(\A),\NRs(\A)]=0.
\ea
$$}
\proof
To prove part (i) use \no{2.15} and \prop{2.6'}, and for part
(ii) use Eqs. \no{3.7} - \no{3.10}.
For $\s=L$ part (iii) follows from
$\Delta(ab)=(a\o b)\Delta(\onne)=\Delta(ba),\ a\in\NLL(\A),\
b\in\NRL(\A)$,
by applying $\e\o\id$. The argument for $\s=R$ is analogous.
\qed

\bsn
Our next aim is to show that
$\es$ maps $\Nsps(\A)$
(anti-)isomorphically onto $\Nssp(\hA)$.
To this end we first need the following

\Lemma{3.3}
{Let $(\A,\onne,\Delta,\e)$ be a weak bialgebra. Then the
following equivalencies hold
$$
\ba{rrcccl}
\mbox{i)}& a\in \NLL(\A) &\Longleftrightarrow &
\phi\arr a=a\heR(\phi),\ \forall\phi\in\hA &\Llra&
ab=b\arl \eL(a),\ \forall b\in\A
\\&&\Llra& \phi\arl a=\eL(a)\phi,\ \forall\phi\in\hA
\\
\mbox{ii)} &a\in\NLR(\A)&\Longleftrightarrow &
\phi \arr a=\heR(\phi)a,\ \forall\phi \in\hA
&\Longleftrightarrow& ba=b\arl\eR(a),\ \forall b\in\A
\\&&\Llra& a\arr\phi =\eR(a)\phi,\ \forall\phi\in\hA
\\
\mbox{iii)}& a\in\NRL(\A) &\Longleftrightarrow&
a\arl\phi =a\heL(\phi),\ \forall\phi\in\hA
&\Longleftrightarrow&  ab=\eL(a)\arr b,\ \forall b\in\A
\\&&\Llra& \phi\arl a =\phi\eL(a),\ \forall\phi\in\hA
\\
\mbox{iv)}& a\in\NRR(\A) &\Longleftrightarrow&
a\arl\phi =\heL(\phi) a,\ \forall \phi\in\hA
&\Longleftrightarrow &  ba =\eR(a)\arr b,\ \forall b\in\A
\\&&\Llra& a\arr\phi=\phi\eR(a),\ \forall\phi\in\hA
\ea
$$
}
\proof
The equivalences (ii), (iii) and (iv) reduce to (i) in
$\A_{op},\ \A^{cop}$ and $\A_{op}^{cop}$, respectively.
To prove (i) first note
that the equivalence
$a\in\NLL(\A)\Lra(\phi\arr a)=a(\phi\arr\onne)\equiv a\heR(\phi)$ is
obvious from the definition \no{3.3}, see also \Eq{3.7}.
Let now $a\in\NLL(\A)$ then for all
$b\in\A$
$$
ab= \e(a\1b\1)a\2b\2=\e(ab\1)b\2=\bra\eL(a)\mid b\1\ket b\2=b\arl\eL(a).
$$
Pairing both sides of this condition with $\phi\in\hA$ we further get
$$
\bra\phi\arl a\mid b\ket =\bra\eL(a)\phi\mid b\ket,\ \forall b\in\A
$$
and therefore $\phi\arl a=\eL(a)\phi$.
Finally, if this holds for all $\phi\in\hA$,
then we get for all $\phi,\psi\in\hA$
\beanon
\bra\phi\o\psi\mid a\onne\1\o\onne\2\ket &\equiv&
\he((\phi\arl a)\psi) =\he (\eL(a)\phi\psi)\\
=\he((\phi\psi)\arl a) &\equiv&
\bra\phi\psi\mid a\ket \equiv \bra\phi\o\psi\mid\Delta (a) \ket
\eeanon
implying $a\in\NLL(\A)$.
\qed

\Theorem{3.4}
{For any weak bialgebra $\A$ and for all
$\sigma,\sigma'\in\{L,R\}$ we have
\\
i) $\es$ maps $\Nsps(\A)$
bijectively onto $\Nssp(\hA)$ with inverse
given by $\hesp : \Nssp(\hA)\to\Nsps (\A)$.
Moreover, these are algebra
isomorphisms, if $\sigma\not= \sigma'$, and algebra anti-isomorphisms, if
$\sigma=\sigma'$.
\\
ii)Defining $\Cs(\A):=\C(\A)\cap\NsL(\A)=\C(\A)\cap\NsR(\A)$
as in \cor{3.2}ii) we have
\bea\lb{3.11}
\es(\NLs(\A)\cap\NRs(\A)) &=& \Cs(\hA)
\\ \lb{3.12}
\heL(\Cs(\hA))=\heR(\Cs(\hA)) &=& \NLs(\A)\cap\NRs(\A)
\\\lb{3.13}
\es(\CL(\A)\cap\CR(\A)) &=& \CL(\hA)\cap\CR(\hA).
\eea
}
\proof
(i) By passing to $\A_{op},\ \A^{cop}$ or $\A_{op}^{cop}$,
respectively, and noting
$\widehat{\A_{op}} = (\hA)^{cop}$ and $\widehat{\A^{cop}} =(\hA)_{op}$,
it suffices to consider the case $\sigma=\sigma'=L$.
By \cor{3.2}i) $\heL\circ\eL|_{\NLL(\A)}=\id$ and it is enough to show
$\eL(\NLL(\A))=\NLL(\hA)$. To show $\eL(\NLL(\A))\subset \NLL(\hA)$
let $a\in\NLL(\A)$ and put $\psi =\eL(a)$.
Then $a=\heL(\psi)$  by \cor{3.2}i) and the last equivalence in
\lem{3.3}i) implies
$$
\psi\phi = \phi\arl \heL(\psi) ,\ \forall\phi\in\hA
$$
By the second equivalence of \lem{3.3}i) this implies $\psi\in\NLL(\hA)$.
Interchanging the role of $\A$ and $\hA$ we also get
$\heL(\NLL(\hA))\subset \NLL(\A)$ and therefore
equality. Finally, putting $\phi=\eL(b)$ in the last
equivalence of \lem{3.3}i) implies for $a,b\in \NLL (\A)$
$$
\eL(a)\eL(b) = \eL(b)\arl a=\eL(ba)
$$
by \Eq{2.15a}, and therefore $\eL|_{\NLL(\A)}$ is an algebra antimorphism.
\\
(ii) Similar as above, by passing to $\A_{op}^{cop}$ it
suffices to consider the case $\s=L$.
If $a\in\NLL(\A)\cap\NRL(\A)$ then by \lem{3.3}(i) and (iii)
$$
\phi\arl a = \eL(a)\phi =\phi\e_L(a),\qquad\forall\phi\in\hA
$$
implying $\eL(a)\in\C(\hA)\cap\NLL(\hA)=\CL(\hA)$
by part (i).
To show that $\eL : (\NLL(\A)\cap\NRL(\A)) \to\CL(\hA)$ is
surjective pick $\phi\in\CL(\hA)$, then by part (i)
$$
\phi=\eL(\hes(\phi))\in\eL(\NLL(\A)\cap\NRL(\A)),
$$
where we have used $\heL(\phi)=\heR(\phi)$ for all
$\phi\in\C(\hA)$. This proves \no{3.11} and by part (i) also the inverse
relation \no{3.12}.\\
Finally, \Eq{3.13} follows by using \cor{3.2}(ii) to get for
$\s=L$ or $\s=R$
$$
\CL(\A)\cap\CR(\A)=\NLs(\A)\cap\NRs(\A)\cap\C(\A)
$$
(and the same formula with $\A$ replaced by $\hA$) and applying
Eqs. \no{3.11}, \no{3.12} and their dual versions.
\qed

\bsn
The algebra
$\Z:=\C_L(\A)\cap\C_R(\A)\cong\C_L(\hA)\cap\C_R(\hA)$ has been called the
``hyper center" of $\A$ (and $\hA$) in [Sz] and it appears as
$C(\M^\A)\cap\C(\M)=\C(\M)\cap\C(\M\cros\A)$ in the crossed
product theory of [NSW].
If $p\in\Z$ is an idempotent, then $\A_p:=p\A\subset\A$ is a weak
sub-bialgebra and by \no{3.13} its dual is given by
$\widehat{\A_p}=\hat p\hA$, where $\hat p=\eL(p)\equiv\eR(p)$.

There is an alternative insight into \lem{3.3} and \thm{3.4} by
considering $\A$ and $\hA$ as subalgebras of $\EndK\A$.
Let $Q_\s:\A\to \EndK\A$ and  $P_\s:\hA\to \EndK\A$, $\s=L,R$,
be given by
$$
\ba{rclcrcl}
Q_L(a)b & := & ab\quad &,& \quad Q_R(a)b & := & ba
\\
P_L(\phi)b &:=& b\arl\phi \quad &,& \quad P_R(\phi)b &:=& \phi\arr b
\ea
$$
where $a,b\in\A$ and $\phi\in\hA$. Then we have the following

\Lemma{3.4}
{For any pair of dual weak bialgebras $\A$ and $\hA$ and for all
$\sigma,\sigma'\in\{L,R\}$ we have
\bea\lb{3.4a}
Q_\s(\Nsps(\A))& = & P_{\s'}(\Nssp(\hA)) = Q_\s(\A)\cap P_{\s'}(\hA)
\\\lb{3.4b}
\es|_{\Nsps(\A)} &=& P_{\s'}^{-1}\circ Q_\s|_{\Nsps(\A)}
\\\lb{3.4c}
\hesp|_{\Nssp(\hA)}&=& Q_\s^{-1}\circ P_{\s'}|_{\Nssp(\hA)}
\eea
}
\proof
\lem{3.3} immediately gives $Q_\s(\Nsps(\A)) =
 P_{\s'}(\Nssp(\hA)) \subset Q_\s(\A)\cap P_{\s'}(\hA)$, as well as
the identities \no{3.4b} and \no{3.4c}.
Conversely, if $Q_L(a)=P_L(\psi)$ then
$\phi\arr(ab)=\phi\arr b\arl\psi=a(\phi\arr b)$ for all
$b\in\A,\ \phi\in\hA$, implying $a\in\N_{LL}(\A)$ by \no{3.7}.
The remaining cases are analogous.
\qed

\bsn
We now show that a weak bialgebra
$(\A,\onne,\Delta,\e)$ is comonoidal if and only
if $\Nssp(\A)=\As$ for all $\sigma,\sigma'\in\{L,R\}$.
Again, this statement may be divided into two pieces.

\Theorem{3.5}
{For  a weak bialgebra $(\A,\onne,\Delta,\e)$ the following
equivalencies hold
\\[.1cm]
i) $\A$ is left-comonoidal
$\Llra\AL=\NLL(\A)\Llra\AR=\NRR(\A)$.\\
If this holds then we also have
$\Nss(\hA)=\hAss$ and $\Nss(\A)=\Ass$,
for $\sigma=L$ and $\sigma=R$.
\\
ii) $\A$ is right-comonoidal $\Llra\AL=\NLR(\A)\Llra\AR=\NRL(\A)$.\\
If this holds then we also have
$\Nssp(\hA)=\hAssp$ and $\Nssp(\A)=\Assp$ for
$(\s,\s')=(L,R)$ and $(\sigma,\sigma')=(R,L)$.
}
\proof
Part (ii) reduces to part (i) in $\A^{cop}$.
To prove part (i) for $\s=L$ observe that
$\NLL(\A)=\AL$ is equivalent to
$\Delta(\heL(\phi))=\heL(\phi)\onne\1\o\onne\2,\
\forall\phi\in\hA,$
and therefore to
$$
\bra\phi\mid\onne\1\ket\onne\2\o\onne\3 =
\bra\phi\mid\onne\1\ket\onne\2\onne_{(1')}\o\onne_{(2')},
\qquad\forall\phi\in\hA
$$
which is the left-comonoidality property \no{2.4a} for $\A$.
Next, we use that
$\Nssp(\A)\subset\Assp\subset \As$ always holds by
\cor{3.2}i). Hence $\NLL(\A)=\AL$ implies
$\NLL(\A)=\ALL$.
To get the dual statement we use that
if $\A$ is left-comonoidal then $\hA$ is
left-monoidal implying $\NLL(\hA)\supset\heLL(\hA)\equiv\hALL$ by
part (1i, left) of \prop{2.6'}.
The case $\sigma=R$ follows by passing to $\A_{op}^{cop}$.
\qed

\Corollary{3.5}
{Let $\A$ be monoidal. Then for $\s\in\{L,R\}$ the restrictions
$\eLs|_{\ARs}$ provide algebra anti-isomorphisms
$\eLs:\ARs\to\ALs$
with inverse $\eRs:\ALs\to\ARs$.
}
\proof
This follows from \thm{3.4} and the dual of \thm{3.5}, implying
$\NsL(\hA)=\NsR(\hA)=\hAs$.
\qed

\bsn
Note that for comonoidal weak bialgebras $\A$ \thm{3.5} still
allows for the possibility\\
$\hAssp\mathop{\subset}\limits_{\not=}\hAs$.
Also note that if $\A$ is comonoidal, then by
\thm{3.5} and \cor{3.2}iii) $\A_L$ and $\A_R$ are commuting
subalgebras of $\A$.
More precisely, we even have

\Corollary{3.3}
{A weak bialgebra $\A$ is comonoidal if and only if it is left-
(or right-) comonoidal and $[\A_L,\A_R]=0$.
In this case we also have
$\A_L\cap\A_R\cong\CL(\hA)\cong\CR(\hA)$ and
$C_{L/R}(\A)=\C(\A)\cap\A_{L/R}$.
}
\proof
By \thm{3.5} and \cor{3.2}iii) comonoidality implies
$[\A_L,\A_R]=0$. Conversely, if $[\A_L,\A_R]=0$, then
$\Nssp(\A)\subset\As$ (\cor{3.2}i)) gives
$\NsL(\A)=\NsR(\A)$ for $\s\in\{L,R\}$ by the
definitions \no{3.3} - \no{3.6}. In this case, by \thm{3.5},
$\A$ is left-comonoidal iff it is right-comonoidal.
The remaining statements follow from \thm{3.4}.
\qed

\bsn
Somewhat surprisingly, under the condition $[\AL,\AR]=0$ also
left- and right-{\em monoidality} become equivalent.

\Lem{3.6}
{Let $\A$ be a weak bialgebra and suppose $[\AL,\AR]=0$.
Then $\A$ is left-monoidal if and only if it is right-monoidal.
}
\proof
By passing to $\A_{op}$ it suffices to prove one direction.
If $\A$ is left-monoidal, then by the duals of \thm{3.5} and
\cor{3.2}i) $\hAL=\NLL(\hA)\supset\NLR(\hA)$ and
$\hAR=\NRR(\hA)\supset\NRL(\hA)$. Since $[\AL,\AR]=0$ implies
$\NsL(\A)=\NsR(\A),\ \s\in\{L,R\}$, we may now use \thm{3.4} to
conclude
$$
\ba{rcccccccl}
\hAL &=&\NLL(\hA)&\cong&\NLL(\A)_{op}&=&\NLR(\A)_{op}&\cong&\NRL(\hA)_{op}
\\
\hAR &=&\NRR(\hA)&\cong&\NRR(\A)_{op}&=&\NRL(\A)_{op}&\cong&\NLR(\hA)_{op}\ .
\ea
$$
Hence, $\dim\hAL=\dim\hAR$, implying also $\hAL=\NLR(\hA)$ and
$\hAL=\NLR(\hA)$, and by \thm{3.5} $\A$ is right-monoidal.
\qed

\bsn
Putting \cor{3.3} and \lem{3.6} together with their dual
versions we arrive at

\Cor{3.7}
{A weak bialgebra $\A$ is bimonoidal $\Lra\A$ is comonoidal and
left- (or right-) monoidal $\Lra\A$ is monoidal and left- (or
right-) comonoidal.
}
As already remarked, \cor{3.7} together with \lem{2.8} imply
that finite dimensional (actually $\dim\A_{L/R}<\infty$ is
sufficient) face algebras in the sense of Hayashi [Ha] are
bimonoidal weak bialgebras.

Let us summarize our findings for comonoidal weak
bialgebras $\A$. Among the algebras
$\Assp\cong(\hat\Asps)_{(op)}$ we are actually left with only
two equivalence classes
\beq\lb{3.14}
\ba{rcccccccl}
\AL &=& \ALL &=&\ALR &\cong &\hARL&\cong &(\hALL)_{op}\\
\AR &=& \ARR &=&\ARL &\cong &\hALR&\cong &(\hARR)_{op}
\ea
\eeq
the isomorphisms being given by the following diagram
\beq\lb{3.14'}
\begin{array}{rcccl}
\Nsps(\hA) &= &\hAsps &\subset  &\hAsp\\
	  &  &\esp \uparrow\downarrow \hes   &  &\downarrow \hes \\
\Nssp(\A)&= & \Assp &=& \As
\end{array}
\eeq
We have $\hes\esp\hes =\hes$ and therefore
$(\hessp)^2 =\hessp$ and $(\esps)^2 =\esps$.
However we may possibly have $\Ker \hes\cap \hAsp\not= 0$ in which case
$\esp \hes\esp \not= \esp$ and
$\hAsps \mathop{\subset}\limits_{\not=}\hAsp$.
This precisely reflects the possibility that comonoidal weak bialgebras
$\A$ may not be monoidal.
Moreover, by the dual of \no{2.7.6}
\bleq{3.18}
\dim\AL=\dim\AR=\dim\hAssp\,,\quad\forall \s,\s'\in\{L,R\}
\eeq
If $\A$ is bimonoidal, then we also
have $\hAs=\hAssp=\Nssp(\hA)$ for all $\s,\s'\in\{L,R\}$ and the above
diagram also holds with $\A$ and $\hA$ interchanged, i.e.
$$
\AL\cong\hAR\cong(\AR)_{op}\cong(\hAL)_{op}
$$
where the isomorphisms are given by $\es:\Asp\to\hAs$ with
inverse $\hesp:\hAs\to\Asp$.


\section{The Comodule Picture}

In this Section we study the category of (right) comodules $\com\A$ of a
weak bialgebra $\A$ and describe its monoidal structure in case $\A$ is
bimonoidal. Of course, since $\com\A =\Rep\hA$,
this could be traced back to the results of Sect. 2.
However, it turns out that in the bimonoidal case the tensor functor in
$\com\A$ is more naturally described by an ``amalgamated" tensor
product, which will then be shown to be equivalent to the
constructions in Sect. 2.
We also generalize a result of [Sz] by showing that for comonoidal $\A$ the self--intertwiner algebra of the ``trivial" $\A$-comodule $\AR$ is given by $\AL\cap\AR$.

As usual, by a right $\A$-comodule we mean a linear space $V$
together with a coaction $\r_V:V\to V\o\A$ satisfying
\bea\lb{3.15}
(\rV\o\idA)\circ\rV &=& (\idV\o\Del)\circ\rV
\\\lb{3.16}
(\idV\o\e)\circ\rV &=& \idV.
\eea
For $v\in V$ we also use the shorthand notation $\rV(v)\equiv
v\0\o v\1$, omitting as usual summation indices and a summation
symbol.
As for ordinary finite dimensional bialgebras, we recall the
one-to-one correspondence between right $\A$-comodules and left
$\hA$-modules given by
\beq\lb{3.17}
\phi\re v:=(\idV\o\phi)(\rV(v))\equiv v\0\bra\phi\mid v\1\ket,
\quad v\in V,\ \phi\in\hA.
\eeq
Based on this observation we get the following
\Proposition{3.7}
{Let $\AAA$ be a comonoidal weak bialgebra. Then any right
$\A$-comodule $V$ naturally becomes an $\AR$-bimodule via
\bea\lb{3.18a}
a_R\cdot v &:=& v\0\e(a_Rv\1)\equiv\eL(a_R)\re v
\\\lb{3.18b}
v\cdot a_R &:=& v\0\e(v\1a_R)\equiv\eR(a_R)\re v,
\eea
where $a_R\in\AR$ and $v\in V$.
Moreover, with respect to this biaction we have
 for all $a_R\in\AR,\ v\in V$
\bea\lb{3.19}
\rV(a_R\cd v)=\Del(a_R)\cd\rV(v)
\\\lb{3.20}
\rV(v\cd a_R)=\rV(v)\cd\Del(a_R)
\\
\lb{3.21}
\eRR(v\1)\cd v\0 = v = v\0\cd\eRL(v\1).
\eea
}
\proof
By \thm{3.4} and \thm{3.5} $\eL:\A_R\to\NLR(\hA)$ is an
algebra isomorphism and $\eR:\A_R\to\NRR(\hA)$ is an algebra
anti-isomorphism. Hence \no{3.18a} and \no{3.18b} provide a
left and a right $\AR$-action, respectively, on $V$, which
commute with each other due to \cor{3.2}(iii).
To prove the identities \no{3.19} and \no{3.20} first note that
they make sense, since in the comonoidal case
$$
\Del(a_R)=\onne\1\o\onne\2a_R=\onne\1\o a_R\onne\2\in\A_R\o\A.
$$
From this \no{3.19} follows by computing
\beanon
\Del(a_R)\cd\rV(v) &=& v\0\e(\onne\1v\1)\o a_R\onne\2v\2=v\0\o a_R v\1
\\
&=&v\0\o\eL(a_R)\arr v\1 = v\0\o v\1\e(a_Rv\2)
\\ &=&\rV(a_R\cd v),
\eeanon
where in the third equation we have used \lem{3.3}(iii) and
$\AR=\NRL(\A)$.
\Eq{3.20} follows analogously from $\AR=\NRR(\A)$ and
\lem{3.3}(iv).
To prove \no{3.21} we compute
$$
\eRR(v\1)\cd v\0=v\0\e(\eRR(v\2)v\1)=v\0\e(v\1)=v
$$
where we have used $\eRR(a\2)a\1=a,\ \forall a\in\A$, by
\no{2.16}. Similarly, using $\eLR(a\1)a\2=a$,
$$
v\0\cd\eRL(v\1)=v\0\e(v\1\eRL(v\2))=v\0\e(v\1)=v\,.\quad\qed
$$
Given two right $\A$-comodules $V,W$, we define
$\r_{V\o W}:V\o W\to V\o W\o\A$ by
\beq\lb{3.22}
\r_{V\o W}(v\o w):=\rV^{13}(v)\rW^{23}(w)\equiv v\0\o w\0\o v\1w\1.
\eeq
One immediatetely checks, that $\rVoW$ again satisfies
\no{3.15}, however it fails \no{3.16} unless $\e$ is
multiplicative. If $\A$ is bimonoidal this may be repaired by using the
$\A_R$-bimodule property to define
\beq\lb{3.23}
\ba{rcl}
VW &:=& V\o_{\AR} W
\\
\rVW(v\o_{\AR} w) &:=&(P_{VW}\o\idA)(\r_{V\o W}(v\o w)),
\ea
\eeq
where $P_{VW}:V\o W\to VW$ is the canonical projection.
Then due to \no{3.19} and \no{3.20}
$$
\rVW: VW\to VW\o\A
$$
is well defined and still satisfies \no{3.15}.
Moreover, we have

\Lemma{3.8}
{Let $\AAA$ be bimonoidal and for two right $\A$-comodules $V,W$
let $\rVW:VW\to VW\o\A$ be given by \no{3.23}. Then
\beq\lb{3.24}
(\idVW\o\e)\c\rVW = \idVW
\eeq
}
\proof
Using $\eL=\eL\c\eRL$ we compute
\beanon
(\idVW\o\e)(\rVW(v\o_{\AR} w))&=& P_{VW}(v\0\o w\0)\e(v\1 w\1)
\\
&=& P_{VW}(v\0\o w\0)\e(\eRL(v\1)w\1)
\\
&=& P_{VW}(v\0\o \eRL(v\1)\cd w)
\\
&=& P_{VW}(v\0\cd\eRL(v\1)\o w)
\\
&=& v\o_{\AR} w,
\eeanon
where in the last line we have used \no{3.21}.
\qed

\bsn
Thus, $\rVW$ is again a right $\A$-coaction. Next, observe that
any $\A$-comodule morphism $f:V\to V'$  (i.e. satisfying
 $\r_{V'}\c f=(f\o\idA)\c\rV$) is also an $\AR$-bimodule map
and therefore the tensor product of two such maps, $f:V\to V'$
and $g:W\to W'$, naturally passes down to an $\A$-comodule morphism
$f\o_{\AR} g :VW\to V'W'$.
In this way $\com\A$ becomes a monoidal category with unit
object given by $\AR$, where
$\r_{\AR}:\AR\to\AR\o\A$ is given by $\r_{\AR}=\Del|_{\AR}$.

Let us now see how, under the identification
$\com\A=\Rep\hA$, this description coincides with the tensor
functor obtained by the dual version of \no{2.0i}--\no{2.0iii}.
To this end we put
\beq\lb{3.25}
V\x W :=(\idV\o\idW\o\e)(\r_{V\o W}(V\o W))\subset V\o W
\eeq
as in \no{2.0i} and correspondingly
\beq\lb{3.26}
\r_{V\x W}:=\r_{V\o W}|_{V\x W}.
\eeq
Then $\r_{V\x W}:V\x W\to(V\x W)\o\A$ is a well defined coaction
satisfying \no{3.15} \underline{and} \no{3.16} and we have

\Lemma{3.9}
{Under the conditions of \lem{3.8} the restriction
$P_{VW}|_{V\x W}:V\x W\to VW$ provides an isomorphism of
$\A$-comodules.
}
\proof
By the definitions \no{3.23} and \no{3.26} $P_{VW}|_{V\x W}$ is a
comodule morphism, which by \lem{3.8} is surjective. To prove
that it is also injective we just have to note that according
to the proof of \lem{3.8} its inverse is given by
$$
VW\ni v\o_{\AR} w\mapsto v\0\o w\0 \e(v\1 w\1)\in V\x W,
$$
which is indeed well defined due to \no{3.19} and \no{3.20}.
\qed

\bsn
We conclude this Section with picking up an observation of [Sz],
who has noticed that in the weak Hopf algebra setting $\A$ is
pure (i.e. the ``trivial" $\A$-module $\E\equiv\hAR$ is
irreducible), if and only if $\C_L(\A)=\C_R(\A)=\CC$.
More generally, the commuant of $\pi_\e(\A)$ in $\EndK\E$ is
given by $\pi_\e(\C_\s(\A))$ [Sz].
It turns out, that this already holds in our setting of
monoidal weak bialgebras. To see this we prove a dual
statement in $\com\A$. First, we need

\Lemma{4.0}
{Let $\AAA$ be a weak bialgebra. Then
$\Delta(\onne)\in(\A_R\o\A)\cap(\A\o\AL)$. If $\A$ is left- or
right-comonoidal, then $\Delta(\onne)\in\AR\o\AL$.
}
\proof
Pick a basis $e_i\in\A$ with dual basis $e^i\in\hA$. Then
$\Delta(\onne)=e^i\arr\onne\o e_i=e_i\o\onne\arl e^i
\in(\A_R\o\A)\cap(\A\o\AL)$.
If $\A$ is left-comonoidal, then
$$
\Delta(\onne)=\onne\1\o\e(\onne\2\onne_{(1')})\onne_{(2')}
=\honne\1\arr\onne\o\onne\arl\honne\2\in\AR\o\AL
$$
The argument for right-comonoidal $\A$ is similar.
\qed

\bsn
Let us now denote the intertwiner spaces in $\com\A$ by
$$
\EndcA V :=\{T\in\EndK V\mid\rV\circ T=(T\o\id)\circ\rV\}
$$
Recall that if $\A$ is comonoidal then $[\AL,\AR]=0$ and $\AR$
is the unit object in $\com\A$, where $\r_{\AR}=\Del|_{\AR}$.

\Lemma{3.10}
{Let $\AAA$ be a comonoidal weak bialgebra and $T\in\EndK\AR$. Then
$T\in\EndcA\AR$ if and only if there exists $z\in\AL\cap\AR$
such that $T(a)= az\equiv za,\ \forall a\in\AR$.}
\proof
We have $T\in\EndcA\AR\Lra\Del(T(a))=(T\o\id)(\Del(a)),\
\forall a\in\A_R.$
Putting $a=\onne$, applying $\e\o\id$ and using \lem{4.0} we
conclude
$$
z:=T(\onne)=\e(T(\onne\1))\onne\2\in\AL\cap\AR.
$$
Now comonoidality implies $\AR=\NRL(\A)$ and therefore
$\Del(a)=(\onne\o a)\Del(\onne),\ \forall a\in\AR$.
Hence, $T(a)=(\e\o\id)(\Del(T(a)))=\e(T(\onne\1))a\onne\2=az,\
\forall a\in\AR$.
Conversely, if $z\in\AL\cap\AR$ then
$\Del(az)=\Del(a)(z\o\onne)$, proving that the map $T:a\mapsto
az$ is in $\EndcA\AR$.
\qed

\bsn
\lem{3.10} may now immediately be dualized. For left
$\A$-modules $V$ denote $\EndA V$ the space of $\A$-linear
endomorphisms of $V$. The following generalizes [Sz, Eq.(3.3)] to our
setting.

\Proposition{3.11}
{Let $\AAA$ be monoidal and denote $\pi_\e:\A\to\EndK\E,\
\E\equiv\hAR$, the ``trivial" representation,
$\pi_\e(a)\phi:=a\arr\phi,\ a\in\A,\,\phi\in\E$.
Then
$$
\EndA\E=\pi_\e(\C_R(\A))=\pi_\e(\C_L(\A))\cong\C_{R/L}(\A)\cong\hAL\cap\hAR
$$
}
\proof
By the dual of \lem{3.10} $T\in\EndA\E$ iff there exists
$\xi\in\hAL\cap\hAR$ such that $T(\phi)=\phi\xi\equiv\xi\phi,\
\forall\phi\in\E$. Using $\E=\NRR(\hA)=\NRL(\hA)$ by \thm{3.5}
we conclude
$$
T(\phi)=\heL(\xi)\arr\phi=\heR(\xi)\arr\phi,\quad\forall\phi\in\E
$$
from the dual versions of \lem{3.3}(iii) and (iv).
Hence, the claim follows since by
the dual of \no{3.11} $\hes(\hAL\cap\hAR)=\C_\s(\A)$ and since
by \thm{3.4}(i) the restriction of $\pi_\e$ to $\NsR(\A)$ -- and
therefore to $\C_\s(\A)\equiv\NsR(\A)\cap\C(\A)$ -- is faithful.
\qed

\bsn
\prop{3.11} in particular implies $\pi_\e(\A)\subset\EndhAA\E$.
In \cor{C3'}i) of Appendix C we will see that equality holds if and only
if $\hAL\o_{\hAL\cap\hAR}\hAR\cong\hAL\hAR$ as a
subalgebra of $\hA$.



\sec{Bimonoidal Weak Bialgebras and Face Algebras}

In this Section we generalize an observation of [Sz] by showing
that in bimonoidal weak bialgebras $\A$ the subalgebras
$\A_{L/R}$ are separable. This will also prove that the
bimonoidal weak bialgebras with abelian $\A_{L/R}$ are
precisely the face algebras of [Ha]. To this end
let us introduce the maps
$S_\sigma:\As\to \A_{-\sigma}$
and $\bar S_\sigma :\As\to \A_{-\sigma}$ given by
\beq\lb{4.6}
\ba{rclcrcl}
S_L &:=& \he_R\circ \e_L |_{\A_L} &\qquad&
S_R &:=& \he_L \circ \e_R |_{\A_R}\\
\bar S_L &:=& \he_R\circ \e_R |_{\A_L} &&
\bar S_R &:=&\he_L\circ \e_L|_{\A_R}\,.
\ea
\eeq
By \thm{3.4}i) and \thm{3.5}, if $\A$ is comonoidal these maps
are algebra anti-isomorphisms and $\bar S_{L/R}=S_{R/L}^{-1}$.
We will see in \cor{9.3'} that if $\A$ is a weak Hopf algebra with
antipode $S$ then $S_{L/R}=S|_{\A_{L/R}}$.

\Lemma{C7}
{Let $\A$ be monoidal or comonoidal and let $a_L,b_L \in\A_L$ and
$a_R, b_R\in \A_R$. Then
\bea
\lb{4.7} &&\e(a_Lb_L) = \e(S_L(a_L)b_L) = \e(a_L\bar S_L(b_L))\\
\lb{4.8} &&\e(a_Rb_R) = \e(\bar S_R(a_R)b_R) = \e(a_RS_R(b_R))
\eea
}
\proof
If $\A$ is monoidal use \lem{2.7}. If $\A$ is comonoidal
use the counit property and $\As=\Nssp(\A)$ to get for
$\sigma =L$ or $\sigma=R$
$$\e(a_\sigma b_\sigma)=\e(a_\sigma\onne\1)\e(\onne\2 b_\sigma)
= \e(a_\sigma\onne\2)\e(\onne\1 b_\sigma)
$$
from which the statements follow by the formulas \no{2.15}.
\qed

\bsn
Next, given a nondegenerate functional $\om:\M\to K$ on a
finite dimensional algebra $\M$ let
$\sum x_i\o y_i \in\M\o\M$ denote the
form-inverse of $\M\o\M\ni(m_1\o m_2)\mapsto\om(m_1m_2)\in K$,
i.e. the unique solution of (summation convention)
\beq\lb{quasibasis}
\om(mx_i)y_i = m = x_i\om(my_i),\quad\forall m\in\M\,.
\eeq
Note that this implies the identity
%
\beq\lb{4.12}
mx_i \o y_i = x_i\o y_i m ,~~~\forall m\in\M.
\eeq
In the terminology of Watatani [Wa] the collection
$\{(x_i,y_i)\}$ would be called a ``quasi-basis" for $\om$.
Generalizing the index notion for conditional expectations we denote
\beq\lb{4.14}
\Ind\om:= x_i y_i\in\C(\M)
\eeq
and call this the {\em Index} of $\om$.
Also recall that
for  a finite dimensional Frobenius algebra $\M$ over a field $K$
 the {\em modular automorphism}
of a non-degenerate functional $\om\in\hat \M$ is defined to be the
unique $\theta_\om\in\Aut \M$ such that
\beq
f(xy)\ =\ f(y\,\theta_\om(x))\ ,\qquad \forall x,y,\in \M\ .
\eeq

\Prop{C6}
{{\rm [BNS]} Let $\A$ be a bimonoidal weak bialgebra.
Then for $\s=L$ and $\s=R$
\\
i)  $\e|_{\As}$ is nondegenerate and $\Ind\e|_{\As}=\onne$.
\\
ii) The quasi-basis
$x^i_\s\o y^i_\s\in\As\o\As$ of $\e|_{\As}$ is given by
\beq\lb{4.10}
 x_L^i \o y_L^i = S_R (\onne\1)\o \onne\2
\qquad\qquad
x_R^i \o y_R^i = \onne\1 \o S_L(\onne\2)
\eeq
iii) The modular automorphisms of $\e|_{\As}$ are given by
$S_R\c S_L\in\Aut\AL$ and
$\bar S_L\c \bar S_R\in\Aut\AR$.
%
\\
iv)  $\AL$ and $\AR$ are separable $K$-algebras, whence semi-simple.
}
\proof
$\e|_{\As}$ is nondegenerate by
\cor{3.1}. To prove (ii) we use
$\Delta(\onne)\in\AR\o\AL$
to compute from  \lem{C7}
\beq
\ba{rcccccl}
\e(a_LS_R(\onne\1))\onne\2
&=&\e(a_L(\bar S_L\c S_R)(\onne\1))\onne\2
&=&\e(a_L\onne\1)\onne\2 &=& a_L\\
\onne\1\e(S_L(\onne\2)a_R) &=&\onne\1\e((\bar S_R\c S_L)(\onne\2)a_R)
&=&\onne\1\e(\onne\2 a_R)&=&a_R
\ea
\eeq
for all $a_L\in\AL$ and $a_R\in\AR$. This proves (ii).
Since by \no{2.16} and the definitions \no{4.6}
$S_R(\onne\1)\onne\2 = \onne\1 S_L(\onne\2)=\onne$,
we also conclude $\Ind\e|_{\As}=\onne$.
Hence, by \no{4.12},
$x^i_\s\o y^i_\s\in\As\o(\As)_{op}$ provides a separating
idempotent, proving part (iv).
Finally, part (iii) also follows from \lem{C7},
since $\A_L$ and $\AR$ commute.
\qed

\bsn
\prop{C6} in particular implies that for abelian
$\A_{L/R}$ our bimonoidal weak bialgebras reproduce the
face algebras in the sense of T. Hayashi [Ha].

\sec{Rigid Weak Bialgebras}

In this Section, adapting ideas of Drinfeld [Dr] for quasi-Hopf
algebras, we propose axioms for
a so-called
{\em rigidity structure} on a monoidal weak bialgebra $\A$,
such that $\Rep\A$
becomes a rigid
monoidal category. In the sequel this will also motivate our antipode
axioms in Section 7.
Unless noted differently, throughout this Section we
suppose $\A$ to be monoidal.

Let us recall from \cor{2.8'} that the unit
representation $\pi_\e$ of a monoidal
weak bialgebra $\A$ may also be realized on $\ALR\equiv\eLR(\A)$,
considered as an $\A$-module via
\beq \lb{7.1}
\pLR (a) b := \eLR (ab),\quad a\in \A,\,b\in\ALR.
\eeq
The equivalence $\pLR\cong\pi_\e$ follows from \no{pLR} and
\no{pLR'}.
In this way, in this Section we
identify $\E\equiv\ALR$.
Moreover, by \cor{3.1} we may identify $\ARL\equiv\hat\E$ as
the dual of $\E$,
with nondegenerate pairing given by
\beq\lb{7.2}
\hat\E\o\E \ni (a\o b) \mapsto \e (ab) \in K\,,
\eeq
Also recall that in the monoidal case $\Assp=\Nssp (\A)$.
%
%
With the present identifications the morphisms
$L_V :\E \times V\to V$ and $R_V:V\to V\times\E$ in $\Rep\A$ introduced in
\no{2.5} and \no{2.6} with inverses \no{2.11} and
\no{2.12} now take the form
\bea\lb{7.3a}
L_V(v) =\eLR (\onne\1)\o\onne\2 \cdot v
&\quad,\quad& R_V(v)=\onne\1\cdot v\o \eLR (\onne\2)
\\
\lb{7.3b}
 L_V^{-1}(a\o v) = \eLR (a) \cdot v\equiv a\cdot v
&\quad,\quad& R_V^{-1} (v\otimes a) =\eRR(a) \cdot v
\eea
where $v\in V$ and $a\in\ALR$.
After these identifications we are prepared to give the

\Definition {7.1}
{i) A {\em pre-rigidity structure}
$(S,\AA,\BB)$ on a monoidal weak bialgebra $\A$ consists of
an anti-algebra map $S:\A\to\A$ and elements $\AA\in\A\o\E$ and
$\BB\in\hat \E\o\A$
satisfying for all $a\in\A$
\beq\lb{7.4}
S(a\1) A^i a\2 \o e^i = A^i \o \eLR (ae^i)
\eeq
\beq\lb{7.5}
\hat e^i \o a\1 B^i S(a\2) = \eRL(\hat e^i a) \o B^i
\eeq
where $\AA\equiv A^i \o e^i$ and $\BB\equiv\hat e^i\o B^i$
and where summations over $i$
are understood.
\\
ii) A pre-rigidity structure $(S,\AA,\BB)$ is called a {\em
rigidity structure} if in addition the elements
$\alpha:= (id\o \e)(\AA)$ and
$\beta := (\e\o id)(\BB)$ satisfy
\beq\lb{7.6}
\onne\1 \beta S(\onne\2)\alpha \onne\3 =\onne
\eeq
\beq\lb{7.7}
S(\onne\1)\alpha\onne\2\beta S(\onne\3) =S(\onne)
\eeq
}
We also call $(\A,\onne,\Del, \e,S,\AA,\BB)$
a {\em rigid weak bialgebra}.
We point out that these axioms are somewhat reminiscent of --
and also motivated by --
Drinfel`d's antipode axioms for quasi-Hopf algebras [Dr]. Also, one
should maybe call this a left rigidity structure,
and one may similarly define
a right rigidity structure on $\A$ as a left rigidity structure on $\Aop$.
Note that if $(S,\AA,\BB)$ is a rigidity structure on $\A$
then $(S,\BB_{op},\AA_{op})$ is a rigidity structure on
$\A_{op}^{cop}$.
%

Given a pre-rigidity structure $(S,\AA,\BB)$ on $\A$
one obtains on $\Rep\A$ a contravariant conjugation
functor $V\to \bar V$ as follows.
Let $V$ be a left $\A$-module with
dual right $\A$-module $\hat V$ and define
$\bar V:= \hat V\cdot S(\onne)$.
Then $\bar V$ becomes a left $\A$-module via
$$
a \cdot u := u\cdot S(a),\quad u\in\bar V,\,a\in \A.
$$
In this way the assignment $V\to \bar V$ provides
a contravariant conjugation functor in $\Rep\A$,
where for $\A$-linear morphisms
$T: V\to W$ we put
$\bar T := T^t|_{\bar W} :\bar W\to \bar V,\ T^t$
being the transpose of $T$.
The terminology ``conjugation" is justified by the following
Lemma, where for left $\A$-modules $V$
we also use the notation
$\piV(a) v \equiv a\cdot v,\ a\in\A,\,v\in V$, where
$\piV:\A \to \End_K V$ denotes the representation homomorphism.

\Lemma {7.2}
{A pre-rigidity structure $(S,\AA,\BB)$ on $\A$
provides in $\Rep\A$ a family of $\A$-linear morphisms
$A_V:\bar V\times V\to \E$ and $B_V :\E\to V\times \bar V$
given by
\beq\lb{7.8}
A_V(u\o v):=\bra u\mid A^i\cdot v\ket e^i
\equiv  \bra u\mid S(\onne\1)\alpha
\onne\2 \cdot v\ket \onne\3\in\ALR\equiv\E
\eeq
\beq\lb{7.9}
B_V(a) := \e (\hat e^j a) \piV(B^j) \equiv
\piV(a\1\beta S(a\2)),\ a\in\ALR\equiv\E,
\eeq
where $\alpha :=(id\o \e) (\AA),\ \beta :=(\e \o id)(\BB)$,
and where in \no{7.9} we have identified
$V\o \hat V\cong \EndK V$.
}
\proof
The second identity in \no{7.9} follows from \no{7.5} and
$\e=\e \c\eRL$ and the second identity in \no{7.8} follows from
\lem{7.2.1} below. The fact that $A_V$ and $B_V$ are
$\A$-linear follows immediately from \no{7.4}, \no{7.5} and the
identities $\e(abc) = \e(\eRL (ab)c) = \e(a\eLR (bc))$ for all
$a,b,c\in\A$, see \lem{2.7}.
\qed

\bsn
In order that the family of morphisms
$A_V,B_V$ indeed provides a rigidity structure
on $\Rep \A$ we also need the axioms \no{7.6} and \no{7.7}.

\Proposition{7.3}
{A pre-rigidity structure  $(S,\AA,\BB)$ on $\A$ is rigid, if
and only if,
under the setting of \lem{7.2}, we have for all $V$ in $\Rep\A$
the {\bf rigidity identities}
\bea \lb{7.10}
R_V^{-1} \c(\onne_V\times A_V)\c(B_V\times \onneV) \c L_V &=&\onneV
\\
\lb{7.11}
L_{\bar V} ^{-1} \c(A_V\times \onne_{\bar V})\c
(\onne_{\bar V}\times B_V)\c R_{\bar V} &=&\onne_{\bar V}\ .
\eea
}
\proof
Using \no{7.9}, \no {7.3a} and $\eR\c\eLR=\eR$ and identifying
$V\o \hat V \cong\End_K V$ we have for
$v\in V$
$$
(B_V\times \onne_V)(L_V(v)) = \piV(\onne\1 \beta S(\onne\2))\o\onne\3 \cdot v.
$$
Similarly, for $v,w\in V$ and $u\in\bar V$ we get
$$
[R_V^{-1} \c(\onne_V\o A_V)] (v\o u\o w)=
\e(\onne\2 e^j) \bra u\mid A^j\cdot w\ket \onne\1 \cdot v
= \bra u\mid S(\onne\2)\alpha \onne\3\cdot w\ket \onne\1\cdot v
$$
by \no{7.3b} and \no{7.4}. Hence
\beanon
[R_V^{-1} \c(\onneV\times A_V)\c(B_V\times\onneV)\c L_V](v)
&=& \onne\1 (\onne_{(1')} \beta S(\onne_{(2')}))
(S(\onne\2)\alpha\onne\3)\onne_{(3')}\cdot v
\\
& =& \onne\1\beta S(\onne\2)\alpha \onne\3 \cdot v =v.
\eeanon
by \no{7.6}. Using \no{7.7} the identity \no{7.11} follows similarly. Putting $V=\A$
we see that the axioms \no{7.6} and \no{7.7} are also necessary.
\qed

\bsn
We leave it to the reader to check that for $\A$-linear
morphisms $T:V\to W$ and $\bar T\equiv T^t|_{\bar W}:\bar W\to\bar V$
the definitions \no{7.8} and \no{7.9} imply
\beanon
\bar T &=& L_{\bar V}^{-1} \c(A_V\times \onne_{\bar V})
 \c(\onne_{\bar V} \times T\times\onne_{\bar V}) \c
(\onne_{\bar V} \times B_V)\c R_{\bar V}\ .
\\
T &=&R_V^{-1} \c(\onne_V\times A_V)\c
(\onne_V\times \bar T\times \onne_V)\c(B_V\times \onne_V)
\c L_V
\eeanon
expressing the standard isomorphism
$\HomA(V,W)\cong\HomA(\bar W,\bar V)$
in rigid monoidal categories.

Next, we point out that for any (pre)rigidity structure
$(S,\AA,\BB)$ on $\A$ the elements
$\AA$ and $\BB$ are already uniquely
determined by $\alpha\equiv(\id\o \e)(\AA)$ and
$\beta \equiv(\e \o \id) (\BB)$.
\Lemma{7.2.1}
{Let $(S,\AA,\BB)$ be a pre-rigidity structure on $\A$. Then
\bea
\lb{7.7c}
\AA&=& S(\onne\1)\alpha \onne\2 \o \onne\3 \\
\lb{7.7d}
\BB&=& \onne\1\o \onne\2\beta S(\onne\3)
\eea
}
\proof
Using $\e=\e\c\eLR =\e\c\eRL$ we compute
\beanon
\AA&=&(id \o \e_{LR})(\AA)= A^i \o \onne\2 \e(\onne\1 e^i)
= S(\onne\1)\alpha\onne\2\o \onne\3
\\
\BB&=& (\eRL \o id)(\BB)
= \onne\1 \e(\hat e^i \onne\2)\o B^i = \onne\1 \o \onne\2\beta S(\onne\3)
\ \qed
\eeanon
\lem{7.2.1} implies that the axioms for rigidity structures
$(S,\AA,\BB)$ may equivalently be reformulated in terms of the data
$(S,\al,\be)$.
To this end, for $x\in\A$ and $S:\A\to\A$
an anti-algebra map let us introduce
the left and right ``$S$-adjoint" actions
\beq\lb{7.7e}
x\liS a := S(a\1)x a\2\quad , \quad a\reS x:= a\1 x S(a\2)
\eeq

\Proposition{7.2.2}
{Let $S:\A\to\A$ be an anti-algebra map and let $\alpha,\beta \in\A$ satisfy
$\alpha\liS  \onne =\alpha$ and $\onne \reS \beta =\beta$.
Put $\AA:=\alpha\liS  \onne\1 \o \onne\2$ and
$\BB:=\onne\1\o\onne\2 \reS\beta$ as in \no{7.7c}
and \no{7.7d}. Then $(S,\AA,\BB)$
provides a pre-rigidity structure on $\A$ if and
only if for all $a\in\A$
\beq\lb{7.7f}
\alpha\liS a=\alpha \liS  \eRL(a)~~~,~~~ a \reS\beta =\eLR(a) \reS \beta
\eeq
If in addition $\al$ or $\be$ are invertible, then
$S(\onne)=\onne$.
}
\proof
By \defi{7.1}, $(S,\AA,\BB)$
provides a pre-rigidity structure if and only
if for all $a\in\A$
\bea
\lb{7.7g}
\alpha \liS (\onne\1 a)\o\onne\2
&=& \alpha\liS \onne\1\o \eLR(a\onne\2)
\\
\lb{7.7h}
\onne\1\o a\onne\2 \reS \beta &=& \eRL(\onne\1 a)\o\onne\2 \reS\beta\ .
\eea
Applying $id\o\e$ and $\e\o id$, respectively, yields
\no{7.7f}. Conversely, using the identity
\beq\lb{7.7i}
\eRL(\onne\1 a)\o\onne\2 = \onne_{(1')}\e(\onne\1 a\onne_{(2')})\o\onne\2
=\onne\1\o\eLR(a\onne\2),
\eeq
\Eq{7.7f} implies \no{7.7g} and \no{7.7h}.
Finally, we have $\al=S(\onne\1)\al\onne\2=S(\onne)\al$ and
similarly $\be=\be S(\onne)$. Hence, if $\al$ or $\be$ are
invertible, then $S(\onne)=\onne$.
\qed

\bsn
Note that the normalization conditions
$\alpha\liS  \onne =\alpha$ and $\onne \reS \beta =\beta$ in
\prop{7.2.2} are imposed to reproduce the original identities
$\al=(\id\o\e)(\AA)$ and $\be=(\e\o\id)(\BB)$.
In view of \lem{7.2.1} and \prop{7.2.2} we will from now
on equivalently talk of (pre)rigidity structures  on $\A$
given by the data $(S,\alpha,\beta)$.

Next, we recall from [Dr] that there is a natural notion of
twist equivalence for rigidity structures $\SSS$.
Let $u, \bar u\in\A$ satisfy
\bleq{u}
\bar uu=S(\onne)\quad,\quad u\bar uu=u\quad,\quad\bar uu\bar u=\bar u,
\eeq
and put
\bleq{twist}
S'(a):=uS(a)\bar u,\quad,\quad
\al':=u\al\quad,\quad\be':=\be\bar u\,.
\eeq
Then \prop{7.2.2} assures that $(S',\,\al',\,\be')$ again
provides a rigidity structure. The inverse transformation is
given by interchanging $u$ and $\bar u$.
One also checks that this indeed provides an equivalence relation.
In \prop{B1} of
Appendix B we will show that any two rigidity structures on
a monidal weak bialgebra are twist equivalent in this sense.

In ordinary bialgebras $S:\A\to\A$ is an
antipode if and only if $(S,\,\al\equiv\onne,\,\be\equiv\onne)$
is a rigidity structure. Thus, to approach
and motivate our antipode axioms in Sect. 7, we now study
rigidity structures satisfying $\al=\be=\onne$.

\Definition{7.10'}
{A (pre-)rigidity structure $\SSS$ is called {\em
normalizable}, if $\al=\be^{-1}$, and it is called {\em normal},
if $\al=\be=\onne$. In this case $S:\A\to\A$ is called a
{\em normal rigidity map}.
}
In Example 2 of Appendix D we will construct a rigid weak bialgebra
whose rigidity structure is {\em not} normalizable.

Clearly, a rigidity structure is normalizable, if and only if
it can be twisted into a normal one.
Thus, by \prop{B1} a normal rigidity map $S$
on a monoidal weak bialgebra is uniquely determined, provided
it exists.
We refrain from calling such an $S$ an antipode, since in
general the dual $\hat S\equiv S^t:\hA\to\hA$ will not be of
the same type.
Our antipode axioms in Sect. 7 will
be symmetric under the duality flip $S\leftrightarrow\hat S$.
We will see in Sect. 8 (\cor{9.3}),
that on {\em bimonoidal} weak bialgebras $S$ is a
normal rigidity map if and only if it is an antipode.

To approach these results let us now introduce, following
[BSz], the linear maps $\PLR_S:\A\to\A$ given by
\bleq{7.30}
\PL_S(a):=a\reS\onne\equiv a\1S(a\2)\quad,
\quad\PR_S(a):=\onne\liS a\equiv S(a\1)a\2\ .
\eeq
Then \prop{7.2.2} immediately implies
\Corollary{7.11}
{An algebra antimorphism $S:\A\to\A$ on a monoidal weak
bialgebra $\A$ is a normal pre-rigidity map if and only if
\beq\lb{7.31}
\ba{rclcrcl}
\PL_S\circ\eLR &=&\PL_S&\quad , \quad&
\PR_S\circ\eRL&=&\PR_S
\\
\PL_S(\onne) &=&\onne &,& \PR_S(\onne) &=&\onne\ .
\ea\eeq
}
Let us next observe that the
antipode axioms of [BSz] would imply (see Appendix A)
\bleq{7.31'}
\PL_S=\eLR\quad,\quad \PR_S=\eRL
\eeq
from which the identities \no{7.31} would follow by \lem{2.7}.
We now show, that converseley \no{7.31} implies \no{7.31'} if
and only if $\AsL=\AsR$, for $\s=L$ and $\s=R$.
To this end let us introduce the linear spaces
$$
\A^{L/R}:=\PLR(\A)\subset\A\,,
$$
where from now on we simplify our notation by writing
$\PLR\equiv\PLR_S$.
Then $\A^{L/R}$ naturally becomes a left (right) $\A$-module
under the left (right) $S$-adjoint action
\bleq{7.34}
\ba{ccc}
\reS:\A\o\A^L\to\A^L &\quad ,\quad &\liS:\A^R\o\A\to\A^R
\\
a\reS x:=a\1\,x\,S(a\2) &,& y\liS a:=S(a\1)\,y\,a\2.
\ea\eeq
It turns out that these $\A$-modules are isomorphic to $\AsR$
and $\AsL$, respectively, with left (right) $\A$ actions
\footnote
{By \lem{2.7} and \cor{3.1} the right $\A$-modules $\AsL$ are dual
to the left $\A$-modules $\AsR$, where the nondegenerate
pairing is given by \no{2.7.1}.
}
\bea\lb{7.35}
\psR(a)x &:=&\esR(ax),\quad a\in\A,\,x\in\AsR.
\\
\lb{7.36}
y\psL(a)&:=&\esL(ya),\quad a\in\A,\,y\in\AsL\,.
\eea

\Theorem{7.12}
{Let $\A$ be monoidal and let $S$ be a normal pre-rigidity map
on $\A$. Then for $\s\in\{L,R\}$
\bea\lb{7.37}
\A^L=\ALL &\qquad &\A^R=\ARR
\\\lb{7.38}
\PL\circ \esR=\PL &&\PR\circ\esL=\PR
\\\lb{7.39}
\PL\circ \eLL=\eLL &&\PR\circ\eRR=\eRR
\\\lb{7.40}
\PL\circ \eRL=S\circ\eRL &&\PR\circ\eLR=S\circ\eLR
\\\lb{7.41}
\eR\circ\PL=\eR &&\eL\circ\PR=\eL\,,
\eea
and we have the following commuting diagrams
of left (right) $\A$-module isomorphisms.
\beq\lb{diag}
\ba{lr}

\setlength{\unitlength}{0.0007in}%
\begingroup\makeatletter\ifx\SetFigFont\undefined%
\gdef\SetFigFont#1#2#3#4#5{%
  \reset@font\fontsize{#1}{#2pt}%
  \fontfamily{#3}\fontseries{#4}\fontshape{#5}%
  \selectfont}%
\fi\endgroup%
\begin{picture}(2955,2310)(1400,-2900)

\thicklines

\put(1201,-2686){\vector( 1, 0){2400}}

\put(3601,-2836){\vector(-1, 0){2400}}

\put(2314,-1027){\vector(-3,-4){1188}}

\put(1033,-2522){\vector( 3, 4){1188}}

\put(3676,-2611){\vector(-3, 4){1188}}

\put(2581,-938){\vector( 3,-4){1188}}

\put(730,-2836){\makebox(0,0)[lb]{$\ALL$}}

\put(3730,-2836){\makebox(0,0)[lb]
{\smash{\SetFigFont{12}{14.4}{\rmdefault}{\mddefault}{\updefault}$\ARR$}}}

\put(2210,-840){\makebox(0,0)[lb]
{\smash{\SetFigFont{12}{14.4}{\rmdefault}{\mddefault}{\updefault}$\ALR$}}}

\put(2326,-3090){\makebox(0,0)[lb]
{\smash{\SetFigFont{12}{14.4}{\rmdefault}{\mddefault}{\updefault}$\PL$}}}

\put(1801,-1936){\makebox(0,0)[lb]
{\smash{\SetFigFont{12}{14.4}{\rmdefault}{\mddefault}{\updefault}$\PL$}}}

\put(1190,-1711){\makebox(0,0)[lb]
{\smash{\SetFigFont{12}{14.4}{\rmdefault}{\mddefault}{\updefault}$\eLR$}}}

\put(2720,-1936){\makebox(0,0)[lb]
{\smash{\SetFigFont{12}{14.4}{\rmdefault}{\mddefault}{\updefault}$\eLR$}}}

\put(3230,-1711){\makebox(0,0)[lb]
 {\smash{\SetFigFont{12}{14.4}{\rmdefault}{\mddefault}{\updefault}$\eRR$}}}

\put(2251,-2570){\makebox(0,0)[lb]
{\smash{\SetFigFont{12}{14.4}{\rmdefault}{\mddefault}{\updefault}$\eRR$}}}

\end{picture}

&

\setlength{\unitlength}{0.0007in}%
\begingroup\makeatletter\ifx\SetFigFont\undefined%
\gdef\SetFigFont#1#2#3#4#5{%
  \reset@font\fontsize{#1}{#2pt}%
  \fontfamily{#3}\fontseries{#4}\fontshape{#5}%
  \selectfont}%
\fi\endgroup%
\begin{picture}(2955,2310)(700,-2900)

\thicklines

\put(1201,-2686){\vector( 1, 0){2400}}

\put(3601,-2836){\vector(-1, 0){2400}}

\put(2314,-1027){\vector(-3,-4){1188}}

\put(1033,-2522){\vector( 3, 4){1188}}

\put(3676,-2611){\vector(-3, 4){1188}}

\put(2581,-938){\vector( 3,-4){1188}}

\put(730,-2836){\makebox(0,0)[lb]{$\ARR$}}

\put(3730,-2836){\makebox(0,0)[lb]
{\smash{\SetFigFont{12}{14.4}{\rmdefault}{\mddefault}{\updefault}$\ALL$}}}

\put(2210,-840){\makebox(0,0)[lb]
{\smash{\SetFigFont{12}{14.4}{\rmdefault}{\mddefault}{\updefault}$\ARL$}}}

\put(2326,-3090){\makebox(0,0)[lb]
{\smash{\SetFigFont{12}{14.4}{\rmdefault}{\mddefault}{\updefault}$\PR$}}}

\put(1801,-1936){\makebox(0,0)[lb]
{\smash{\SetFigFont{12}{14.4}{\rmdefault}{\mddefault}{\updefault}$\PR$}}}

\put(1190,-1711){\makebox(0,0)[lb]
{\smash{\SetFigFont{12}{14.4}{\rmdefault}{\mddefault}{\updefault}$\eRL$}}}

\put(2720,-1936){\makebox(0,0)[lb]
{\smash{\SetFigFont{12}{14.4}{\rmdefault}{\mddefault}{\updefault}$\eRL$}}}

\put(3230,-1711){\makebox(0,0)[lb]
 {\smash{\SetFigFont{12}{14.4}{\rmdefault}{\mddefault}{\updefault}$\eLL$}}}

\put(2251,-2570){\makebox(0,0)[lb]
{\smash{\SetFigFont{12}{14.4}{\rmdefault}{\mddefault}{\updefault}$\eLL$}}}

\end{picture}
\\
\ &\ 
\ea
\eeq
%

%
%
%
Here, in the left diagram we consider $\ALR,\,\ARR$ and
$\ALL\equiv\A^L$ as left $\A$-modules with $\A$-actions
\no{7.35} and (\ref{7.34}Left), respectively, and in the right
diagram we consider $\ARL,\,\ALL$ and
$\ARR\equiv\A^R$ as right $\A$-modules with $\A$-actions
\no{7.36} and (\ref{7.34}Right), respectively.
}
\proof
By passing to $\A_{op}^{cop}$ it is enough to prove the left
statements.
Eq. \no{7.38} follows from \no{7.31} and the identities
$\eLR=\eLR\,\eRR$ and $\eRL=\eRL\,\eLL$, see \lem{2.7}.
Eqs. \no{7.39} and \no{7.40} follow from $\Assp=\Nssp(\A)$.
In particular, this also gives $\ALL\subset\A^L$.
Together with $\A^L=\PL(\ARR)$ by \no{7.38} this implies
\bleq{7.44}
\dim\ALL\le\dim\A^L\le\dim\ARR
\eeq
and hence equality by \no{2.7.6}, thus proving (\ref{7.37}Left).
Let us now turn to the left diagram in \no{diag}.
First, by \cor{2.8'} $\eRR:\ALR\to\ARR$ is an $\A$-linear
bijection with inverse $\eLR:\ARR\to\ALR$.
Second, by \no{7.38} and \no{7.35}
$\PL:\AsR\to\A^L\equiv\ALL$ is $\A$-linear and surjective,
whence bijective by \no{2.7.6}.
Third, by \no{7.38} and \no{7.39}
$$
\ba{rcccl}
\PL\circ\eRR|_{\ALL} &=&\PL|_{\ALL}&=&\id_{\!\ALL}
\\
\PL\circ\eLR|_{\ALL} &=&\PL|_{\ALL}&=&\id_{\!\ALL}
\ea
$$
and therefore
\bleq{7.45}
(\PL|_{\AsR})^{-1} = \esR|_{\ALL}\,.
\eeq
Finally, the diagram commutes, since $\eLR\,\eRR=\eLR$ by
\lem{2.7}.
We are left to prove (\ref{7.41}Left), which follows since
$\eR=\eR\,\eRR$ by \lem{2.7}, and since \no{7.45} implies
$\eRR=\eRR\circ\PL\circ\eRR$, whence
$$
\eR\circ\PL = \eR\circ\eRR\circ\PL\circ\eRR = \eR\circ\eRR =\eR
$$
by \no{7.38}.
\qed

\Corollary{7.13}
{Under the setting of \thm{7.12} we have
\bea\lb{7.42}
\PL=\eLR &\Llra & \ALL=\ALR
\\\lb{7.43}
\PR=\eRL &\Llra & \ARR=\ARL\,.
\eea
}
\proof
By passing to $\A_{op}^{cop}$ it suffices to prove the first
statement.
If $\PL=\eLR$ then $\ALL\equiv\PL(\A)=\ALR$. Converseley, if
$\ALL=\ALR$, then $\eLR|_{\ALL}=\id$ and (\ref{diag}Left) implies
$\PL|_{\ARR}=\eLR|_{\ARR}$.
Hence, by \lem{2.7},
$$
\PL\equiv\PL\circ\eRR=\eLR\circ\eRR=\eLR\,.
\quad\qed
$$
Note that the conditions of \cor{7.13} are in particular satisfied if
$\A$ is bimonoidal, yielding $\AL=\ALL=\ALR$ and
$\AR=\ARR=\ARL$.
In the next Section we will take the left hand side of Eqs.
\no{7.42} and \no{7.43} as the defining relations for a
pre-antipode S.



\sec{The Antipode Axioms}

Let us first understand why for general weak
bialgebras the ordinary antipode axioms
\bleq{8.1}
S(a\1)a\2=a\1S(a\2)=\e(a)\onne,\quad a\in\A
\eeq
would be too restrictive.
Call a linear map $S:\A\to\A$ satisfying \no{8.1} a {\em Hopf antipode}.
As for ordinary Hopf algebras, a Hopf antipode $S$ would be the
inverse of $\idA$ in the convolution algebra $(\EndK\A,*)$,
where $(S*T)(a):=S(a\1)T(a\2)$. Also, if
$S$ is a Hopf antipode on $\A$, then $\hat S:=S^t$ is a Hopf
antipode on $\hA$.

\Lem{8.1}
{Let $\AAA$ be a weak bialgebra with Hopf antipode $S$. Then
\\
i) $\eLR(a)=\eRL(a)=\e(a)\onne,\ \forall a\in\A$.
\\
ii)
$\A$ is right-monoidal iff $\e$ is multiplicative,
in which case $\A$ is also monoidal.
\\
iii)
$\A$ is right-comonoidal iff $\Del(\onne)=\onne\o\onne$, in which case
$\A$ is also comonoidal.
}
\proof
The unit  in $(\EndK\A,*)$ is given by $a\mapsto\e(a)\onne$.
Hence, if $\idA$ has a convolution inverse $S$, part (i)
follows from the identities
\bleq{8.2}
\eLR*\idA=\idA=\idA*\eRL
\eeq
which one gets by putting $b=\onne$ in \no{2.16}.
Part (ii) follows by applying $\e$ to (\ref{2.19}right)
to get $\e(ab)=\e(a)\e(b)$, and
part (iii) follows by duality.
\qed

\bsn
\lem{8.1} shows, that in general the Hopf antipode axioms
\no{8.1} are too restrictive. Instead,
motivated by our analysis of rigidity structures in Sect. 6 and
in particular by \cor{7.13} we now define

\Definition{8.2}
{A {\em pre-antipode} $S$ on a weak bialgebra $\A$ is a
linear map $S:\A\to\A$ satisfying for all $a\in\A$
\bleq{8.4}
a\1S(a\2)=\eLR(a)\quad ,\quad S(a\1)a\2=\eRL(a)
\eeq
A pre-antipode $S$ is called an {\em antipode}, if
\bleq{8.5}
S(a\1)a\2S(a\3)=S(a),\quad \forall a\in \A
\eeq
}

\bsn
Note that by \no{8.2} a pre-antipode always satisfies
\bleq{8.5'}
a\1S(a\2)a\3=a,\quad\forall a\in\A\ .
\eeq
In $(\EndK\A,*)$ the identities \no{8.4},
\no{8.5} and \no{8.5'}
 can be rewritten, respectively, as
\bea\lb{8.6}
\idA*S=\eLR &\quad , \quad& S*\idA=\eRL
\\\lb{8.7}
S*\idA*S = S &,& \idA*S*\idA=\idA\,.
\eea
Hence,
in $(\EndK\A,*)$ an antipode $S$ may be considered as a
``quasi-inverse" of $\idA$.
Also note that if $S$ is a (pre-)antipode on $\A$, then it is
also a (pre-) antipode on $\A_{op}^{cop}$ and by
\no{2.17} its transpose $\hat S$ is a (pre-)antipode on $\hA$.
By \lem{8.1}i) a Hopf antipode is always an antipode, and a
pre-antipode is a Hopf antipode iff
$\eLR(a)=\eRL(a)=\e(a)\onne,\ \forall a\in\A$.
Moreover, we have

\Lem{8.3}
{i) A weak bialgebra $\A$ has at most one antipode $S$.
If $\A$ has a preantipode $S_p$ then $\eLR*\eLR=\eLR,\ \eRL*\eRL=\eRL$
and $S:=S_p*\idA*S_p$ provides an antipode.
\\
ii) If a
pre-antipode $S$ is anti-multiplicative, then $S(\onne)=\onne$,
and if it is anti-comultiplicative, then $\e\circ S=\e$.
}
\proof
(i) If $S_1$ and $S_2$ are antipodes, then
$
S_1=S_1*\idA*S_1=S_1*\idA*S_2=S_2*\idA*S_2=S_2\ .
$
If $S_p$ is a preantipode, then by \no{8.2}
$
\eLR*\eLR=\eLR*\idA*S_p=\idA*S_p=\eLR
$
and similarly $\eRL*\eRL=\eRL$. Hence $S:=S_p*\idA*S_p$ is an
antipode.
(ii)
If $S$ is a pre-antipode satisfying
$S(ab)=S(b)S(a)$ then, using $\eLR(\onne)=\onne$,
$S(\onne)=S(\onne)S(\onne_1)\onne\2=S(\onne_1)\onne\2=\onne$.
The statement for anti-comultiplicative $S$ follows by duality.
\qed

\Lemma{8.3'}
{A pre-antipode $S$ on a right-monoidal or right-comonoidal
weak bialgebra $\A$ is an antipode, if $S$ is
anti-multiplicative or anti-comultiplicative.
}
\proof
By duality it is enough to consider the case of $S$ being
anti-multiplicative. If $\A$ is right-comonoidal, then
$a\1\o\eLR(a\2)=\onne\1a\o\onne\2$ by the dul of
(\ref{2.23}right). Hence,
$$
S(a\1)a\2S(a\3)=S(a\1)\eLR(a\2)=S(a)S(\onne\1)\onne\2=S(a)\ .
$$
If instead $\A$ is right-monoidal, then for all $a,b\in\A$
\beanon
\eLR(a)S(b) &=&a\1S(ba\2)=a\1\e(b\1a\2)S(b\2a\3)
\\
&=&\eRL(b\1)a\1S(b\2a\2)
\\
&=& S(b\1)b\2a\1S(b\3a\2)\,,
\eeanon
where in the second line we have used (\ref{2.20}right).
Putting $a=\onne$ we conclude that $S$ is an antipode.
\qed

\Corollary{8.3'}
{An algebra antimorphism $S$ on a monoidal weak bialgebra $\A$
is a pre-antipode (and therefore an antipode) if and only if
$S$ is a normal pre-rigidity map and $\AsL=\AsR$, for $\s=L$
and $\s=R$.
}
\proof
This follows from \lem{8.3'}
and \cor{7.13}, since \no{8.4} is the same as
\no{7.31'}, and therefore implies \no{7.31}.
\qed

\bsn
In ordinary bialgebras an antipode is always a bialgebra
antimorphism and hence a normal rigidity map. In weak bialgebras,
the following Theorem analyses necessary and sufficient
conditions for an antipode $S$ to be anti-multiplicative and/or
anti-comultiplicative.

\Theorem{8.4}
{Let $\A$ be a weak bialgebra with pre-antipode $S$ and
consider the following additional properties:
\begin{center}
\begin{tabular}{rlrl}
1a)& $S$ is anti-multiplicative.
\qquad&
1b)& $\A$ is right-monoidal.
\\
1c)& $[\ALR,\ARL]=0$.
\qquad &
1d)& $S$ is an antipode.
\end{tabular}
\end{center}
Then the following implications hold
$$
\ba{rrcl}
1i)&\qquad 1a)+1b)&\Longrightarrow &1c)+1d)
\\
1ii)& 1a)+1c)&\Longrightarrow& 1b)+1d)
\\
1iii)&\qquad 1b)+1c)+1d)&\Longrightarrow &1a)
\ea
$$
Similarly, consider the following properties
\begin{center}
\begin{tabular}{rlrl}
2a)& $S$ is anti-comultiplicative.
\qquad&
2b)& $\A$ is right-comonoidal.
\\
2c)& $[\hALR,\hARL]=0$.
\qquad&
2d)& $S$ is an antipode.
\end{tabular}
\end{center}
Then the following implications hold
$$
\ba{rrcl}
2i)& 2a)+2b)&\Longrightarrow& 2c)+2d)
\\
2ii)& 2a)+2c)&\Longrightarrow &2b)+2d)
\\
2iii)&\qquad 2b)+2c)+2d)&\Longrightarrow &2a)
\ea
$$
}
\proof
Part 2.) is the dual of part 1.).
Le us first prove 1iii). If $\A$ is right-monoidal, then by
(\ref{2.19}right)
\beanon
a\1b\1S(b\2)S(a\2) &=&a\1\eLR(b)S(a\2)=\e(a\1b)a\2S(a\3)
\\
&=&\eLR(\e(a\1b)a\2)= \eLR(a\eLR(b))
\\
&=&\eLR(ab)\equiv a\1b\1S(a\2b\2)\,,
\eeanon
where in the last line we have used \lem{2.7}{\sc r}. The same
argument in $\A_{op}^{cop}$ gives
$$
S(a\1b\1)a\2b\2=S(b\1)S(a\1)a\2b\2\ .
$$
Hence, using $[\ALR,\ARL]=0$ and \no{8.5}
\beanon
S(ab) &=& S(a\1b\1)a\2b\2S(a\3b\3)
\\
&=&S(b\1)S(a\1)a\2b\2S(b\3)S(a\3)
\\
&=& S(b\1)b\2S(b\3)S(a\1)a\2S(a\3)
\\
&=&S(b)S(a)\,,
\eeanon
proving part 1iii). To prove 1.ii) let
 $[\ALR,\ARL]=0$ and $S$ anti-multiplicative, then by \no{8.5'}
\beanon
a\eLR(b) &=& a\1S(a\2)a\3b\1S(b\2)
\\
&=& a\1b\1S(b\2)S(a\2)a\3
\\
&=&a\1b\1S(a\2b\2)a\3
\\
&=&\eLR(a\1b)a\2\ .
\eeanon
Hence, by (\ref{2.19}right), $\A$ is right-monoidal and by
\lem{8.3'} $S$ is an antipode.
Finally, to prove 1.i) assume $S$ anti-multiplicative and $\A$
right-monoidal. Then $S$ is an antipode by \lem{8.3'}, and
Eq. (\ref{2.19}right) implies for all $a,b\in\A$
\bleq{8.9}
a\eLR(b)=a\1b\1S(a\2b\2)a\3=a\1\eLR(b)\eRL(a\2)\,.
\eeq
Hence,
\beanon
\eRL(a)\eLR(b) &=& S(a\1)a\2\eLR(b)\eRL(a\3)
\\
&=&\eRL(a\1)b\1S(a\2b\2)a\3
\\
&=&b\1\eRL(a\1b\2)S(a\2b\3)a\3
\\
&=&b\1S(a\1b\2)a\2
\\
&=&\eLR(b)\eRL(a)\,,
\eeanon
where we have used \no{8.4} and \no{8.9} in the first line,
\no{8.4} in the second line,
(\ref{2.20}right) in the third line and the antipode identity
$S(a)=\eRL(a\1)S(a\2)$ \no{8.5} in the fourth line. Hence,
$\ALR$ and $\ARL$ commute.
\qed

\Cor{8.5}
{Let $\A$ be a monoidal weak bialgebra with antipode $S$.
Then $[\ALR,\ARL]=0$ if and only if $\ALL=\ALR$ and $\ARR=\ARL$,
and in this case $S$ is a normal rigidity map.
}
\proof
By part 3.) of \prop{2.6'} $[\ALR,\ARL]=0$ follows from
$\ALL=\ALR$ and $\ARR=\ARL$. Conversley,
if $[\ALR,\ARL]=0$ then $S$ is anti-multiplicative by
 \thm{8.4}(1iii) and \cor{8.5} follows from \cor{8.3'}.
\qed

\bsn
Note that the conditions of \cor{8.5} in particular hold if
$\A$ is bimonoidal.

Next, we provide conditions under which an antipode $S$
is invertible by using an invertibility result for rigidity
maps proven in \thm{7.8} of Appendix B. To this end we
need the counit to be $S$-invariant.

\Lemma{8.6}
{Let $S$ be an antipode on $\A$ and assume the axioms ({\sc r})
of \lem{2.7}, i.e. $\e(ab)=\e(a\onne\2)\e(\onne\1b),\ \forall
a,\,b\in\A$. Then $\e\circ S=\e$ and, more generally,
\bleq{8.11}
\eL\circ S=\eL\circ\eLR \quad,\quad \eR\circ S=\eR\circ\eRL
\eeq
}
\proof
Using \lem{2.7} and \no{8.5} we compute
$$
\eL(S(a))=\eL(\eRL(a\1)S(a\2))=\eL(a\1S(a\2))=\eL(\eLR(a))\,.
$$
The second identity in \no{8.11} follows by passing to
$\A_{op}^{cop}$.
Pairing these equations with $\onne\in\A$ and using
$\e\circ\essp=\e$ we get $\e\circ S=\e$.
\qed

\bsn
Since the condition ({\sc r}) of \lem{2.7} in particular holds if $\A$ is
right-monoidal we arrive at

\Corollary{8.7}
{Let $\A$ be a weak bialgebra with pre-antipode $S$.
If $\A$ is monoidal (comonoidal) and $S$ is anti-multiplicative
(anti-comultiplicative), respectively, then $S$ is bijective.
}
\proof
If $\A$ is monoidal and $S$ is anti-multiplicative, then $S$ is
an antipode by \lem{8.3'}, $\e\circ S=\e$ by \lem{8.6},
$S(\onne)=\onne$ by \lem{8.3} and
$S$ is a normal rigidity map by \cor{8.3'}.
Hence, the invertibility of $S$ follows from
\thm{7.8}.
The ``co"-statement follows by duality.
\qed

\bsn
Finally, we show that under the conditions of \cor{8.7}
$S^{-1}$ is an antipode on $\A_{op}$ and $\A^{cop}$.
As for ordinary Hopf algebras, we call such a map a {\em pode}.

\Definition{9.5}
{A {\em pre-pode} $\bar S$ on a weak bialgebra $\A$ is a
linear map $\bar S:\A\to\A$ satisfying for all $a\in\A$
\bleq{9.1}
a\2\bar S(a\1)=\eRR(a)\quad ,\quad \bar S(a\2)a\1=\eLL(a)\,.
\eeq
A pre-pode $S$ is called a {\em pode}, if
\bleq{9.2}
\bar S(a\3)a\2\bar S(a\1)=\bar S(a),\quad \forall a\in \A\,.
\eeq
}
One immediately checks that the above axioms are precisely the
antipode axioms in $\A_{op}$ and $\A^{cop}$.
For ordinary Hopf algebras the inverse of an antipode is always
a pode. In our setting we have

\Lemma{9.6}
{Let $\A$ be a weak bialgebra with invertible pre-antipode $S$.
\\
1.) Assume $S$ anti-multiplicative. Then $S^{-1}$ is a pre-pode if
and only if
\bleq{9.3a}
\eLR=S\c\eRR \quad , \quad \eRL=S\c\eLL\,.
\eeq
In this case $S$ is an antipode and $S^{-1}$ is a pode.
\\
2.) Assume $S$ anti-comultiplicative. Then $S^{-1}$ is a pre-pode if
and only if
\bleq{9.3b}
\eLR=\eLL\c S \quad , \quad \eRL=\eRR\c S\,.
\eeq
In this case $S$ is an antipode and $S^{-1}$ is a pode.
}
\proof
Part 2.) is the transpose of the dual version of 1.).
To prove 1.) assume $S(ab)=S(b)S(a)$ and apply $S^{-1}$ to
\no{8.4} to obtain
$$
a\2S^{-1}(a\1)=S^{-1}(\eLR(a))\quad ,\quad
S^{-1}(a\2)a\1=S^{-1}(\eRL(a))\,.
$$
Hence, \no{9.3a} is equivalent to $S^{-1}$ being a pre-pode.
In this case $S^{-1}$ is also a pode and $S$ is also an
antipode, since
\beanon
S(a\1)a\2S(a\3)&=&S(\eRR(a\2)a\1)=S(a)
\\
S^{-1}(a\3)a\2S^{-1}(a\1)&=&S^{-1}(a\1\eRL(a\2))=S^{-1}(a)
\eeanon
by the identities \no{2.16} for $b=\onne$.
\qed

\bsn
The conditions of \lem{9.6} in particular hold if $\A$ is
monoidal or comonoidal, respectively.

\Proposition{9.7}
{Let $\A$ be (co-)monoidal and $S:\A\to\A$ a (co-)algebra
anti-automorphism, respectively.
Then $S$ is a pre-antipode (and therefore an antipode) if and
only if $S^{-1}$ is a pre-pode (and therefore a pode).
}
\proof
If $\A$ is monoidal we have $\Assp=\Nssp(\A)$, and if $S$ is
an anti-multiplicative pre-antipode then it is an antipode by
\lem{8.3'} and $\eRR(a)\in\NRL(\A)$ by \cor{8.3'}. Hence, by
\lem{2.7}{\sc r}),
$$
\eLR(a)=\eLR(\eRR(a))=\eRR(a)\1S(\eRR(a)\2)=
\onne\1S(\eRR(a)\onne\2)=S(\eRR(a))\,.
$$
The same argument in $\A_{op}^{cop}$ yields $\eRL=S\c\eLL$,
whence $S^{-1}$ is a pode by \lem{9.6}. Repeating these arguments in
$\A_{op}$ yields the inverse implication.
Finally, the ``co"-statements follow by duality.
\qed


\sec{Weak Hopf Algebras}

In order to obtain an explicitely selfdual notion of weak Hopf
algebras the following Definition will be appropriate.

\Definition{9.1}
{A weak Hopf algebra $\A$ is a bimonoidal weak bialgebra with
antipode $S$.
}
Hence, if $(\A,S)$ is a weak Hopf algebra, then its dual
$(\hA,\hat S)$ is also a weak Hopf algebra.
In this Section we will show that in a weak Hopf algebra $\A$
the antipode $S$ is always a bialgebra anti-automorphism and
that $S^{-1}$ is a pode.
If $S$ is already known to be
anti-(co)multiplicative, then part of the bimonoidality axioms
may be dropped or replaced altogether by the requirement
$[\AL,\AR]=0$.
We also show that on bimonoidal weak bialgebras $S$ is an
antipode if and only if it is a normal rigidity map.
As an application, the generalized Kac algebras of [Ya] will
then be shown to be weak Hopf algebras with an involutive
antipode.
The relation of \defi{9.1} with the axioms of [BSz,Sz] will be
clarified in Appendix A.
Let us now first observe

\Lemma{9.3}
{A weak Hopf algebra $\A$ with antipode $S$ is an ordinary Hopf
algebra, if and only if  one of the following conditions hold.
\\
i) $\e$ is multiplicative
\\
ii) $\Del(\one)=\one\o\one$
\\
iii) $S$ is a Hopf antipode.
}
\proof
Clearly, (i)+(ii) $\Rightarrow$ (iii), and in this case $\A$ is
an ordinary Hopf algebra. Conversely, (iii) $\Rightarrow$ (i)+(ii)
by \lem{8.1} and if $\A$ is bimonoidal then (i) $\Lra$ (ii) by
\no{2.26} and \no{2.28}.
\qed

\Proposition{9.2}
{For a weak bialgebra $\A$ with pre-antipode $S$ the following
equivalencies hold:
\\
i) $\A$ is comonoidal and $S$ is anti-multiplicative
\\
ii) $\A$ is monoidal and $S$ is anti-comultiplicative
\\
iii) $\A$ is bimonoidal and $S$ is an antipode.
}
\proof
If $\A$ is bimonoidal, then
$[\AL,\AR]=0$ and $[\hAL,\hAR]=0$ and in this case, by
\thm{8.4}, an antipode $S$ is a
bialgebra anti-homomorphism. This proves (iii) $\Rightarrow$
(i)+(ii). Assume now (i). Then
comonoidality implies $[\AL,\AR]=0$ and by part 1ii) of
\thm{8.4} $S$ is an antipode and $\A$ is right-monoidal. Thus,
$\A$ is bimonoidal by \cor{3.7}. This proves
(i) $\Rightarrow$ (iii). The implication (ii) $\Rightarrow$ (iii)
follows by duality.
\qed

\Corollary{9.3'}
{In a weak Hopf algebra we have $S|_{\A_{L/R}}=S_{L/R}$ as
given in \no{4.6}.
}
\proof
If $S$ is a pre-antipode then
$
S(\onne\1)\onne\2=\onne=\onne\1 S(\onne\1)
$
implying for $a\in\AL$
$$
S_L(a)\equiv\eRL(a)=S(a\1) a\2=S(\onne\1a)\onne\2=S(a)S(\onne\1)\onne\2=S(a)
$$
and similarly $S_R(b)=S(b)$ for $b\in\AR$.
\qed

\Corollary{9.3}
{A linear map $S:\A\to\A$ on a bimonoidal weak bialgebra $\A$
is an antipode if and only if $S$ is a normal pre-rigidity
map. In this case $S$ is always a bialgebra anti-automorphism
and $S^{-1}$ is a pode.
}
\proof
This follows from \prop{9.2}, \cor{8.3'}, \cor{8.7} and
\prop{9.7}.
\qed

\bsn
In specific examples the (co)monoidality axioms are typically
much harder to verify then anti-(co)multiplicativity of the
antipode.
In this light the following Theorem is very useful.

\Theorem{9.4}
{Let $\A$ be a weak bialgebra with pre-antipode $S$ and
assume $S$ to be a bialgebra anti-morphism. Then the following
statements are equivalent:
\\
i) $S$ is an antipode and $\A$ is a weak Hopf algebra.
\\
ii) $\A$ is bimonoidal
\\
iii) $\A$ is right-monoidal and right-comonoidal
\\
iv) $[\AL,\AR]=0$
\\
v) $[\hAL,\hAR]=0$
}
\proof
If $\A$ is bimonoidal, then $S$ is an antipode by \lem{8.3'},
thus proving (i) $\Lra$ (ii).
The implication (ii) $\Rightarrow$ (iii) being trivial
let us next prove (iii) $\Rightarrow$ (iv)+(v).
By \thm{3.5} part (iii) implies $\As=\Assp$ and $\hAs=\hAssp$
for $\s\neq\s'\in\{L,R\}$.
Thus, (iv) and (v) follow from 1i) and 2i) of \thm{8.4}.
Finally, given  $[\AL,\AR]=0$, part 1ii) of \thm{8.4} implies
right-monoidality, and by \lem{3.6} also left-monoidality.
The dual of \cor{3.3} then gives  $[\hAL,\hAR]=0$. Dualizing
this argument we also conclude that  $[\hAL,\hAR]=0$ implies
$[\AL,\AR]=0$ and comonoidality. This proves (iv) $\Lra$ (v) and
(iv)+(v) $\Rightarrow$ (ii).
\qed

\bsn
We close this Section by proving that the generalized Kac
algebras of T. Yamanouchi [Ya] are special kinds of weak Hopf
algebras.
Following [Ya] a generalized Kac algebra is a finite dimensional von
Neumann algebra $\A$ equipped with a coassociative non-unital $*$-algebra
map $\Del:\A\to\A\o\A$, a $*$-preserving involutive bialgebra
antiautomorphism $S:\A\to\A$ and a positive faithful
$S$-invariant trace $\l$ on $\A$ satisfying
\beq\lb{8.1'}
a\1\l(ba\2)=S(b\1)\l(b\2a),\quad\forall a,b\in\A.
\eeq
It follows [Ya] that $\A$ also admits a counit $\e$. Hence,
$\A$ is in fact a weak Hopf algebra, since we have more
generally

\Theorem{8.7}
{Let $\AAA$ be a weak bialgebra and $S:\A\to\A$ a bialgebra
antiautomorphism. Assume there exists a nondegenerate $\l\in\hA$
satisfying \no{8.1'}. Then $S$ is an antipode and $\A$ is a weak
Hopf algebra.
}
\proof
First we prove that $\A$ is monoidal. Let $l,r\in\A$ be the
unique solutions of
\beq\lb{8.2'}
l\arr\l=\e=\l\arl r.
\eeq
Then for all $a\in\A$
$$
a=a\1\e(a\2)=\left\{
\ba{rcl}
a\1\l(ra\2) &=&S(r\1)\l(r\2a)
\\
a\1\l(a\2l) &=& \l(al\2) S^{-1}(l\1)
\ea
\right.
$$
Hence, for $a,b\in\A$,
\beanon
a\1\e(ba\2) &=& a\1\l(rba\2)=S(b\1)S(r\1)\l(r\2b\2a)=S(b\1)b\2a
\\ &=&\one\1\e(b\one\2)a\equiv\eRL(b)a\,,
\eeanon
where the second line follows by putting $a=\one$ in the first
line. Thus, by (\ref{2.20}right), $\A$ is right monoidal.
Similarly,
\beanon
a\1\e(a\2b) &=&
a\1\l(a\2bl)=S^{-1}(l\1)S^{-1}(b\1)\l(ab\2l\2)=ab\2S^{-1}(b\1)
\\ &=&a\one\1\e(\one\2b)\equiv a\eRR(b)
\eeanon
and, by (\ref{2.19}left), $\A$ is right monoidal.
In particular, we also get $S(b\1)b\2=\eRL(b)$ and therefore
also, putting $b=S(a)$,
$$
b\1S(b\2)=S(S(a\1)a\2)=(S\circ\eRL\circ S^{-1})(b)=\eLR(b)\,.
$$
Thus, $S$ is a pre-antipode and the remaining claims follow
from \prop{9.2}ii).
\qed

\bsn
A deeper investigation of weak Hopf algebras including a
theory of integrals and $C^*$-structures will be given in [BNS].
In particular, there we will see that \no{8.1'} is one of the
defining relations of a left integral $\l\in\hA$ and that the
elements $l,r$ defined in \no{8.2'} are nondegenerate left and
right integrals, respectively, in $\A$ satisfying $l=S(r)$.
This generalizes well known results for ordinary finite
dimensional Hopf algebras by [LS].


\begin{appendix}

\renewcommand{\thesection}{\,}

\renewcommand{\theequation}{\mbox{\Alph{subsection}.\arabic{equation}}}
\renewcommand{\thesubsection}{\Alph{subsection}}
\renewcommand{\thetheorem}{\Alph{subsection}\arabic{theorem}}
\ncm{\subsec}{\setc{0}\setcounter{theorem}{0}\subsection}

\sec{Appendix}

\subsec{The B\"ohm-Szlach\'anyi Axioms}

In this Appendix we relate our \defi{9.1} of weak Hopf algebras
to the setting of G. B\"ohm and K. Szlach\'anyi.
In [BSz,Sz] the authors required $\A$ to be a weak bialgebra
satisfying the ``almost--monoidality" axioms {\sc l} and {\sc r} of
\lem{2.7}. Moreover, the antipode $S$ was
required to be a bialgebra anti-automorphism satisfying the
following two equivalent relations for all $a\in\A$
\bleq{9.4}
S(a\1)a\2\o a\3=(\one\o a)\Del(\one)
\qquad ,\qquad
a\1\o a\2S(a\3)=\Del(\one)(a\o\one)\,.
\eeq
Note that if $S$ is not required to be invertible, then the two
relations in \no{9.4} are independent of each other.
We now show that the BSz-axioms are equivalent to our
\defi{9.1}.

\Lemma{9.8}
{A linear map $S:\A\to\A$ on a weak bialgebra $\A$
satisfies \no{9.4} if and only if $\A$ is
right-comonoidal and $S$ is a pre-antipode.
In this case $\A$ is a weak Hopf algebra and $S$ is an antipode
(whence invertible by \cor{9.3}) if and only if $S$ is
a bialgebra anti-morphism and the counit axiom
\no{R} holds, i.e.
\bleq{9.5}
\e(ab) =\e(a\onne\2)\e (\onne\1 b),\ \forall a,b\in\A\,.
\eeq
}
\proof
Applying $(\id\o\e)$ and $(\e\o\id)$, respectively, to \no{9.4}
proves that $S$ is a pre-antipode.
Hence, \no{9.4} also implies
\bleq{9.6}
\eRL(a\1)\o a\2=(\one\o a)\Del(\one)\qquad ,\qquad
a\1\o\eLR(a\2)=\Del(\one)(a\o\one)
\eeq
and by the duals of (\ref{2.23}right) or (\ref{2.24}right) $\A$ is
right-monoidal.
Conversely, if $\A$ is right-monoidal
then \no{9.6} holds, implying \no{9.4} for any pre-antipode $S$.
Next, if $\A$ is right-monoidal then \no{9.5} is equivalent to
$\A$ being also right-comonoidal by \lem{2.8}iv). In this case
$S$ is an antipode and $\A$ is a weak Hopf algebra if and only
if $S$ is a bialgebra anti-morphism, see \prop{9.2} and
\thm{9.4}.
\qed

\bsn
Next, we remark that if $K=\CC$ and $\A$ is a weak
$*$-bialgebra (i.e. a $*$-algebra such that the coproduct is a
$*$-algebra homomorphism), then $\e(a^*)=\overline{\e(a)},\
a\in\A$, and therefore \no{9.5} is equivalent to
\bleq{9.7}
\e(ab) =\e(a\onne\1)\e (\onne\2 b),\ \forall a,b\in\A\,,
\eeq
which is actually the axiom postulated in [BSz].
One readily verifies that in a weak $*$-bialgebra also
our left- and right- (co)monoidality axioms are equivalent.
Based on this observation we now show that in the $*$-algebra
setting of [BSz] the axioms \no{9.5} and \no{9.7} are in fact
redundant, as well as the BSz-requirements $S(a^*)^*=S^{-1}(a)$
and $\Del\c S=(S\o S)\c\Del_{op}$.

\Lemma{9.9}
{Let $\A$ be a weak $*$-bialgebra and let $S:\A\to\A$ be an algebra
anti-morphism. Then $S$ satisfies \no{9.4} if and only if $S$ is
an antipode and $\A$ is a weak Hopf algebra. Moreover, in this
case $S$ is a bialgebra anti-automorphism
 and $S(a^*)^*=S^{-1}(a),\ \forall a\in\A$.
}
\proof
By \lem{9.8} \Eq{9.4} is equivalent to $S$ being a pre-antipode and $\A$
being right-comonoidal, whence also left-comonoidal. Thus the
first statement follows from \prop{9.2}.
In this case one readily checks that $\bar S(a):=S(a^*)^*$
provides a pode and therefore $\bar S=S^{-1}$ by \cor{9.3} and
the uniqueness of (anti)podes.
\qed


\subsec{More on Rigidity Structures}

In this Appendix, inspired by ideas of Drinfel`d [Dr],
we show that rigidity structures $\SSS$ are unique up to
equivalence and that for any rigid weak bialgebra
$\A$ there exists a twisted coproduct $\Del'$ on $\A':=S(\A)$ given by
$\Del'(S(a))=(S\o S)(\Del(a))$.
Under the conditions $S(\onne)=\onne$ and $\e\circ S=S$ this will
further imply $S$ to be invertible.
First we need the following two Lemmas

\Lemma{7.5}
{In a pre-rigid weak bialgebra $\AAAA$ the following identities
hold for all $a,b\in\A$.
\bea
\lb{7.18}
ab\1 \o S(b\2)\,\al\, b\3 &=& a\1 b\1\o S(a\2 b\2)\,\al\, a\3 b\3
\\
\lb{7.19}
S(b\1)\,\al\, b\2 \o ab\3
&=& S(a\1 b\1) \,\al\, a\2 b\2 \o a\3 b\3
\\
\lb{7.20}
a\1\,\be\, S(a\2)\o a\3 b
&=& a\1 b\1 \,\be\, S(a\2 b\2)\o a\3 b\3
\\
\lb{7.21}
a\1 b\o a\2 \,\be\, S(a\3) &=& a\1 b\1 \o a\2 b\2\,\be\, S(a\3 b\3)
\eea
}
\proof
The last two equations follow from the first two by
passing to $\A_{op}^{cop}$.
Also, it is enough to prove \no{7.18} and \no{7.19} for
$b=\onne$. Using \no{7.7f}
and (\ref{2.21}left) we compute
\beanon
a\onne\1\o S(\onne\2)\,\al\,\onne\3 &=& a\onne\1 \o\,\al\,\liS \onne\2\\
&=& a\onne\1\o\,\al\,\liS \eRL(\onne\2)\\
&=& a\1\o S(a\2)\,\al\, a\3
\eeanon
Using (\ref{2.21}right) and \no{7.7f}, \Eq{7.19} for $b=\onne$ (and
therefore for all b) follows similarly.
\qed

\Lemma {7.6}
{In a rigid weak bialgebra $\AAAA$ we have for all $a\in\A$
\bea
\lb{7.5.1}
a&=&a\1 \,\be\, S(a\2)\,\al\, a\3
\\
\lb{7.5.2}
S(a) &=& S(a\1)\,\al\, a\2 \,\be\, S(a\3)
\\
\lb{7.5.3}
\Del(a) &=& a\1\,\be\, S(a\4)\,\al\, a\5 \o a\2\,\be\, S(a\3)\,\al\, a\6
\\
\lb{7.5.4}
S(a\2)\o S(a\1) &=&
S(a\2)\,\al\, a\3\,\be\, S(a\6)\o S(a\1)\,\al\, a\4 \,\be\, S(a\5)
\eea
}
\proof
Eqs.\no{7.6} and \no{7.18} imply
$$
a=a\onne\1 \be S(\onne\2)\alpha\onne\3=a\1\be S(a\2)\alpha a\3\ .
$$
Similarly \no{7.7} and \no{7.21} give
$$
S(a)=S(\onne\1 a)\alpha \onne\2\be S(\onne\3) = S(a\1)\alpha a\2\be S(a\3)\ .
$$
Using \no{7.5.1} we now compute
\beanon
\Del(a) &=&
a\1\onne\1\beta S(\onne\2)\alpha\onne\3 \o a\2 \beta S(a\3)\alpha a\4
\\
&=& a\1\onne\1\beta S(\onne\4)\alpha\onne\5\o
a\2\onne\2\beta S(a\3\onne\3)\alpha a\4
\\
&=& a\1\onne\1\beta S(a\4\onne\4)\alpha a\5 \onne\5\o
a\2\onne\2 \beta S(a\3\onne\3)\alpha a\6
\\
&=& a\1 \beta S(a\4) \alpha a\5 \o a\2 \beta S(a\3) \alpha a\6
\eeanon
where in the second line we have used \no{7.21}
and in the third line \no{7.18}.
Similarly, we get
\beanon
 S(a\2)\o S(a\1)& =&
S(a\2)\alpha a\3\beta S(a\4)\o S(\onne\1 a\1)\alpha\onne\2\beta S(\onne\3)
\\
&=&S(\onne\2 a\2) \alpha\onne\3a\3\beta S(a\4)
\o S(\onne\1 a\1)\alpha\onne\4 \beta S(\onne\5)
\\
&=& S(\onne\2a\2)\alpha\onne\3a\3\beta S(a\6)
\o S(\onne\1a\1)\alpha\onne\4a\4\beta S(\onne\5 a\5)
\\
&=& S(a\2)\alpha a\3\beta S(a\6)\o S(a\1)\alpha a\4\beta S(a\5)
\eeanon
where in the first line we have used \no{7.5.2} and \no{7.7},
in the second line \no{7.18} and in the third line \no{7.21}.
\qed

\bsn
We now show that similar as in [Dr]
rigidity structures on a monoidal weak bialgebra are unique up
to equivalence.

\Prop{B1}
{Let $\A$ be monoidal and let
$(S_1, \alpha_1, \beta_1)$ and $(S_2,\alpha_2, \beta_2)$ be two
rigidity structures on $\A$. Define $u,\bar u\in \A$ by
\bea\lb{B1}
u &=& S_2 (\onne\1) \alpha_{2} \onne\2 \beta_1 S_1 (\onne\3)
\\\lb{B2}
\bar u &=& S_1 (\onne\1)\alpha_1 \onne\2\beta_2 S_2(\onne\3)
\eea
Then the following identities hold
for all $a\in\A$
$$
\begin{array}{rclcrcl}
uS_1(a) &=& S_2(a) u  &\qquad& \bar u S_2(a) &=& S_1(a)\bar u
\\
\al_2 &=& u\al_1 && \al_1 &=& \bar u\al_2
\\
\be_2 &=&\be_1\bar u && \be_1 &=&\be_2u
\\
u\bar u &=& S_2 (\onne) &\qquad& \bar u u &=&S_1(\onne)
\\
u\bar u u&=&u && \bar u u\bar u &=&\bar u
\end{array}
$$
Conversely, under these identities $u$ and $\bar u$ are
necessarily given by \no{B1} and \no{B2}.
}
\proof
By interchanging $1\leftrightarrow 2$ and $u\leftrightarrow\bar u$
it is enough to prove the left identities.
Using \no{7.19} and \no{7.21} we get for all $a\in\A$
$$
uS_1(a)=S_2(a\1)\al_2 a\2\be_1S_1(a\3) =S_2(a) u\,.
$$
From \no{7.5.1} and \no{7.19} one computes
$$
\al_2\equiv S_2(\onne\1)\al_2\onne\1 =
S_2(\onne\1)\al_2\onne\2\be_1 S_1(\onne\3)\al_1\onne_4
=uS_1(\onne\1)\al_1\onne\2\equiv u\al_1\,.
$$
Using \no{7.21} the identities $\be_2=\be_1\bar u$
follow similarly.
Moreover,
$$
u\bar u =
S_2 (\onne\1) \alpha_{2} \onne\2 \beta_1 S_1(\onne\3)
\al_1\onne\4\be_2 S_2(\onne\5)
=S_2 (\onne\1) \alpha_{2} \onne\2 \beta_2 S_2 (\onne\3)
=S_2(\onne)
$$
 Finally, $u\bar uu=S_2(\onne)u=u$.
Conversely, if $u,\bar u\in\A$ intertwine
$(S_1, \alpha_1, \beta_1)$ and $(S_2,\alpha_2, \beta_2)$ as
above, then $S_2(a) = uS_1(a)\bar u$, whence
$$
S_2(\onne\1)\al_2\onne\2\be_1S_1(\onne_3)
=uS_1(\onne\1)\al_1\onne\2\be_1S_1(\onne_3) =u\,,
$$
where we have used $\bar uu\al_1=S_1(\onne)\al_1=\al_1$.
\Eq{B2} follows similarly.
\qed

\bsn
This proof may of course be traced
back to the fact that in monoidal categories any
two rigidity structures are naturally equivalent.
In the same spirit, the following proposition reflects the natural
equivalence $\overline{V\x W}\cong \bar W\x\bar V$ in rigid
monoidal categories, see also [Dr].

\Prop{7.4}
{Let $\AAAA $ be a rigid weak bialgebra and let
$\FF,\bar\FF\in\A\o\A$ be given by
\bea
\lb{7.12}
\FF:= [S(\onne\2)\alpha \o S(\onne\1)\alpha] \Del(\onne\3\beta S(\onne\4))
\\
\lb{7.13}
\bar \FF := \Del (S(\onne\1)\alpha \onne\2)[\beta S(\onne\4)\o \beta S(\onne\3)]
\eea
Then the following identities hold for all $a\in\A$
\bea
\lb{7.14}
\FF\,\Del (S(a)) &=& (S\o S)(\Del_{op}(a))\,\FF\,
\\
\lb{7.15}
\Del(S(a))\,\bar \FF &=& \bar \FF\, (S\o S)(\Del_{op} (a))
\\
\lb{7.16}
\bar \FF\,\FF = \Del(S(\onne)) &,&
 \FF\,\bar \FF = (S\o S)(\Del_{op} (\onne))
\\
\lb{7.17}
\FF\,\bar \FF\,\FF=\FF &,& \bar \FF\, \FF\, \bar \FF =\bar \FF
\eea
}
\proof
To prove \no{7.14} we compute
\beanon
(S\o S)(\Del_{op} (a))\,\FF
&=&[S(\onne\2a\2)\alpha \o S(\onne\1a\1)\alpha] \Del(\onne\3\beta S(\onne\4))
\\
& =&[S(a\2)\alpha\o S(a\1)\alpha]\Del (a\3\beta S(a\4))
\\
& =&[S(a\2)\alpha a\3\o S(a\1)\alpha a\4] \Del(\beta S(a\5))
\\
& =&[S(\onne\2)\alpha\onne\3 \o S(a\1\onne\1)\alpha a\2] \Del(\beta S(a\3))
\\
& =&[S(\onne\2)\alpha\onne\3\o S(\onne_{(1')} \onne\1)\alpha\onne_{(2')}]
\Del(\beta S(a\onne_{(3')}))
\\
& =& [S(\onne\2)\alpha\onne\3 \o S(\onne\1)\alpha \onne\4]
\Del (\beta S(a\onne\5))
\\
&=& \FF\, \Del (S(a))
\eeanon
Here we have used \no{7.21}
in the second line, \no{7.19} in the fourth and the
fifth line and \no{7.18} in the sixth line.
Interchanging $\alpha$ and $\beta$
and repreating this proof in $\A_{op}^{cop}$ yields \no{7.15}.
To prove \no{7.16} we compute
\beanon
\bar \FF\, \FF
&=& \Del(S(\onne\1)\alpha\onne\2)
[\beta S(\onne_{(2')}\onne\4)\alpha\o\beta S(\onne_{(1')}\onne\3)\alpha]
\Del(\onne_{(3')}\beta S(\onne_{(4')}))
\\
&=& \Del (S(\onne\1)\alpha\onne\2)
[\beta S(\onne\4)\alpha\o\beta S(\onne\3)\alpha]
\Del(\onne\5\beta S(\onne\6))
\\
&=& \Del (S(\onne\1)\alpha)\Del (\onne\2)\Del (\beta S(\onne\3))
\\
&=& \Del (S(\onne))
\eeanon
where in the second line we have used
\no{7.21} and in the third line \no{7.5.3}.
Next, using \no{7.14} we get
\beanon
\FF\,\bar \FF
&=& [S(\onne\2)\o S(\onne\1)] \,\FF\,\Del(\alpha\onne\3)
[\beta S(\onne\3)\o\beta S(\onne\4)]
\\
&=& [S(\onne_{(2')}\onne\2)\alpha\o S(\onne_{(1')}\onne\1)\alpha]
\Del(\onne_{(3')}\beta S(\onne_{(4')} ) \alpha\onne\3)
[\beta S(\onne\5)\o \beta S(\onne\4)]
\\
&=& [S(\onne_{(2')} \onne\2)\alpha\o S(\onne_{(1')} \onne\1)\alpha]
\Del(\onne_{(3')}\onne\3 \beta S(\onne_{(4')}\onne\4)\alpha \onne\5)
[\beta S(\onne\7)\o\beta S(\onne\6)]
\\
&=& [S(\onne\2)\alpha\o S(\onne\1)\alpha]
\Del(\onne\3\be S(\onne\4)\al\onne\5)[\beta S(\onne\7)\o\beta S(\onne\6)]
\\
&=& [S(\onne\2)\alpha\o S(\onne\1)\alpha]
\Del(\onne\3)[\beta S(\onne\5)\o\beta S(\onne\4)]
\\
&=& S(\onne\2)\o S(\onne\1).
\eeanon
Here we have used \no{7.21} in the third line,
\no{7.5.1} in the fifth line and
\no{7.5.4} in the last line.
Finally, \no{7.17} follows from the obvious identities
\beq\lb{7.22}
\FF\,\Del (S(\onne))=\FF\quad ,\quad
\Del (S(\onne))\,\bar \FF=\bar \FF\ .
\qed
\eeq
If a rigidity structure satisfies $\e\c S=\e$, then
\prop{7.4} allows to define
a rigid monoidal weak bialgebra structure on $\A':= S(\A)$,
such that $S:\A\to\A'$ becomes a bialgebra homomorphism. First,
we consider $\A' \subset \A$ with opposite multiplication, i.e.
as a subalgebra of $\A_{op}$ with unit $\onne':= S(\onne)$. The
coproduct $\Del':\A'\to \A'\o \A'$ is given by
$$
\Del' (S(a)) := \FF_{op}\, \Del_{op} (S(a)) \bar\FF_{op}
\equiv (S\o S)(\Del(a))\ ,
$$
which is clearly a coassociative and multiplicative. Moreover,
if $\e \c S=\e$ then $\e':=\e|_{\A'}$ is a counit for $\Del'$ and
therefore  $(\A',\onne',\Del', \e')$ becomes
a monoidal weak bialgebra%
\footnote
{Presumably, if $\e$ is not $S$-invariant, there
still may exist a transformed counit $\e'$ on $\A'$ satisfying
$\e'\c S=\e$.
}.

\Lemma{7.7}
{Let $\AAAA$ be a rigid weak bialgebra satisfying $\e\c S=\e$.
Put $S':=S|_{\A'},\ \al':=S(\al)$ and $\be':=S(\be)$.
Then $(S',\al',\be')$ provides a
rigidity structure on $\A'$.
}
The proof of \lem{7.7} is straightforward and therefore omitted. Using this
result we are now able to prove

\Theorem{7.8}
{Let $\AAAA$ be a rigid weak bialgebra satisfying $\e\c S=\e$.
Then $S$ is bijective if and only if $S(\onne)=\onne$.
}
\proof
The identity \no{7.6} requires $S$ to be nonzero. Iterating
\lem{7.7} and using $\dim \A<\infty$ we conclude
$S^{n+1}(\A)=S^n(\A)$ for some $n\in \NN$.
We show that if $S(\onne)=\onne$,
then this implies $S^n(\A)=S^{n-1}(\A)$, whence
$S(\A)=\A$ by induction, thus proving bijectivity of $S$.
Replacing $\A$ by $S^{n-1}(\A)$ it is
enough to consider the case $n=1$.
Thus, assume $S^2(\A)=S(\A)$ and therefore
$\Ker S\cap S(\A)=0$. Let $S'=S|_{S(\A)}$ and define
$$P:= S'^{-1}\c S:\A\to S(\A)$$
Then $P$ is a multiplicative projection satisfying
$P\c S=S$ and therefore
$$
P(aS(b))=P(a)S(b),\quad\forall a, b\in\A
$$
By \lem{7.9} below, if $S(\onne)=\onne$
there exists $p\in\A$ such that $P(a)=pa$
for all $a\in\A $.
Hence $p=P(\onne)=\onne$, implying $\Ker S= \Ker P=0$ and $\A=S(\A)$.
\qed

\bsn
\Lemma {7.9} {Let $\AAAA$ be a rigid weak bialgebra.
Denote $L$ the left multiplication
of $\A$ on itself and consider $\A\equiv\A_{S(\A)}$ as a right
$S(\A)$-module. If $S(\onne)=\onne$, then
\beq
\EndSA (\A_{S(\A)}) = L(\A)
\eeq
}
\proof
By standard arguments for rigid monoidal categories
(see e.g. [Ka]), for any left $\A$-module $V$ we have an
anti-multiplicative isomorphism
\beq\lb{conj}
\EndA V\cong \EndA \bar V,
\eeq
given by
$\EndA V\ni T\mapsto \bar T\in\EndA\bar V$, where
$$
\bar T:= L_{\bar V}^{-1} \c(A_V\times \onne_{\bar V})
 \c(\onne_{\bar V} \times T\times\onne_{\bar V}) \c
(\onne_{\bar V} \times B_V)\c R_{\bar V}\ .
$$
The inverse of the assignment $T\mapsto\bar T$ is given by
$$
T=R_V^{-1} \c(\onne_V\times A_V)\c
(\onne_V\times \bar T\times \onne_V)\c(B_V\times \onne_V)\c L_V
$$
In our setting one
straightforwardly checks that $\bar T$ coincides with the
restriction of the transpose $T^t$ to
$\bar V\equiv \hat V\cdot S(\onne)$.
We now apply this to $V=\hA$ with canonical
left $\A$-action $\pi_V(a)\phi := a\arr \phi$. If
$S(\onne)=\onne$, the
conjugate left $\A$-module is given by
$\bar V\equiv \hat V=\A$, with left $\A$-action
$\pi_{\bar V}(a) b:= bS(a)$. Hence, the
isomorphism \no{conj} gives
$$
\EndSA (\A_{S(\A)}) \equiv \EndA ({_\A} \bar V) \cong
\EndA ({_\A} V)\equiv\EndA ({_\A}\hA)\ .
$$
On the other hand, under the transposition $T\mapsto T^t$ we clearly have
$$
\EndA({_\A} \hA)\cong \EndA(\A_\A) = L(\A)
$$
where $\A_\A$ denotes the regular right $\A$-module,
being the natural dual of the left $\A$-module ${_\A}\hA$.
\qed

\bsn
As already remarked, the condition $\e\c S=\e$
in \thm{7.8} may presumably be
dropped, if there is a twisted counit $\e'$ for $\Del'$ on
$\A'\equiv S(\A)$ satisfying $\e'\circ S=\e$.
Also recall from \prop{7.2.2}
that the condition $S(\onne)=\onne$ holds
if $\alpha$ or $\beta$ are invertible, and therefore in
particular for normal rigidity maps.



\subsec{Minimal and Cominimal Weak Bialgebras}

In this Appendix we introduce a special ``minimal" class of
comonoidal weak bialgebras $\A$ as well as their ``cominimal"
duals $\hA$.
As a motivation recall that if $\A$ is comonoidal, then
$[\AL,\AR]=0$ by \cor{3.3}.
Moreover, by Eqs.\ \no{3.3} - \no{3.6} and \thm{3.5} we have
\beq\lb{5.0}
\Delta(ab)=(a\o b)\Delta(\onne) = \Del(\onne)(a\o b),~~~a\in\AL,b\in\AR.
\eeq
Hence, if $\A$ is comonoidal then $\B:=\AL\AR\subset \A$
provides a weak comonoidal sub-bialgebra, since $\Delta(\onne)\in\B\o\B$
by \lem{4.0}.
Also, since $\e$ restricts to the counit on $\B$,
if $\A$ is bimonoidal then so is $\B$, and if $\A$ is a weak Hopf
algebra then $\B$ is a weak Hopf subalgebra.
This observation motivates the following

\Definition{4.1}
{A comonoidal weak bialgebra $\A$ is called {\em minimal},
if $\A=\AL\AR$. A monoidal weak bialgebra $\A$ is called
{\em cominimal}, if its dual $\hA$ is minimal. If
$\A=\A_1\A_2$ is an algebra generated by two commuting subalgebras $\A_1$
and $\A_2$, then we call a minimal weak bialgebra structure $(\Del,\e)$
on $\A$ {\em adapted} (to $\A_1$ and $\A_2$), if $\A_1=\AL$ and $\A_2=\AR$.
}
It will turn out that an adapted weak bialgebra structure is uniquely
determined by $\e$ or by $P:=\Del(\one)$. Since every
comonoidal weak bialgebra contains a minimal one,
these results will be very useful when constructing general
examples of weak comonoidal (or Hopf) bialgebras, see e.g.
Examples 1-3 in Appendix D.
The results of this section will also be needed when
constructing weak Hopf algebra structures on a large class of
quantum chains known from physical models, where $\A_{L/R}$
will be the left/right wedge algebras of these models
([N3], see also Appendix D).

\bsn
Let us start with
preparing some useful formalism.
Given two $K$-vector spaces $\A_{1/2}$ of equal finite
dimension,
a bilinear form $Q:\A_1\o\A_2\to K$ is called nondegenerate, if the map
$Q_L:a\mapsto Q(a\o \cdot)$
(equivalently $Q_R:b\mapsto Q(\cdot \o b))$ is an isomorphism
$A_1\to \hat\A_2$
(isomorphism $A_2\to\hat\A_1)$. This holds if and only if $Q$ has
a {\em form-inverse} $P\equiv\sum u_i\o v_i\in\A_2 \o\A_1$ satisfying
(throughout we drop again summation symbols)
\bea\lb{4.2}
Q(a\o u_i) v_i &=& a , \quad\forall a\in \A_1\nonumber \\
u_i Q(v_i\o b) &=& b, \quad\forall b\in\A_2
\eea
Clearly,   $P$ as a functional $\hat\A_2\o\hat\A_1 \to K$
 is also nondegenerate and its form-inverse is given by $Q$.
Form-inverses are of course
uniquely determined if they exist.
If $\A_{1/2}\subset\A$ are two commuting subalgebras and
$\A=\A_1\A_2$, then with any
$\phi \in\hat\A$ we associate the bilinear functional
$Q_\phi:\A_1\o\A_2\to K$ given by
\beq\lb{4.1}
Q_\phi (a\o b):= \bra\phi\mid ab\ket
\eeq
If $Q_\phi$ is nondegenerate we denote its
form-inverse by $P_\phi$. Note that $\phi$ need not be nondegenerate as a
functional on $\A$ in order for $Q_\phi$ to be nondegenerate.

\Proposition{4.2}
{Let $\A=\A_1\A_2$ be generated by two commuting subalgebras
$\A_1$ and $\A_2$.
\\
i) If there exists an adapted minimal weak bialgebra
structure $(\Del,\e)$ on $\A$, then
$Q_\e:\A_1\o\A_2\to K$ is nondegenerate and its form-inverse $P_\e \in
\A_2\o\A_1$ is idempotent. Moreover, $\Del$ is uniquely fixed
by $\e$ via
\beq\lb{4.4}
\Del(ab)=(a\o b)P_\e\equiv P_\e(a\o b),~~~ a\in \A_1,\ b\in\A_2.
\eeq
ii) If $\A\cong \A_1\o_{\A_1\cap\A_2}\A_2$ then the relation
\no{4.4} provides a one-to-one correspondence between adapted
minimal weak bialgebra structures $(\Del,\e)$ on $\A$ and
idempotents $P\in \A_2\o \A_1$ which are nondegenerate as
functionals $\hA_2\o\hA_1\to K$ and satisfy
\beq\lb{4.5}
(z\o\onne) P=(\onne \o z) P,~~\forall z\in \A_1\cap \A_2.
\eeq
}
\proof
By \Eq{2.7.4} of \cor{3.1} $Q_\e :\AL\o\AR\to K$ is nondegenerate for any
comonoidal weak bialgebra. Thus part (i) follows from \Eq{5.0} provided
$P_\e:=\Del(\onne)$ is the form-inverse of $Q_\e$. However this follows
from the definitions \no{4.1} and \no{4.2},
since the counit property of $\e$ gives
\beq\lb{4.5'}
\ba{rcl}
Q_\e(a\o\onne\1)\onne\2 &=& \e(a\1) a\2 =a\\
\onne\1 Q_\e (\onne\2\o b) &=& b\1 \e (b\2) =b
\ea
\eeq
where $a\in\AL=\NLL(\A)$ and $b\in\AR=\NRR(\A)$, see \no{3.3}
and \no{3.6}.
To prove part (ii) first note that given $(\Del,\e)$ we may put
$P:=\Del(\onne_\A)$,
which is idempotent and nondegenerate by part (i).
For $z\in\A_1\o\A_2$
\Eq{4.4} then gives $\Del(z)=(z\o\one)P=(\one\o z)P$.
Conversely, let $P\in\A_2\o\A_1$ be idempotent and nondegenerate with
form-inverse $Q:\A_1\o\A_2\to K$. Then \Eq{4.5} implies
$$Q(az\o b) = Q(a\o zb)$$
for all $a\in\A_1,b\in\A_2$ and $z\in\A_1\cap\A_2$.
Hence, if $\A\cong \A_1\o_{\A_1\cap\A_2}\A_2$
the functional $\e\in\hA$
$$\e(ab):=Q(a\o b),~~a\in \A_1,b\in\A_2$$
is well defined and we have $Q=Q_\e$.
Moreover, by \Eq{4.5} $\Del:\A\to\A\o\A$ given by \no{4.4}
is also a well defined algebra map.
Clearly, $\Del$ is coassociative and comonoidal,
since $(\onne_\A\o P)$ commutes with $(P\o\onne_\A)$,
and $\e$ is the counit for $\Del$,
since $Q$ is the form-inverse of $P$.
\qed

\bsn
More generally, \prop{4.2} shows that for any comonoidal
weak bialgebra $\A$ we may put $P:=\Del(\one)\in\AR\o\AL$ to
get a sequence of minimal weak bialgebras
\bleq{C8}
\AL\o\AR\longrightarrow\AL\o_{\AL\cap\AR}\AR\longrightarrow
\AL\AR\subset\A\,,
\eeq
where the arrows are the natural projections, being also
weak bialgebra homomorphisms.
To describe the dual cominimal weak bialgebras observe that
$Q_\e:\AL\o\AR\to K$ being nondegenerate with  form inverse
$P\equiv\Del(\one)\in\AR\o\AL$ we have the natural $K$-linear
isomorphism
$
\HomK(\AL\o\AR,K)\ni\Phi\mapsto T_\Phi\in\EndK\AR
$
with inverse
$
\EndK\AR\ni T\mapsto\Phi_T\in\HomK(\AL\o\AR,K)
$
given by
\bea
T_\Phi(b)&:=&\onne\1\Phi(\one\2\o b),\quad b\in\AR
\\
\Phi_T(a\o b)&:=&\e(a\,T(b)),\quad a\in\AL\,,\ b\in\AR\,.
\eea

\Prop{C3}
{Let $\A$ be comonoidal and $\pi_{\he}:\hA\to\EndK\AR$ the unit
representation. Then
\\
i)
Due to the isomorphism
$
\EndK\AR\cong\HomK(\AL\o\AR,K)
$
we obtain
 \bleq{C10}
\EndK\AR\supset\EndAA\AR\supset\pi_{\he}(\hA)\,,
\eeq
as an inclusion of cominimal weak bialgebras dual to \no{C8},
with coproduct $\delta:\EndK\AR\to\EndK\AR\o\EndK\AR$ given by
\bleq{C9}
\delta T(a\o b):=\one\1\o\one_{(1')}\,\e(\one\2\one_{(2')}\,T(ab)),
\quad a,b\in\AR\,,T\in\EndK\AR\,.
\eeq
ii) $\pi_{\he}:\hA\to\EndK\AR$ provides a bialgebra homomorphism and
$\Ker\pi_{\he}=(\AL\AR)^\perp$, i.e. $\Ker\pi_{\he}\subset\hA$ is the
annihilator of $\AL\AR\subset\A$.
}
\proof
Clearly, $T\in\EndAA\AR\Lra\Phi_T\in\HomK(\AL\o_{\AL\cap\AR}\AR,K)$.
To show $T\in\pi_{\he}(\hA)\Lra\Phi_T\in\HomK(\AL\AR,K)$ we compute for
$\psi\in\hA,\ b\in\AR$ and $a\in\AL=\NLL(\A)$
\bleq{C7}
\Phi_{\pi_{\he}(\psi)}(a\o b)=\e(a(\psi\arr b))=\e(\psi\arr(ab))=
\bra\psi\mid ab\ket
\eeq
thus proving (i). Part (ii) follows, since $\AL\AR\subset\A$ is a weak
subbialgebra and $(\AL\AR)^\perp=\Ker\pi_{\he}$ by \no{C7}. Hence
$
\pi_{\he}:\hA\to\pi_{\he}(\hA)\cong\hA/(\AL\AR)^\perp\cong\HomK(\AL\AR,K)
$
is a weak bialgebra epimorphism.
\qed

\Corollary{C3'}
{Let $\A$ be monoidal. Then
\\
i)
$\pi_\e(\A)=\EndhAA\hAR$ if and only if
$\hAL\hAR\cong\hAL\o_{\hAL\cap\hAR}\hAR$.
\\
ii)
$\A$ is cominimal if and only if $\pi_\e$ is faithful.
}
These results generalize
the weak Hopf algebra structure
on $\Mat\!_K(N)$ given by [BSz]. In fact, the dual of the
BSz-construction is obtained by putting in \prop{4.2}
$\A_1=\A_2=K^N$ (i.e. the abelian
algebra of diagonal $N\x N$-matrices), $\A=\A_1\o\A_2$ and
$Q_\e(a\o b)=\sum a_ib_i,\ a,b\in K^N$, yielding $P_\e=\sum
e_i\o e_i$, where $e_i\in K^N$ are the minimal orthogonal
projections, $e_ie_j=\delta_{ij}e_j$.

\bsn
Next, we look at  minimal
weak {\em Hopf} algebras and recall our definition of the maps
$S_\sigma:\As\to \A_{-\sigma}$
and $\bar S_\sigma :\As\to \A_{-\sigma}$ given in \no{4.6}.
By \thm{3.4}i) and \thm{3.5}, if $\A$ is comonoidal these maps
are algebra anti-isomorphisms satisfying $\bar S_{L/R}=S_{R/L}^{-1}$,
and if $\A$ is a weak Hopf algebra with
antipode $S$, then  by \cor{9.3'} $S_{L/R}=S|_{\A_{L/R}}$.
Using this, \prop{4.2} now generalizes to a complete
characterization of minimal weak Hopf algebras of the form
$\A\cong\A_1\o_{\A_1\cap \A_2}\A_2$.

\Theorem{4.5}
{Let $\A=\A_1\A_2\cong\A_1 \o_{\A_1\cap \A_2}\A_2$, where
$\A_1$ and $\A_2$ commute. Then the relation
\bea\lb{4.17}
\e (a_L b_R) &=& \om(a_L S_R(b_R)),~~ a_L\in \A_L \equiv \A_1,\
b_R\in\A_R \equiv \A_2
\\\lb{4.20}
S(a_Lb_R) &=& S_R(b_R)S_L(a_L)
\eea
provides a one-to-one correspondence between adapted
weak Hopf algebra structures $(\Del,\e,S)$ on $\A$ and pairs
$(\om,S_R)$, where $\om:\A_1\to K$ is a nondegenerate
functional satisfying $\Ind\om=\onne$,
$S_R:\A_2\to \A_1$ is an algebra anti-isomorphism
restricting to the identity on $\A_1\cap \A_2$
and where $S_L=S_R^{-1}\c\theta_\om,\ \theta_\om:\A_1\to\A_1$
being the modular automorphism of $\om$.
}
\proof
If $(\Del,\e)$ is adapted and bimonoidal, then by
\prop{C6}ii) $\ \om:=\e|_{\A_L}$ is nondegenerate with
$\Ind\om=\onne$ and by \lem{C7} \Eq{4.17} holds
with $S_R:=\eLR|_{\AR}$, implying also
$S_L\equiv\eRL|_{\AL}=S_R^{-1}\c\theta_\om$ by part
(iii) of \prop{C6}.
Moreover, in this case $S$ in \no{4.20} is well defined and
anti-multiplicative, since for comonoidal weak bialgebras $\eRL|_{\A_R}
= \id\!_{\AR}$ and $\eLR|_{\AL} =\id\!_{\AL}$  implying
\beq\lb{4.24}
S_L|_{\AL\cap\AR} = S_R|_{\AL\cap\AR} = \id_{\AL\cap\AR}
\eeq
Using $\Delta(\onne)\in\AR\o\AL$ and
$S_R(\onne\1)\onne\2\equiv\eLR(\one\1)\one\2=\onne$ by
\no{2.15} we then compute for $a_L\in\AL,\ a_R\in\AR$
and $a=a_La_R$
$$
S(a\1)a\2=S(\onne\1 a_L)\onne\2 a_R =
S_L(a_L) S_R(\onne\1)\onne\2 a_R=\eRL(a_L)a_R= \eRL(a_La_R),
$$
where in the last equation we have used $a_R=\eRL(a_R)$ and
part (2ii,left) of \prop{2.6'}. Similarly, using $\onne\1
S_L(\onne\2)=\onne $ and $a_L=\eLR(a_L)$ we get
$$
a\1S(a\2)=a_L\onne\1 S(a_R\onne\2) =
a_L\onne\1 S_L(\onne\2) S_R(a_R)=a_L\eLR (a_R)=\eLR (a_La_R)
$$
by part (2ii,right) of \prop{2.6'}.
Hence, $S$ is a pre-antipode and therefore an antipode  by \prop{9.2}i).

Converseley, we now reconstruct $(\Del,\e)$ from $(\om,S_R)$.
First, since $S_R$ restricts
to the identity on $\A_1\cap \A_2$, the functional $\e$ is well
defined on $\A$ by \Eq{4.17} and
$Q_\e :\A_1\o\A_2\to K$ is nondegenerate. By \prop{4.2} we have
to show that its
form-inverse $P_\e$ is idempotent. Clearly, if
$x_i \o y_i \in\A_1\o\A_1$ is the form-inverse of
$\om$, then $P_\e =S_R^{-1} (x_i)\o y_i$.
Hence,
$$
P_\e^2 = S_R^{-1} (x_ix_j)\o y_jy_i = S_R^{-1} (x_j)\o y_j x_i y_i =P_\e
$$
where we have used \no{4.12} and $\Ind \om=\onne$.
Thus, by \prop{4.2} we get a uniquely
determined adapted comonoidal weak bialgebra structure
$(\Del,\e)$ on $\A$.
Since $\A_1\cap\A_2\subset\C(\A_1)$ we have
$\theta_\om|_{\A_1\cap\A_2}=\id$ and therefore \no{4.20} provides a
well defined algebra anti-automorphism $S:\A\to\A$.
Moreover, for $a_L\in\A_1$ and $b_R\in\A_2$ we get
$$
\ba{rcccccccl}
\eLR(b_R) &\equiv&\e(\one\1 b_R)\one\2
&=&\om(S_R(b_R)S_R(\one\1))\one\2 &=& \om(S_R(b_R)x_i)y_i &=&S_R(b_R)
\\
\eRL(a_L) &\equiv& \one\1\e(a_L\one\2)
&=&S_R^{-1}(x_i)\om(y_i\theta_\om(a_L)) &=&
S_R^{-1}(\theta_\om(a_L)) &\equiv& S_L(a_L)\,.
\ea
$$
Hence, by the above arguments, $S$ is an antipode and $\A$ is a weak
Hopf algebra.
\qed



\subsec{Examples}

{\bf Example 1:}\\
This example provides a minimal comonoidal weak bialgebra  $\A =
\A_L\otimes \A_R$ which is not monoidal. We choose
$\A_L:=K^3$, i.e. the commutative algebra of diagonal $(3\x
3)$-matrices, and $\A_R$ the algebra of upper
triangular $(2\x 2)$ matrices
\beq\lb{5.1}
\A_R:= \left\{ \left( \matrix{x&y\cr 0& z\cr} \right) \mid x,y,z\in
K\right\} \subset \Mat_K(2).
\eeq
Let $e_i,i=1,2,3$, denote the pairwise orthogonal minimal projections
in $\A_L$, and let $b_i,i=1,2,3$ be the basis in
$\A_R$ given by
\beq \lb{5.2}
b_1= \left(\matrix{1 &0\cr 0 &0\cr}\right) , \ b_2 = \left(
\matrix{0&1\cr0&0\cr} \right) , \ b_3 = \left( \matrix{0&0\cr0&1\cr}
\right)
\eeq
Following the lines of \prop{4.2}(ii) we define the coproduct
$\Delta:\A\to\A\o\A$ for $a_L\in\A_L$ and $a_R\in\A_R$ by
$$\Delta (a_La_R) =(a_L\o a_R)P$$
where the idempotent $P\equiv \Delta(\onne_\A) \in \A_R \o \A_L$  is defined to be
\beq\lb{5.3}
P=b_1\o (e_1 + e_2) +b_2\o e_2 + b_3\o e_3
\eeq
Using the relations $e_ie_j =\delta_{ij} e_i$ and
\bea
\lb{5.4}
b_1^2 = b_1 , \ b_1 b_2 = b_2 b_3 =b_2 , \ b_3^2=b_3\\
\lb{5.5}
b_2^2=b_2b_1 =b_3b_2 =b_1b_3=b_3b_1=0
\eea
one immediately verifies $P^2 =P$. Also, as a functional
$(\widehat{\A_R}) \otimes (\widehat{\A_L}) \to K$, $P$ is
nondegenerate. Hence, the counit $\e:\A_L \o \A_R\to K$ is given as
the form-inverse of $P$, i.e.
\beq\lb{5.6}
\e = e^1\o (b^1-b^2)\o (e^2\o b^2) \o (e^3 \o b^3)
\eeq
where $e^i$ and $b^j$ denote the dual basises. According to part (ii)
of \prop{4.2} these data define a comonoidal weak
bialgebra structure on $\A\cong \A_L\A_R$, which by \cor{3.5}
cannot be monoidal, since $\A_L$ is commutative
and $\A_R$ is noncommutative.
This example does not admit an antipode, since for $x=e_1b_1$ one
easily computes $\eLR(x)=e_1+e_2$, whence $e_2\eLR(x)=e_2$,
wheras $x\1S(x\2)=e_1\one\1S(\one\2b_1)$ implying $e_2x\1S(x\2)=0$.
Nevertheless, the dual of this example admits a
(non-normalizable) rigidity structure,
see Example 2.
\qed

\bsn
{\bf Example 2}\\
In this example we construct  rigidity structures $\SSS$ on cominimal
weak bialgebras of the form $\EndK\AR$ (or $\EndAA\AR$), where $\A$ is
comonoidal, see \prop{C3}. This will in particular show that
the dual of Example 1 is rigid, although it does not admit an
antipode.

It is more convenient to perform the construction on the dual
$\B:=\AL\o\AR$ (or $\B:=\AL\o_{\AL\cap\AR}\AR$, respectively).
Thus, dualizing \no{7.7f}, \no{7.6} and \no{7.7} we seek for a
map $\SB:\B\to\B$ and functionals $\al,\be\in\hat\B$ satisfying
for all $x\in\B$
\bea\lb{a}
&\SB(x\1)\al(x\2)x\3\in\AR\qquad , \qquad
x\1\be(x\2)\SB(x\3)\in\AL&
\\\lb{b}
& \eB\left(x\1\be(x\2)\SB(x\3)\al(x\4)x\5\right) = \eB(x)
\\\lb{c}
& \eB\left(\SB(x\1)\al(x\2)x\3\be(x\4)\SB(x\5)\right) =
\eB(\SB(x))\,.
\eea
Moreover, the normalization conditions $\al\liS\eB=\al$ and
$\eB\reS\be=\be$ of \prop{7.2.2} become
\bleq{d}
\eB\left(\SB(x\1)\al(x\2)x\3\right)=\al(x)
\qquad , \qquad
\eB\left(x\1\be(x\2)\SB(x\3)\right)=\be(x)\,.
\eeq
Choose now $S_R:\AR\to\AL$ an arbitrary linear bijection and let
$S_L:\AL\to\AR$ be the transpose of $S_R^{-1}$ with respect to
the pairing $Q_\e$, i.e.
$$
Q_\e(a\o S_L(b)):=Q_\e(b\o S_R^{-1}(a)),\quad a,b\in\AL\,.
$$
Note that $P_\e\equiv\one\1\o\one\2$ being the form inverse of
$Q_\e$ implies
\bleq{e}
S_L(\one\2)\o S_R(\one\1)=\one\1\o\one\2\,.
\eeq
For $x=a_La_R,\ a_L\in\AL,\,a_R\in\AR$ define
\bea\lb{r1}
&\SB(x):=S_L(a_L)S_R(a_R)&
\\\lb{r2}
&\al(x):=Q_\e(\one\o S_L(a_L)a_R) \qquad ,  \qquad
\be(x):=Q_\e(a_LS_R(a_R)\o\one)\,. &
\eea
In case we want $(\SB,\al,\be)$ to be well defined on
$\AL\o_{\AL\cap\AR}\AR$ we also have to require $S_R$ (and therefore
$S_L$) to be $(\AL\cap\AR)$-linear.
Using $\Del_\B(x)=a_L\one\1\o\one\2 a_R$ and
the comonoidality property we now have
$$
x\1\o x\2\o x\3 =a_L\one\1\o\one\2\one_{(1')}\o\one_{(2')}a_R
$$
and similarly for higher coproducts.
Hence, the identities \no{a} - \no{d} are immediately verified,
provided we have
\bleq{f}
\SB(\one\1)\al(\one\2)\one\3 =\one
=\one\1\be(\one\2)\SB(\one\3)\,.
\eeq
To check \no{f} use \no{e} to compute
\beanon
&
\SB(\one_{(1')})\al(\one_{(2')}\one\1)\one\2
= Q_\e(\one\o S_L(\one_{(2')})\one\1)\,S_R(\one_{(1')})\one\2
=Q_\e(\one\o\one\1)\one\2 =\one
&
\\
&
\one\1\be(\one\2\one_{(1')})\SB(\one_{(2')}=
\one\1S_L(\one_{(2')})Q_\e(\one\2 S_R(\one_{(1')})\o\one)
=\one\1 Q_\e(\one\2\o\one) =\one\,.
&
\eeanon
This proves \no{f} and therefore the rigidity identities \no{a}
- \no{d}.
Finally, we also have
$$
\Del_\B(\SB(x)) =
\one\1S_R(a_R)\o S_L(a_L)\one\2
=\SB(\one\2a_R)\o \SB(a_L\one\1)
=(\SB\o\SB)(\Del_\B^{op}(x))
$$
and therefore $\SB$ is anti-comultiplicative.
Hence, $(\hat\SB,\al,\be)$ provides a rigidity structure on
$\hat \B$.

\smallskip\noindent
In Example 1 one may choose
$$
\ba{rclcrclcrcl}
S_R(b_1) &=&e_1+e_2 &\quad,\quad& S_R(b_2) &=&e_2 &\quad , \quad&
S_R(b_3) &=&e_3
\\
S_L(e_1) &=&b_1-b_2 &\quad,\quad& S_L(e_2) &=&b_2 &\quad , \quad&
S_L(e_3) &=&b_3
\ea
$$
to obtain $\be=\e\equiv\hat\one$ given by \no{5.6} and
$\al=e^1\o b^1 + e^3\o b^3$.
In particular, this rigidity structure is not normalizable.
\qed

\bsn
{\bf Example 3}\\
In this example we extend minimal weak Hopf algebras of the form
$\B=\A_L\o \A_R$ to weak Hopf algebras $\A= \G \reli \B \cong
\A_L\>cros\G\<cros \A_R$ by a {\em two-sided crossed product construction}
with a Hopf algebra $\G$.

Let $\B=\A_L\o \A_R$ be a minimal weak bialgebra, with counit $\e_B$
and coproduct $\Delta_\B$.
Let $\GGG$ be a finite dimensional Hopf algebra and assume a left Hopf
module $\G$-action $\re:\G \o \A_L\to \A_L$ and a
right Hopf module $\G$-action $\li :\A_R\o \G\to \A_R$. Following [HN1]
we define the { \em two-sided crossed product}
$\A:=\A_L\>cros \G\<cros \A_R$ to be the vector space $\A_L\o \G \o
\A_R$ with multiplication structure
\beq\lb{5.7}
(a_L\>cros g\<cros a_R)(b_L \>cros h\<cros b_R):= \Big( a_L(g\1 \re
b_L)\>cros g\2 h\1 \<cros (a_R\li h\2)b_R \Big)
\eeq
Equivalently, $\A$ may be identified with the {\em diagonal crossed
product} $\G\reli \B$ via
\beq\lb{5.12}
 \G\reli \B\ni (g\, (a_L\o a_R)) \mapsto ((g\1 \re a_L)
\>cros g\2 \<cros a_R) \in\A.
\eeq
 Here, the multiplication in $\G\reli\B$
is fixed by either of the equivalent relations [HN1]
\beanon
g\,(a_L\o a_R) &=& \Big( (g\1 \re a_L)\o (a_R\li S^{-1} (g\3)) \Big)\, g\2
\\
(a_L\o a_R)\,g &=& g\2\,\Big( (S^{-1}(g\3) \re a_L)\o (a_R\li (g\1) \Big)
\eeanon
where $g\in\G$ and $a_{L/R} \in \A_{L/R}$. Since as a linear space $\G
\reli \B=\G\o \B$, it comes equipped with the natural
tensor product coalgebra structure from $\G$ and $\B$. With respect
to the identification \no{5.12} this induces a coalgebra
structure $(\Delta_\A,\e_\A)$ on $\A$ given by
\beq\lb{5.10}
 \Delta_\A(a_L\>cros g\<cros a_R) := (a_L\>cros g\1 \<cros \onne\1)\o
((g\2\re \onne\2) \>cros g\3 \<cros a\3 )
\eeq
\beq\lb{5.11}
\e_\A(a_L\>cros g\<cros a_R) := \e_\B ((S^{-1} (g) \re a_L) \o a_R)
\eeq
where $\onne\1 \o\onne\2 \equiv \Delta_\B (\onne_\B) \in \A_R \o \A_L$.
Assume now
\beq\lb{5.9}
\e_\B \big( (g\re a_L) \o a_R\big) = \e_\B\big( (a_L \o (a_R \li g) \big)
\eeq
for all $g\in\G$ and $a_\sigma \in\As$. Equivalently, since
$\onne\1\o\onne\2\in\AR\o \AL$  is the form-inverse of $\e_\B:
\AL\o\AR\to K$, this means
\beq\lb{5.8}
 \onne\1\o(g\re\onne\2) = (\onne\1\li g) \o\onne\2,\quad \forall g\in\G.
\eeq
Given this condition one easily verifes that $(\A,\onne_\A,\Delta_\A,\e_\A)$
becomes a comonoidal weak bialgebra extension of $\B$.
Moreover, using \lem{2.8} iii) and iv) one checks that if $\B$ is
left- or right-monoidal then so is $\A$.

Next, note that \no{5.9} implies the restricted functional
$\omega:= \e_\B|_{\AL} $ to be $\G$-invariant. Hence, if
$\B$ is also a weak Hopf algebra and if $S_R:\AR\to\AL$ is the
restriction of the antipode to $\AR$ (see \thm{4.5}),
then \Eq{5.9} implies for all
$a_\sigma \in\As$ and $g\in\G$
\beq\lb{5.13}
\omega \big (a_L S_R(a_R\li g) \big) = \e (a_L\o a_R\li g)= \omega
\big((g\re a_L) S_R(a_R) \big)=
\omega \Big( a_L \big( S_\G (g)\re S_R(a_R) \big) \Big).
\eeq
By the nondegeneracy of $\omega$ we conclude
\beq\lb{5.14}
S_R(a_R\li g) = S_\G(g) \re S_R(a_R)~~ ,~~ \forall a_R\in\AR,\ g\in\G,
\eeq
where $S_\G$ is the antipode on $\G$. Hence, by \thm{4.5} any
 weak Hopf algebra  in the form of a two-sided
crossed product $\A=\AL\>cros \G\<cros\AR$ as above
may uniquely be constructed from a
left Hopf module $\G$-action on $\AL$, a nondegenerate $\G$-
invariant functional $\omega$ on $\AL$ satisfying $\Ind\omega=\onne$
and an anti-isomorphism $S_R:\AR \to \AL$.
The right $\G$-action on $\AR$ is then given by \no{5.14} and the
compatibility condition \no{5.9} follows from \no{4.17}.
Moreover, the antipode $S_B:\B\to\B$ constructed in \no{4.20} extends
to an antipode $S_\A:\A\to\A$ by putting
\beq\lb{5.15}
S_\A (a_L \>cros g\<cros a_R) := S_R(a_R) \>cros S_\G (g) \<cros S_L (a_L)
\eeq
Note that this is indeed an algebra anti-isomorphisms, since one checks,
similarly as in \no{5.14},
\beq\lb{5.16}
S_L(g\re a_L) = S_L(a_L) \li S_\G (g)
\eeq
for all $g\in\G$ and $a_L\in \A_L$.
\qed

\bsn
Interpreting $\A_L$ and $\A_R$ as left and right   ``wedge algebras",
this construction puts a weak Hopf algebra structure
on the Hopf algebraic quantum chains of [NSz]. More general
quantum chains can be treated by allowing $\G$ itself to be a
weak Hopf algebra. Moreover,
using the methods of [HN1], the above example also
generalizes to the case where
$\G$ is the dual of a quasi-Hopf algebra. If in this case also
$\A_L=\A_R=\G$, this will provide a general ``blowing up"
procedure in the spirit of [BSz] from quasi-coassociative Hopf
algebras $\G$ to weak coassociative Hopf algebras
$\G\>cros\hG\<cros\G$ in our sense, with
equivalent representation categories.
More details on this
will be discussed elsewhere.

\bsn
{\bf Example 4:}
Let $\GGG$ be a finite dimensional Hopf algebra and let $\H \subset
\G$ be an $\Ad$-invariant Hopf subalgebra, i.e.
\beq  \lb{5.26}
g\1 \H\, S_\G (g\2) \subset \H~~,~~\forall g\in \G
\eeq
Assume $\H$ semisimple and denote $p=S_\G(p)= p^2\in\H$ the unique
normalized two-sided integral.
Then the crossed product $\A:=\H \>cros_{Ad}\,\G$ becomes a
weak Hopf algebra with
\beq\lb{5.27}
\Delta_\A (h\o g) := (hS_\G(p\1) \o p\2 g\1)\o (p\3 \o g\2)
\eeq
\beq\lb{5.28}
\e_\A (h\o g) := \lambda (h) \e_\G (g)
\eeq
where $h\in\H , g\in\G$, and where $\lambda \in\hH$ is the left
integral dual to $p$, i.e. the unique solution of
\beq\lb{5.29}
\lambda \arr p=\onne_\H.
\eeq
Clearly, $\e_\A$ is a right counit for $\Delta_\A$ and using
\beq\lb{5.30} S_\G (p\1)\o p\2 = p\2 \o S_\G^{-1} (p\1)
\eeq
and the identity $\lambda_L^{-1} = S_\G^{-1} \circ p_R$ [LS], $\e_\A$ is
also a left counit. The coassociativity of $\Delta_\A$ follows
from
\beq\lb{5.31}
p\1 \o p\2 S_\G (p_{(1')} ) \o p_{(2')} = p\1 p_{(1')} \o p\2 \o p_{(2')} .
\eeq
To see that $\Delta_\A$ is an algebra map we compute for $h,k\in\H$
and $f,g\in\G$
\beanon
&& \Delta (h\o g)\Delta (k\o f) =\\
&&= \Big(h g\1 kS_\G (p_{1'}) S_\G (p\1 g\2) \o p\2 g\3 p_{(2')} f\1\Big)
\o \Big(p\3 g\4 p_{(3')} S(g\5) \o g\6 f\2\Big)\\
&&= \Big(hg\1 k S_\G (p\1 g\2) \o p\2 g\3 f\1\Big) \o \Big(p\3\o g\4 f\2\Big)\\
 &&=\Delta (gh\1 kS(g\2)\o g\3 f) .
\eeanon
Here we have used $p\in \C(\G)$, since $ g\1 p S(g\2) =
\e(g) p$,  for all $g\in\G$,
which follows since the l.h.s. is again a two-sided
integral in $\H$. To see that $\A$ is comonoidal let
$P:=\Delta_\A(\onne_\A)$. Then
\beanon
(\onne_\A \o P) (P\o \onne_\A ) &=& (P\o\onne_\A) (\onne_\A\o P))\\
&=& \Big(S_\G (p\1)\o p\2\Big)\o \Big(p\3 S_\G(p_{(1')}) \o p_{(2')}\Big)
\o \Big(p_{(3')}\o\onne_\G\Big)\\
&= &\Big(S_\G (p_{(1')}) S_\G (p\1) \o p\2 p_{(2')} \Big)\o
\Big(p\3\o p_{(3')}\Big)\o \Big(p_{(4')} \o \onne_\G\Big)\\
&=& (\Delta_\A \o id) (P),
\eeanon
where we have used \no{5.31}.
Finally, applying \lem{2.8}iii) and iv) and using \no{5.30} together
with the Fourier transformation identities [LS]
\bea
p\1 \lambda (h p\2) &=& S_\G (h)\\
\lambda (S_\G^{-1} (p\1) h) p\2 &=& h
\eea
for all $h\in\H$
\footnote{For a review of the theory of Fourier
transformations on finite dimensional Hopf algebras see also [N1].},
one also checks that $(\A,\onne_\A,\Delta_\A,\e_\A)$ is monoidal.
In this example we have
\bea
\AL &=& \{ (h\o \onne_\G ) \mid h\in \H \} \cong \H\\
\AR &=& \{ (S(h\1)\o h\2) \mid h\in\H\} \cong \H_{op}
\eea
as well as the identities
$$
\ba{rclcrcl}
\eLL (h\o g) &=& (S^{-1} (g\2) hg\1\o\onne_\G) &,& \eRR (h\o g)&=&\e_\G
(g)  \Big( h\2\o S^{-1} (h\1) \Big)\\
\eLR (h\o g) &=& \e_\G(g) (h\o\onne_\G) &,& \eRL (h\o g)
&=&S^2_\G\Big((h\li g)_{(2)}\Big) \o S_\G  \Big( (h\li g)_{(1)} \Big)
\ea
$$
where $h\in\H,g\in\G$ and $h\li g:= S_\G^{-1} (g\2) h g\1$.
Using these formulas the reader is invited to check that $S_\A:\A\to
\A$ given by
$$S_\A (h\o g) := S_\G^2 (h\2 \li g\2) \o S_\G (h\1 g\1)$$
provides an antipode and therefore $\A$ is a weak Hopf algebra.
\qed

\bsn
Putting $\G = \CC G$ and $\H=\CC H$ for some finite group $G$ with
normal subgroup $H$, the above example
\footnote{possibly deformed by a cocycle} appears as a weak Hopf
symmetry in any Jones triple
$$\M^G\subset \M \subset \M\>cros G$$
where $\M$ is a von-Neumann factor and $G$ is a group of automorphisms
of $\M$ with inner part given by $H$
[NSW].

\bsn
{\bf Acknowledgements}
I thank G. B\"ohm and K. Szlach\'anyi for many useful
discussions and for sharing their knowledge with me.
I am also grateful to L. Vainerman for bringing the papers [M,
V, Ya] to my attention.
\end{appendix}

\end{document}